\documentclass[11pt,twoside]{amsart}

\usepackage{amsmath}
\usepackage{amsfonts}
\usepackage{amssymb}
\usepackage[all]{xy}

\title[Dg algebras with enough idempotents and their derived categories]{Dg algebras with enough idempotents, their dg modules and their derived categories}

\author{Manuel Saor\'\i{}n}

\dedicatory{Dedicated to the memory of Serge Ovsienko}

\newtheorem{Lemma1}{{Lemma}}[section]
\newtheorem{Theo1}[Lemma1]{{Theorem}}
\newtheorem{Def1}[Lemma1]{{Definition}}
\newtheorem{Prop1}[Lemma1]{{Proposition}}
\newtheorem{Claim1}[Lemma1]{{Claim}}
\newtheorem{Rem1}[Lemma1]{{Remark}}
\newtheorem{Cor1}[Lemma1]{{Corollary}}
\newtheorem{Ex1}[Lemma1]{{Example}}
\newtheorem{Qu1}[Lemma1]{{Question}}
\newtheorem{Exer}[Lemma1]{Exercise}

\newenvironment{Lemma}{\begin{Lemma1}}{\end{Lemma1}}
\newenvironment{Def}{\begin{Def1}\em}{\end{Def1}}
\newenvironment{Prop}{\begin{Prop1}}{\end{Prop1}}
\newenvironment{Rem}{\begin{Rem1}\rm}{\end{Rem1}}
\newenvironment{Theorem}{\begin{Theo1}}{\end{Theo1}}

\newenvironment{Cor}{\begin{Cor1}}{\end{Cor1}}

\newenvironment{Example}{\begin{Ex1}\em}{\end{Ex1}}

\newcommand{\lra}{\longrightarrow}

\newcommand{\ra}{\rightarrow}
\newcommand{\sdp}{\times\kern-.2em\vrule height1.1ex depth-.05ex}
\newcommand{\epi}{\lra \kern-.8em\ra}

\newcommand{\ol}{\overline}

 \normalbaselines
\subjclass[2010]{Primary: 16E45, 18E30; Secondary: 16E35, 18E25}
\keywords{ Dg algebra, dg module, dg category,  dg functor, dg adjunction, homotopy category, derived category, derived functor.}

\date{}

\thanks{The author is highly indebted to Alexander Zimmermann for the careful reading of these notes, for his comments and for his help in improving the presentation. This work is backed by reseach projects from the Ministerio de Econom\'ia y Competitividad of Spain (MTM201346837-P) and the Fundaci\'on 'S\'eneca' of Murcia (19880/GERM/15), both with a part of FEDER funds. We thank these institutions for their support.}

\begin{document}

\maketitle

\begin{abstract}
We develop the theory dg algebras with enough idempotents and their dg modules and show their equivalence with that of small dg categories and their dg modules. We introduce the concept of dg adjunction and show that the classical covariant tensor-Hom and contravariant Hom-Hom adjunctions of modules over associative unital algebras are extended as dg adjunctions between categories of dg bimodules. The corresponding adjunctions of the associated triangulated functors are studied, and we investigate when they are one-sided parts of bifunctors which are triangulated on both variables. We finally show that, for a dg algebra with enough idempotents, the perfect left and right derived categories are dual to each other.
\end{abstract}

%
%
%

\section*{Introduction}

All throughout this paper, we fix a commutative ground ring   $K$ with unit and the term `category'
will mean `$K$-linear category',  unless otherwise specified, and all functors will be $K$-linear.

Small differential graded (dg) categories and their dg modules have played a fundamental role in  Mathematics for a long time. In the 70's and 80's of last century, they were the major tool to study matrix problems related to representation theory of algebras (see \cite{Kl-R}, \cite{GOR},  \cite{D1}, \cite{D2},...), which, among other things, led to Drozd's proof of the tame-wild dichotomy theorem (see \cite{D1} and \cite{D2}[Theorem 2]).  
In modern times, their main importance comes from a fundamental result of Keller (see \cite[Theorem 4.3]{K1}) which states that any compactly generated algebraic triangulated category is equivalent to the derived category of a small dg category. That importance grew even bigger when  Tabuada (see \cite{Ta}) showed that the category $Dgcat$ of small dg categories admits a model structure on which the weak equivalences are the so-called quasi-equivalences and To\"en (see \cite{T}) studied in depth the associated homotopy category $\text{Ho}(Dgcat)$, showing in particular that it had an internal $Hom$ and deriving several applications of this fact to Homotopy Theory and Algebraic Geometry.

By definition, a small dg category is a small category with a grading and a differential satisfying certain conditions (see the details in next section). But  from the time of Gabriel's thesis (see \cite{G}) one knows that small categories may be viewed as algebras with enough idempotents (or rings with several objects in the spirit of \cite{M}), and vice versa. Furthermore, if $A$ is such a small category then the  category $[A^{op},\text{Mod}-K]$  of contravariant functors is equivalent to $\text{Mod}-A$, the category of unitary right $A$-modules, when $A$ is viewed as an algebra with enough idempotents. It is natural to expect that the mentioned one-to-one correspondence extends to the dg setting. That requires the development of a theory of dg algebras with enough idempotents and their dg modules. This development is, in some sense, a demand of a part of the mathematical community. Indeed, apart from the unavoidable technicalities concerning the use of signs, the language of small dg categories and their dg modules is very technical and elusive for many people and, although the terminology is sometimes similar, concepts as dg modules or dg bimodules over small dg categories are intuitively far from the traditional concepts of module or bimodule over an associative algebra. This leads  some mathematicians to avoid the topic and others to present results about small dg categories and their derived categories only in terms of dg algebras (equivalently,  dg categories with just one object). This demand is the main motivation for these notes. They were initially thought as an appendix to a joint paper with Alexander Zimmermann (see \cite{SZ2}), where we needed to use some adjunctions between categories of dg bimodules, which we could not find explicit in the literature of small dg categories  and which became excessively unintuitive in that language (see \cite{NS-GTT}). As the notes grew longer than expected, we decided to offer them as a separated paper. Since a thorough development of the topic is out of question, we have concentrated on the basic aspects, with emphasis on those  needed for \cite{SZ2}, leaving aside other  important features of the theory.  

Our goal in the paper is to develop the basics of the theory of dg algebras with enough idempotents and their dg modules, to show its equivalence with the theory of small dg categories and their dg modules, and to revisit dg functors between categories of dg modules and their derived versions.  In particular, we construct explicitly the correspondents in the new setting of the covariant tensor-Hom and the contravariant Hom-Hom adjunctions of (bi)module categories over algebras with unit, together with their derived versions. Since the notes are specially aimed at making the  dg world more accessible to people that work with algebras and modules in the classical way, even at the cost of an excessive length, we have taken care in checking  essentially all the details in proofs. This care has been special on what concerns  signs in equations, which  are  most elusive for beginners and very important in the dg context, but whose associated calculations are rarely found explicit in the literature.

The organization of the paper goes as follows. In Section 1 we recall the definitions of dg category (not necessarily small) and  dg functor. In Section 2 we define what a dg algebra with enough idempotents is and give its category $Dg-A$ of right dg modules, proving in Section 3 that there is a one-to-one correspondence between small dg categories and dg algebras with enough idempotents, and showing a dg equivalence between $Dg-A$ and the category $\mathcal{C}_{dg}A$ of dg modules over the associated small dg category (see Theorem \ref{thm.dg modules versus dg functors}). In Sections 4 and 5 we define left dg modules and dg bimodules, respectively, and show that their corresponding categories can be realized as categories of right dg modules. In Section 6 we introduce the homotopy and derived category of a dg algebra with enough idempotents, and state the corresponding version of the mentioned Keller's theorem (see Corollary \ref{cor.Keller theorem}). In Section 7 we study derived functors of dg functors between categories of dg modules over algebras with enough idempotents and study when they appear as `one-sided part' of a bifunctor which is triangulated on both variables. In our approach, a fundamental role is played by the concept of dg adjunction (see Definition \ref{def.dg adjunction}). In section 8 we define the correspondents of the classical tensor and $Hom$ bifunctors. Concretely,  we show that if $A$, $B$ and $C$ are dg algebras with enough idempotents, then there are canonical dg functors $\overline{HOM}_A(?,?):(C-Dg-A)^{op}\otimes (B-Dg-A)\longrightarrow B-Dg-C$ and $?\otimes_B ?:(C-Dg-B)\otimes (B-Dg-A)\longrightarrow C-Dg-A$, where $\overline{HOM}_A(M,X):=B\text{HOM}_A(M,X)C$ is the `unitarization' of the non-unitary dg $B-C-$bimodule $\text{HOM}_A(M,X)$.  In Section 9, we show that if $X$ is a dg $B-A-$bimodule, then we have   dg adjunctions $(?\otimes_BX:C-Dg-B\longrightarrow C-Dg-A,\overline{\text{HOM}}_A(X,?):C-Dg-A\longrightarrow C-Dg-B)$ and $(\overline{HOM}_{B^{op}}(?,X)^o:B-Dg-C\longrightarrow (C-Dg-A)^{op},\overline{HOM}_{A}(?,X):(C-Dg-A)^{op}\longrightarrow B-Dg-C)$, both of which give rise to adjunctions between the corresponding derived functors (see Theorems \ref{thm.adjunction tensor-HOM} and \ref{thm.adjunction X-reflexive}). In the final Section 10 we use the last contravariant adjunction to prove that, for any dg algebra with enough idempotents $A$ and taking $X=A$, the adjunction $(\overline{HOM}_{A^{op}}(?,A),\overline{HOM}_A(?,A))$ gives rise to  quasi-inverse dualities $\text{per}(A^{op})\stackrel{\cong^o}{\longleftrightarrow}\text{per}(A)$ between the left and right perfect derived categories.

The paper tries to be as self-sufficient as possible, but some classical concepts are used without being explicitly introduced.  For the general theory of modules over algebras, the reader is referred to \cite{AF} and \cite{reptheobuch}, and specifically for modules over nonunital rings and algebras, we refer to  \cite{Wis}.
All right (resp. left) modules $M$ over an  algebra
$A$ will be assumed to be unitary. That is, we will assume that $MA=M$
(resp. $AM=M$). The corresponding category is denoted by $\text{Mod}-A$ (resp. $A-\text{Mod}$). When a non-unitary module eventually appears it will be explicitly mentioned.  On what concerns graded algebras (or rings) and graded modules, the reader is referred to \cite{NVO} for the basic concepts. Although this reference deals with graded unital rings, only a minimal adaptation is needed when passing to graded nonunital algebras. Finally, we freely use some terminology about triangulated categories. Basic references for this are \cite{Neeman} and \cite[Chapter 10+ss]{KS}, but, for a given triangulated category, we denote by $?[1]$ the shift or suspension functor, that was denoted by $\Sigma$ or $T$ in these references. Given a triangulated category $\mathcal{D}$ a subcategory $\mathcal{T}$ is a \emph{thick subcategory} when it is closed under extensions, shifts and direct summands. When $\mathcal{S}$ is a class of objects of $\mathcal{D}$, we shall denote by $\text{thick}_\mathcal{D}(\mathcal{S})$ the smallest thick subcategory of $\mathcal{D}$ containing $\mathcal{S}$. Recall  (see \cite[Definition 2.1.1]{Neeman}) that a functor $F:\mathcal{D}\longrightarrow\mathcal{D}'$ between triangulated categories is a \emph{triangulated functor} when  there is a natural isomorphism $\phi_F :F\circ (?[1])\stackrel{\cong}{\longrightarrow}(?[1])\circ F$ such that, for each triangle $X\stackrel{u}{\longrightarrow}Y\stackrel{v}{\longrightarrow}Z\stackrel{w}{\longrightarrow}X[1]$ in $\mathcal{D}$, the sequence $F(X)\stackrel{F(u)}{\longrightarrow}F(Y)\stackrel{F(v)}{\longrightarrow}F(Z)\stackrel{\phi_{F,X}\circ F(w)}{\longrightarrow}F(X)[1]$ is a triangle in $\mathcal{D}'$. If $F,G:\mathcal{D}\longrightarrow\mathcal{D}'$ are triangulated functors, then a \emph{natural transformation of triangulated functors $\tau :F\longrightarrow G$} is a natural transformation such that, for each triangle in $\mathcal{D}$ as above, one has $\phi_{G,X}\circ\tau_{X[1]}=\tau_X[1]\circ\phi_{F,X}$. It is well-known, and will be frequently used through these notes, that if $\tau :F\longrightarrow G$ is a natural transformation as above, then the class of objects $X\in\text{Ob}(\mathcal{D})$ such that $\tau_X$ is an isomorphism is a thick subcategory of $\mathcal{D}$.

\section{Dg categories and dg functors} \label{section.dg categories dg functors}

In this section we collect some basic notations, mainly taken from \cite{K1} and \cite{K2} (see also \cite{Yek} for an in-preparation textbook), which will be
used all throughout these notes.
Recall that a differential graded (dg) $K$-module is a graded $K$-module $V=\bigoplus_{n\in\mathbb{Z}}V^n$,
together with a graded $K$-linear map $d:V\longrightarrow V$ of degree $+1$ such that $d\circ d=0$. The category which will be
indistinctly denoted by $Dg-K$ or $\mathcal{C}_{dg}K$
has as objects the dg $K$-modules. Moreover each space of morphisms $\text{HOM}_K(V,W)$ has a structure
of dg $K$-module given by the following data:

\begin{enumerate}
 \item[i)] The grading is $\text{HOM}_K(V,W)=\bigoplus_{n\in\mathbb{Z}}\text{HOM}_K^n(V,W)$, where
 $\text{HOM}_K^n(V,W)$ consists of the graded $K$-linear maps $\alpha:V\longrightarrow W$ of degree
 $n$, i.e., such that $\alpha (V^k)\subseteq W^{k+n}$, for all $k\in\mathbb{Z}$.
\item[ii)] The differential $d:\text{HOM}_K(V,W)\longrightarrow\text{HOM}_K(V,W)$, which is a graded
$K$-linear map of degree $+1$ such that $d\circ d=0$, is defined by the rule $d(\alpha)=d_W\circ \alpha-(-1)^{|\alpha|}\alpha\circ d_V$,
where $|?|$ denotes the degree,  whenever $\alpha$ is a homogeneous element of $\text{HOM}_K(V,W)$.
\end{enumerate}

For any dg $K$-module $V$ and for any $n\in\mathbb{Z}$, one puts $d^n:=d_{|V^n}:V^n\longrightarrow V^{n+1}$, and defines  $Z^n(V):=\text{Ker}(d^n)$, $B^n(V):=\text{Im}(d^{n-1})$ and $H^n(V):=Z^n(V)/B^n(K)$, which are called respectively the  \emph{($K$)-module of $n$-cycles}, the \emph{module of $n$-boundaries} and the \emph{$n$-homology module} of $V$, respectively. We say that $V$ is \emph{acyclic} when $H^n(V)=0$, for all $n\in\mathbb{Z}$. 

Note  that if $V$ and  $W$ are dg $K$-modules, the tensor product $V\otimes W:=V\otimes_KW$ also becomes
an object of $Dg-K$, where the grading is given by $(V\otimes W)^n=\oplus_{i+j=n}V^i\otimes W^j$ and the
differential $d:V\otimes W\longrightarrow V\otimes W$ by the rule
$$d_{V\otimes W}(v\otimes w)=d(v)\otimes w+(-1)^{|v|}v\otimes d_W(w),$$
for all homogeneous elements $v\in V$ and $w\in W$. All throughout these notes, we use the unadorned symbol $\otimes$ to denote $\otimes_K$. Given a dg $K$-module $V$,
one has an associated dg $K$-module $V[1]$, where the grading is given by $V[1]^n=V^{n+1}$,
for each $n\in\mathbb{Z}$, and where $d_{V[1]}=-d_V[1]$. That is, $d_{V[1]}(v)=-d_V(v)$, for
each homogeneous element $v\in V$.

The category $Dg-K$ is the prototype of a  \emph{differential graded (=dg) category}. This is
any category $\mathcal{A}$ such that, for each pair $(A,B)$ of its objects, the $K$-module of
morphisms, denoted indistinctly by $\mathcal{A}(A,B)$ or $\text{Hom}_\mathcal{A}(A,B)$, has a
structure of differential graded $K$-module so that the composition map $\mathcal{A}(B,C)\otimes\mathcal{A}(A,B)\longrightarrow\mathcal{A}(A,C)$
($g\otimes f\rightsquigarrow g\circ f$) is a morphism of degree zero of the underlying graded $K$-modules which commutes with the differentials. This means that $d(g\circ f)=d(g)\circ f+(-1)^{|g|}g\circ d(f)$ whenever $g\in\mathcal{A}(B,C)$ and $f\in\mathcal{A}(A,B)$ are homogeneous morphisms.
If $\mathcal{A}$ is a dg category, then the \emph{opposite dg category} $\mathcal{A}^{op}$
has the same class of objects as $\mathcal{A}$ and the differential on morphisms $d:\mathcal{A}^{op}(A,B)=\mathcal{A}(B,A)\longrightarrow\mathcal{A}(B,A)=\mathcal{A}^{op}(A,B)$
is the same as in $\mathcal{A}$, but the composition of morphisms is given as
$\beta^o\circ\alpha^o=(-1)^{|\alpha| |\beta|}(\alpha\circ\beta)^o$, where we use the
superscript $^o$ to emphasize that a morphism is viewed as one in $\mathcal{A}^{op}$.

If $\mathcal{A}$ and $\mathcal{B}$ are dg categories, then the \emph{tensor product
dg category} $\mathcal{A}\otimes\mathcal{B}$ has $\text{Ob}(\mathcal{A})\times\text{Ob}(\mathcal{B})$
as its class of objects and, for all pairs $(A,B),(A',B')\in \text{Ob}(\mathcal{A})\times\text{Ob}(\mathcal{B})$,
we define
$\text{Hom}_{\mathcal{A}\otimes\mathcal{B}}[(A,B),(A',B')]=\mathcal{A}(A,A')\otimes\mathcal{B}(B,B')$,
with its canonical structure of dg $K$-module. The composition of homogeneous morphisms in
$\mathcal{A}\otimes\mathcal{B}$ is given by the rule
$$(\alpha_1\otimes\beta_1)\circ (\alpha_2\otimes\beta_2)=
(-1)^{|\alpha_2| |\beta_1|}(\alpha_1\circ\alpha_2)\otimes (\beta_1\circ\beta_2).$$

When $\mathcal{A}$ and $\mathcal{B}$ are dg categories,  a \emph{dg functor}
$F:\mathcal{A}\longrightarrow\mathcal{B}$ is a graded  functor (i.e.
$F(\mathcal{A}^n(A,A'))\subseteq\mathcal{B}^n(F(A),F(A'))$, for all $n\in\mathbb{Z}$ and
$A,A'\in\text{Ob}(\mathcal{A})$) such that $F(d_\mathcal{A}(\alpha ))=d_\mathcal{B}(F(\alpha ))$,
for each homogeneous morphism $\alpha$ in $\mathcal{A}$. We will frequently use the following
criterion for dg functors from a tensor product dg category.

\begin{Lemma} \label{lem.dg bifunctor}
Let $\mathcal{A}$, $\mathcal{B}$ and $\mathcal{C}$ be dg categories and let
$F:\mathcal{A}\otimes\mathcal{B}\longrightarrow\mathcal{C}$ be an assignment on
objects $(A,B)\rightsquigarrow F(A,B)$ and an assignment on homogeneous morphisms
$\alpha\otimes\beta\rightsquigarrow F(\alpha\otimes\beta)$ such that
$|F(\alpha \otimes\beta)|=|\alpha| +|\beta|$. The following assertions are equivalent:

\begin{enumerate}
\item The given assignments define a dg functor $F:\mathcal{A}\otimes\mathcal{B}\longrightarrow\mathcal{C}$.
\item The following conditions hold:
\begin{enumerate}
\item For any fixed object $A\in\mathcal{A}$, the assignments $B\rightsquigarrow F(A,B)$
on objects and $\beta\rightsquigarrow F(1_A\otimes\beta)$ on morphisms define a dg functor
$\mathcal{B}\longrightarrow\mathcal{C}$.
\item For any fixed object $B\in\mathcal{B}$, the assignments $A\rightsquigarrow F(A,B)$
on objects and $\alpha\rightsquigarrow F(\alpha\otimes 1_B)$ on morphisms define a
dg functor $\mathcal{A}\longrightarrow\mathcal{C}$.
\item For all homogeneous morphisms $\alpha :A\longrightarrow A'$ and $\beta :B\longrightarrow B'$,
in $\mathcal{A}$ and $\mathcal{B}$, respectively, there is the equality
$$(-1)^{|\alpha | |\beta|}F(1_{A'}\otimes\beta)\circ F(\alpha\otimes 1_B)=
F(\alpha\otimes\beta)=F(\alpha\otimes 1_{B'})\circ F(1_A\circ\beta).$$
\end{enumerate}
\end{enumerate}
\end{Lemma}

\begin{proof}
$1)\Longrightarrow 2)$ Since $F$ is a dg functor it commutes with the differentials,
so that $d_\mathcal{C}(F(\alpha\otimes\beta))=F(d_{\mathcal{A}\otimes\mathcal{B}}(\alpha\otimes\beta))$,
for all homogeneous morphisms $\alpha:A\longrightarrow A'$ in $\mathcal{A}$ and
$\beta :B\longrightarrow B'$ in $\mathcal{B}$.
That is, we have an equality
$$d_\mathcal{C}(F(\alpha\otimes\beta))=
F(d_\mathcal{A}(\alpha)\otimes\beta)+(-1)^{|\alpha|}
F(\alpha\otimes d_\mathcal{B}(\beta ))\;\;\;\;\;\;\;\;\;\;\;\; (*). $$
On the other hand, by the definition of composition of morphisms in
$\mathcal{A}\otimes\mathcal{B}$,
we have an equality
$$(\alpha\otimes 1_{B'})\circ (1_A\circ\beta)=
\alpha\otimes\beta =(-1)^{|\alpha| |\beta|}(1_{A'}\otimes\beta)\circ (\alpha\otimes 1_A).$$
Applying $F$ to all members of this equality and using the
functoriality of $F$, we get condition 2.c.

We next check condition 2.a, condition 2.b following by an analogous argument.
The fact that, for fixed $A\in\mathcal{A}$,  the assignments $B\rightsquigarrow F(A,B)$
and $\beta\rightsquigarrow F(1_A\otimes\beta)$ define a $K$-linear graded functor
$F_A:\mathcal{B}\longrightarrow\mathcal{C}$ follows directly from the functoriality
of $F$. (The corresponding construction fixing an object $B$ of $\mathcal B$ is denoted $F^B$).
We just need to check the dg condition of $F_A$. That is, we need to prove
that if $B,B'\in\mathcal{B}$ are any two objects, then the following square is commutative
$$\xymatrix{\mathcal{B}(B,B')\ar[r]^{d_\mathcal{B}}\ar[d]^{F_A}&
\mathcal{B}(B,B')\ar[d]^{F_A}\\
\mathcal{C}(F_A(B),F_A(B'))\ar[r]^{d_\mathcal{C}}&
\mathcal{C}(F_A(B),F_A(B'))\ar@{=}[r]&\mathcal{C}(F(A,B),F(A,B')).} $$
Indeed we have $(F_A\circ d_\mathcal{B})(\beta )=F(1_A\otimes d_\mathcal{B}(\beta ))$
while $$(d_\mathcal{C}\circ F_A)(\beta )=d_\mathcal{C}(F(1_A\otimes\beta ))=
F(d_\mathcal{A}(1_A)\otimes\beta )+(-1)^{|1_A|}
F(1_A\otimes d_\mathcal{B}(\beta ))=F(1_A\otimes d_\mathcal{B}(\beta )),$$
due to the equality $(*)$ above and the fact that $d_\mathcal{A}(1_A)=0$.

\bigskip

$2)\Longrightarrow 1)$ Let $\alpha_1:A_1\longrightarrow A_2$ and
$\alpha_2:A_2\longrightarrow A_3$ be homogeneous morphisms in $\mathcal{A}$
and let $\beta_1:B_1\longrightarrow B_2$ and $\beta_2:B_2\longrightarrow B_3$
be homogeneous morphisms in $\mathcal{B}$. Due to condition 2.c, we  have
\begin{eqnarray*}
F[(\alpha_2\otimes\beta_2)\circ \lefteqn{(\alpha_1\otimes\beta_1)]=}\\
&=&(-1)^{|\alpha_1| |\beta_2|}F((\alpha_2\alpha_1)\otimes (\beta_2\beta_1))\\
&=&(-1)^{|\alpha_1| |\beta_2|}F((\alpha_2\alpha_1)\otimes 1_{B_3})\circ
F(1_{A_1}\circ (\beta_2\beta_1)))\\
&=&(-1)^{|\alpha_1| |\beta_2|}F^{B_3}(\alpha_2\alpha_1)\circ F_{A_1}(\beta_2\beta_1)\\
&=&(-1)^{|\alpha_1| |\beta_2|}F^{B_3}(\alpha_2)\circ F^{B_3}(\alpha_1)\circ
F_{A_1}(\beta_2)\circ F_{A_1}(\beta_1)\\
&=&(-1)^{|\alpha_1| |\beta_2|}F(\alpha_2\otimes 1_{B_3})\circ
F(\alpha_1\otimes 1_{B_3})\circ F(1_{A_1}\otimes\beta_2)\circ  F(1_{A_1}\otimes\beta_1).
\end{eqnarray*}
and
$$F(\alpha_2\otimes\beta_2)\circ F(\alpha_1\otimes\beta_1)=
F(\alpha_2\otimes 1_{B_3})\circ F(1_{A_2}\otimes\beta_2)\circ
F(\alpha_1\otimes 1_{B_2})\circ F(1_{A_1}\otimes\beta_1).$$
We then get that
$$F[(\alpha_2\otimes\beta_2)\circ (\alpha_1\otimes\beta_1)]=
F(\alpha_2\otimes\beta_2)\circ F(\alpha_1\otimes\beta_1)$$
because, by condition 2.c, we have
$$ F(1_{A_2}\otimes\beta_2)\circ F(\alpha_1\otimes 1_{B_2})=
(-1)^{|\alpha_1| |\beta_2|} F(\alpha_1\otimes 1_{B_3})\circ F(1_{A_1}\otimes\beta_2).$$
Moreover, we have $F(1_A\otimes 1_B)=F_A(1_B)=1_{F_A(B)}=1_{F(A,B)}$, for all
$A\in\mathcal{A}$ and $B\in\mathcal{B}$,  due to the functoriality of
$F_A:\mathcal{B}\longrightarrow\mathcal{C}$. Therefore $F$ is a (clearly graded)
$K$-linear functor $\mathcal{A}\otimes\mathcal{B}\longrightarrow\mathcal{C}$.

It remains to check that $F$ is a dg functor, which amounts to prove the equality
(*) above for all $\alpha$ and $\beta$ as there. Indeed, using condition 2.c,  we have
\begin{eqnarray*}
d_\mathcal{C}(F(\alpha\otimes\beta))&=&
d_\mathcal{C}(F(\alpha\otimes 1_{B'})\circ F(1_A\otimes\beta))\\
&=&d_\mathcal{C}(F(\alpha\otimes 1_{B'}))\circ F(1_A\otimes\beta)+
(-1)^{|F(\alpha\otimes 1_{B'})|}F(\alpha\otimes 1_{B'})\circ d_\mathcal{C}(F(1_A\otimes\beta))\\
&=&(d_\mathcal{C}\circ F^{B'})(\alpha )\circ F(1_A\otimes\beta)+
(-1)^{| \alpha|}F(\alpha\otimes 1_{B'})\circ (d_\mathcal{C}\circ F_A)(\beta ).
\end{eqnarray*}
But the fact that $F_A$ and $F^{B'}$ are dg functors implies that
$d_\mathcal{C}\circ F^{B'}=F^{B'}\circ d_\mathcal{A}$ and
$d_\mathcal{C}\circ F_{A}=F_{A}\circ d_\mathcal{B}$. Then, using condition 2.c again,  we have
\begin{eqnarray*}
d_\mathcal{C}(F(\alpha\otimes\beta ))&=&
(F^{B'}\circ d_\mathcal{A})(\alpha )\circ F(1_A\otimes\beta )+
(-1)^{|\alpha |}F(\alpha\otimes 1_{B'})\circ (F_A\circ d_\mathcal{B})(\beta )\\
&=&F(d_\mathcal{A}(\alpha )\otimes 1_{B'})\circ F(1_A\otimes\beta)+(-1)^{|\alpha |}
F(\alpha\otimes 1_{B'})\circ F(1_A\otimes d_\mathcal{B}(\beta ))\\
&=&F(d_\mathcal{A}(\alpha )\otimes\beta)+(-1)^{|\alpha |}F(\alpha\otimes d_\mathcal{B}(\beta )),
\end{eqnarray*}
so that the equality $(*)$ holds.
\end{proof}

\begin{Ex1} \label{ex.the regular dg bifunctor}
If $\mathcal{A}$ is a dg category, then the following data give a dg functor $\mathcal{A}(?,?):\mathcal{A}^{op}\otimes\mathcal{A}\longrightarrow Dg-K$:

\begin{enumerate}
 \item An assignment on objects
 $(A,A')\rightsquigarrow\mathcal{A}(A,A')=\text{Hom}_\mathcal{A}(A,A')$.
\item If $\alpha :A\longrightarrow B$ and $\alpha ':A'\longrightarrow B'$
are homogeneous morphisms in $\mathcal{A}$, then
$\mathcal{A}(\alpha^o\otimes\alpha '):\mathcal{A}(B,A')\longrightarrow\mathcal{A}(A,B')$
takes $f\rightsquigarrow (-1)^{|\alpha |(|\alpha '|+|f|)}\alpha '\circ f\circ\alpha$,
for each homogeneous element $f\in\mathcal{A}(B,A')$.
\end{enumerate}
\end{Ex1}
\begin{proof}
For a fixed object $A$ in $\mathcal{A}$,
$\mathcal{A}(?,A)=A^{\wedge}:\mathcal{A}^{op}\longrightarrow\mathcal{C}_{dg}K=Dg-K$
acts on morphisms as $A^{\wedge}(\alpha )(f)=(-1)^{|\alpha| |f|}f\circ\alpha$
whenever $f$ and $\alpha$ are composable homogeneous morphisms in $\mathcal{A}$.
Then $\mathcal{A}(?,A)$ is what Keller calls the free right dg $\mathcal{A}$-module
associated to $A$ (see \cite[Section 1.1]{K1}), today more commonly known as the
\emph{representable right dg $\mathcal{A}$-module} associated to $A$, and it
is then a dg functor. Dually $A^{\vee}=\mathcal{A}(A,?):\mathcal{A}\longrightarrow{C}_{dg}K=Dg-K$
is the representable left dg $\mathcal{A}$-module, which acts on morphisms as
$A^{\vee}(\alpha )(f)=\alpha\circ f$, and is then a dg functor. So conditions
2.a and 2.b of the last lemma hold.

On the other hand, if $\alpha:A\longrightarrow B$, $\alpha':A'\longrightarrow A'$
and $f:B\longrightarrow A'$ are as in the statement, then one has
$$[\mathcal{A}(\alpha^o\otimes 1_{B'})\circ\mathcal{A}(1_B^o\otimes\alpha ')](f)=
\mathcal{A}(\alpha^o\otimes 1_{B'})(\alpha'\circ f)=(-1)^{|\alpha| (|\alpha '|+|f|)}\alpha '\circ f\circ\alpha$$
while
$$[\mathcal{A}(1_A^o\otimes\alpha ')\circ\mathcal{A}(\alpha^o\otimes 1_{A'})](f)=
(-1)^{|\alpha| |f|}\mathcal{A}(1_A^o\otimes\alpha ')(f\circ\alpha )=
(-1)^{|\alpha| |f|}\alpha '\circ f\circ\alpha.$$
Therefore condition 2.c in last lemma also holds.
\end{proof}

\begin{Example} \label{ex.tensor product of dg functors}
Let $F:\mathcal{A}\longrightarrow\mathcal{A}'$ and $G:\mathcal{B}\longrightarrow\mathcal{B}'$
be dg functor between dg categories. The following data define a dg functor
$F\otimes G:\mathcal{A}\otimes\mathcal{B}\longrightarrow\mathcal{A}'\otimes\mathcal{B}'$:

\begin{enumerate}
\item On objects one defines $(F\otimes G)(A,B)=(F(A),G(B))$.
\item If $\alpha :A_1\longrightarrow A_2$ and $\beta:B_1\longrightarrow B_2$
are morphisms in $\mathcal{A}$ and $\mathcal{B}$, respectively,  then
$$
\xymatrix{
(\mathcal{A}\otimes \mathcal{B})[(A_1,B_1),(A_2,B_2)]\ar@{=}[d]&&
(\mathcal{A}'\otimes\mathcal{B}')[(F\otimes G)(A_1,B_1),(F\otimes G)(A_2,B_2)] \ar@{=}[d]\\
\mathcal{A}(A_1,A_2)\otimes\mathcal{B}(B_1,B_2)\ar[rr]&&
\mathcal{A'}(F(A_1),F(A_2))\otimes\mathcal{B}'(G(B_1),G(B_2))
}$$
is the map given by $(F\otimes G)(\alpha\otimes\beta)=F(\alpha )\otimes G(\beta)$.
\end{enumerate}
\end{Example}

\begin{proof}
We do not need to use Lemma \ref{lem.dg bifunctor},  but  the definition of the
composition of morphisms in the tensor product dg category. Then
a direct proof is easy and left as an exercise.
\end{proof}

\medskip

With each dg category $\mathcal{A}$, one canonically associates its
\emph{$0$-cycle category} $Z^0\mathcal{A}$ and its \emph{$0$-homology category}
$\mathcal{H}^0\mathcal{A}$. Both of them have the same objects as $\mathcal{A}$,
and as morphisms one puts $\text{Hom}_{Z^0\mathcal{A}}(A,A')=Z^0(\mathcal{A}(A,A'))$
and $\text{Hom}_{H^0\mathcal{A}}(A,A')=H^0(\mathcal{A}(A,A'))$. In both cases,
the composition of morphisms is induced from that of $\mathcal{A}$.
If $F:\mathcal{A}\longrightarrow\mathcal{B}$ is any dg functor, the
fact that it induces a morphism
$\mathcal{A}(A,A')\longrightarrow\mathcal{B}(F(A),F(A'))$ of graded $K$-modules which commutes with the differentials
implies that it also induces  a morphism of $K$-modules
$$\text{Hom}_{Z^0\mathcal{A}}(A,A')=Z^0(\mathcal{A}(A,A'))\longrightarrow
Z^0(\mathcal{B}(F(A),F(A')))=\text{Hom}_{Z^0\mathcal{B}}(F(A),F(A'))$$
resp. $$\text{Hom}_{H^0\mathcal{A}}(A,A')=H^0(\mathcal{A}(A,A'))\longrightarrow
H^0(\mathcal{B}(F(A),F(A')))=\text{Hom}_{H^0\mathcal{B}}(F(A),F(A')),$$
for all objects $A,A'\in\text{Ob}(\mathcal{A})$. It immediately follows that
these are the assigments on morphisms of well-defined $K$-linear functors
$Z^0F:Z^0\mathcal{A}\longrightarrow Z^0\mathcal{B}$ and
$H^0F:H^0\mathcal{A}\longrightarrow H^0\mathcal{B}$.

The following concepts will be useful in the sequel.

\begin{Def} \label{def.homological natural transformation}
Let $F,G:\mathcal{A}\longrightarrow\mathcal{B}$ be dg functors between dg categories.
A natural transformation of dg functors $\tau :F\longrightarrow G$ is a natural
transformation of $K$-linear functors such that $\tau_A:F(A)\longrightarrow G(A)$
is a homogeneous morphism of zero degree in $\mathcal{B}$, for each object
$A\in\mathcal{A}$. We will say that $F$ is a \emph{homological natural
transformation of dg functors} when, in addition, $\tau_A\in Z^0\mathcal{B}(F(A),G(A))$,
for each $A\in\mathcal{A}$.

A  \emph{natural isomorphism of dg functors}  is a homological natural transformation
$\tau :F\longrightarrow G$ which is pointwise an isomorphism.
\end{Def}

\section{Dg algebras with enough idempotents and their categories of right dg modules}
\label{section.dg algebras enough idempotents}

With 'algebra' instead of 'ring', the following concept is well-known (see \cite[Chapter 10, Section 49]{Wis}). Note that we use the term 'distinguished family' instead of the term 'complete family' used in this reference.

\begin{Def}
An \emph{algebra with enough idempotents} is an algebra $A$ which admits a family of
nonzero orthogonal idempotents $(e_i)_{i\in I}$ such that $\bigoplus_{i\in I}e_iA=A=\bigoplus_{i\in I}Ae_i$.
This family $(e_i)_{i\in I}$ will be called a \emph{distinguished family of orthogonal idempotents}.
A \emph{graded algebra with enough idempotents} is an algebra with enough
idempotents together with a grading $A=\bigoplus_{n\in\mathbb{Z}}A^n$ on it such that $A$ admits
a distinguished family of orthogonal idempotents consisting of homogeneous elements of degree $0$.
Without further remark, on a graded algebra with enough idempotents we only consider
distinguished families consisting of degree zero homogeneous idempotents.
\end{Def}

Note that, for $A$ as above,  to say that a right (resp. left) $A$-module  is unitary  is equivalent to say that we have an internal decomposition
$M=\bigoplus_{i\in I}Me_i$ (resp. $M=\bigoplus_{i\in I}e_iM$) as $K$-module. Recall that all our modules will be unitary, unless explicitly said otherwise.

The crucial concept for us is the following:

\begin{Def} \label{def.dgalgebra with enough idempotents}
A \emph{differential graded (dg) algebra with enough idempotents} is a pair $(A,d)$, where $A$
is a graded algebra with enough idempotents and $d:A\longrightarrow A$ is a morphism  of degree
$+1$ of graded $K$-modules, called the \emph{differential},
satisfying the following conditions: i)  $d\circ d=0$;  ii) $d(e_i)=0$ for all $i\in I$; and iii) (\emph{Leibniz rule}) $d(ab)=d(a)b+(-1)^{|a|}ad(b)$,
for all homogeneous elements $a,b\in A$.
\end{Def}

Given a dg algebra with enough idempotents $A=(A,d)$, the usual opposite
algebra has a canonical structure of graded algebra. However, the differential
$d$ would not satisfy  Leibniz rule  when viewed as a map
$d^o:A^{op}\longrightarrow A^{op}$. This forces to redefine the concept of
\emph{opposite graded algebra with enough idempotents} $A^{op}$ as the one
having the same underlying graded $K$-module as $A$, but where the multiplication
of homogeneous elements is defined by $a^o\cdot b^o:=(-1)^{|a| |b|}(ba)^o$, for
all $a,b\in A$. Here we use the upper index $^o$ to indicate that we are viewing
the element as one of the opposite graded algebra. The following is now routine:

\begin{Exer}
If $(A,d)$ is a dg algebra with enough idempotents and $A^{op}$ is the opposite graded algebra in
the above sense, then $d^o:A^{op}\longrightarrow A^{op}$ is a differential making the pair
$(A^{op},d^o)$ to be a dg algebra with enough idempotents (with the same distinguished
family of homogeneous idempotents as $A$).
\end{Exer}

The following gives the definition of the \emph{tensor product of two dg algebras
with enough idempotents}.

\begin{Lemma} \label{lem.tensor product dg algebras}
Let $A =(A,d)$ and $B=(B,d)$ be two dg algebras with enough idempotents and let $A\otimes B$ their tensor
product in $Dg-K$.
When one defines the multiplication of homogeneous tensors by the rule that
$(a\otimes b)\cdot (c\otimes d)=(-1)^{|b| |c|}ac\otimes bd$, $A\otimes B$
becomes a dg algebra with enough idempotents, with the same differential as in $Dg-K$.
\end{Lemma}

\begin{proof}
It is routine to check the associativity, so that $A\otimes B$ becomes an
associative graded algebra with the given multiplication. Moreover, if
$(e_i)_{i\in I}$ and $(e'_j)_{j\in J}$ are distinguished families of homogeneous
idempotents of degree 0 in $A$ and $B$, respectively, then
$(e_i\otimes e'_j)_{(i,j)\in I\times J}$ is a distinguished family
of homogeneous orthogonal idempotents of degree $0$ in $A\otimes B$.
On the other hand, the differential $d:A\otimes B\longrightarrow A\otimes B$
vanishes on each $e_i\otimes e'_j$
It remains to check Leibniz rule. It is also routine, but for the convenience
of the reader we explicitly give the calculations:
\begin{eqnarray*}
d[(a_1\otimes b_1)\cdot (a_2\otimes b_2)]&=&(-1)^{|b_1| |a_2|}d[(a_1a_2)\otimes (b_1b_2)]\\
&=&(-1)^{|b_1| |a_2|}[d(a_1a_2)\otimes (b_1b_2)+(-1)^{|a_1|+|a_2|}(a_1a_2)\otimes d(b_1b_2)]\\
&=&(-1)^{|b_1| |a_2|}[(d(a_1)a_2+(-1)^{|a_1|}a_1d(a_2))\otimes (b_1b_2)]\\
&&+(-1)^{|b_1| |a_2|+|a_1|+|a_2|}[(a_1a_2)\otimes (d(b_1)b_2+(-1)^{|b_1|}b_1d(b_2))]\\
&=&(-1)^{|b_1| |a_2|}d(a_1)a_2\otimes b_1b_2\\
&&+ (-1)^{|b_1| |a_2|+|a_1|}a_1d(a_2)\otimes b_1b_2\\
&&+(-1)^{|b_1| |a_2|+|a_1|+|a_2|}a_1a_2\otimes d(b_1)b_2\\
&&+(-1)^{|b_1| |a_2|+|a_1|+|a_2|+|b_1|}a_1a_2\otimes b_1d(b_2).
\end{eqnarray*}
while, on the other side, we have:
\begin{eqnarray*}
d(a_1\otimes b_1)\lefteqn{\cdot (a_2\otimes b_2)+(-1)^{|a_1|+|b_1|}(a_1\otimes b_1)\cdot d(a_2\otimes b_2)=}\\
&=&[d(a_1)\otimes b_1+(-1)^{|a_1|}a_1\otimes d(b_1)]\cdot (a_2\otimes b_2)\\
&&+(-1)^{|a_1|+|b_1|}(a_1\otimes b_1)\cdot [d(a_2)\otimes b_2+(-1)^{|a_2|}a_2\otimes d(b_2)]\\
&=&(-1)^{|b_1| |a_2|}d(a_1)a_2\otimes b_1b_2\\
&&+(-1)^{|a_1|}(-1)^{(|b_1|+1) |a_2|}a_1a_2\otimes d(b_1)b_2\\
&&+(-1)^{|a_1|+|b_1|}(-1)^{|b_1|(|a_2|+1)}a_1d(a_2)\otimes b_1b_2\\
&&+(-1)^{|a_1|+|b_1|}(-1)^{|a_2|}(-1)^{|b_1| |a_2|}a_1a_2\otimes b_1d(b_2).
\end{eqnarray*}
Therefore Leibniz rule holds for the given multiplication in $A\otimes B$.
\end{proof}

\medskip

Associated with any graded algebra with enough idempotents $A$, we have the
category $Gr-A$ of graded right $A$-modules, where the morphisms between two
objects $M$ and $N$ are the homomorphisms of right $A$-modules $f:M\longrightarrow N$
such $f(M^n)\subseteq N^n$, for all $n\in\mathbb{Z}$. The category $Gr-A$
comes with a \emph{shift functor} $?[1]:Gr-A\longrightarrow Gr-A$.
For each graded right $A$-module $M$, $M[1]$ has the same underlying
(ungraded) $A$-module as $M$, but the grading on $M[1]$ is given by
$M[1]^n=M^{n+1}$, for all $n\in\mathbb{Z}$. The action of $?[1]$ on
morphisms is the identity. It is clear that $?[1]$ is an equivalence
of categories, which allows to define the iterated powers $?[n]=(?[1])^n$,
for all $n\in\mathbb{Z}$. We then form the graded category $GR-A$.
Its objects are the same as in $Gr-A$ and, given two graded right $A$-modules
$M$ and $N$, we  define $$\text{HOM}_A(M,N):=\bigoplus_{n\in\mathbb{Z}}\text{Hom}_{Gr-A}(M,N[n])$$
as space of morphisms in $GR-A$. On this space of morphisms we have an obvious
grading given by $\text{HOM}_A^n(M,N):=\text{Hom}_{Gr-A}(M,N[n]),$ for each $n\in\mathbb{Z}$.
The composition $g\circ f$ in $GR-A$ of two homogeneous morphisms
$f:M\longrightarrow N[n]$ and $g:N\longrightarrow P[p]$ is defined as the
composition $g[n]\circ f$ in $Gr-A$.

\begin{Def} \label{def.dg module}
Let $A=(A,d)$ be a dg algebra with enough idempotents. A
\emph{right (resp. left) differential graded (dg) $A$-module}
is a pair $(M,d_M)$ consisting of a graded right (resp. left)
$A$-module $M=\bigoplus_{n\in\mathbb{Z}}M^n$ together with a morphism
$d_M:M\longrightarrow M[1]$ in $Gr-K$ such $d_M\circ d_M=0$ and $d_M(xa)=d_M(x)a+(-1)^{|x|}xd(a)$
(resp. $d_M(ax)=d(a)x+(-1)^{|a|}ad_M(x)$), for all homogeneous
elements $x\in M$ and $a\in A$.
\end{Def}

Suppose that $A$ is a dg algebra with enough idempotents and that $M$
is a right dg $A$-module. The graded right $A$-module $M[1]$ with its differential  $d_{M[1]}=-d_M$ as dg $K$-module (see Section \ref{section.dg categories dg functors}) becomes a right dg $A$-module.
Indeed, if one has $x\in M[1]^n=M^{n+1}$ and $a\in A^p$, then
\begin{eqnarray*}
d_{M[1]}(xa)&=&-d_M(xa)=-[d(x)a+(-1)^{n+1}xd(a)]=-d(x)a+(-1)^nxd(a)\\
&=&d_{M[1]}(x)a+(-1)^{|x|}xd(a),
\end{eqnarray*}
where $|x|=n$ is the degree of $x$ in $M[1]$.  In this way, we get a functor $?[1]:Dg-A\longrightarrow Dg-A$ which is 'almost' a dg functor, in the sense that if $d:\text{HOM}_A(M,N)\longrightarrow\text{HOM}_A(M,N)$ and $\delta :\text{HOM}_A(M[1],N[1])\longrightarrow\text{HOM}_A(M[1],N[1])$ are the respective differentials on $Hom$ spaces, then $\delta (f[1])=-d(f)[1]$, for each homogeneous morphism $f\in\text{HOM}_A(M,N)$. The reader is referred to section \ref{section.right versus left} to see that the corresponding functor for left dg modules produces suprising effects.

\begin{Prop} \label{prop.rightdgmodules as dg category}
Let $A$ be a dg algebra with enough idempotents. The following data give a dg category  $Dg-A$,
the dg category of right dg $A$-modules:
\begin{enumerate}
\item[-] The objects of $Dg-A$ are the right dg $A$-modules (see Definition
\ref{def.dg module});
\item[-] The morphisms in $Dg-A$ and the composition of them is defined as
in the category $GR-A$.
\item[-] For each pair $(M,N)$ of objects, the differential
$d:\text{HOM}_A(M,N)\longrightarrow\text{HOM}_A(M,N)$ on  $Hom$
spaces is defined by the rule $d(f)=d_N\circ f-(-1)^{|f|}f\circ d_M$,
for each homogeneous morphism $f$.
\end{enumerate}
\end{Prop}

\begin{proof}
We first need to check that the differential on $Hom$ spaces is well-defined,
i.e. that $d(f)$ is a homomorphism of graded right $A$-modules, which is
homogeneous of degree $|f|+1$,  whenever $f\in\text{HOM}_A(M,N)$ is
homogeneous. Indeed if $x\in M$ and $a\in A$  are homogeneous elements, then we have:
\begin{eqnarray*}
d(f)(xa)&=&[d_N\circ f-(-1)^{|f|}f\circ d_M](xa)\\
&=&d_N(f(x)a)-(-1)^{|f|}f(d_M(xa))\\
&=&d_N(f(x))a+(-1)^{|f|+|x|}f(x)d(a)-(-1)^{|f|}f(d_M(x)a+(-1)^{|x|}xd(a))\\
&=&(d_N\circ f)(x))a+(-1)^{|f|+|x|}f(x)d(a)-(-1)^{|f|}(f\circ d_M)(x)a\\
&&-(-1)^{|f|+|x|}f(x)d(a)\\
&=&(d_N\circ f)(x))a-(-1)^{|f|}(f\circ d_M)(x)a\\
&=&d(f)(x)a.
\end{eqnarray*}
Then $d(f)$ is a homogeneous morphism in $GR-A$, clearly of degree $|f|+1$.

Given the fact that the differential on $\text{HOM}_A(M,N)$ is the restriction of the differential on $\text{HOM}_K(M,N)$ and that the composition of morphism in $Dg-A$ is defined as in $Dg-K$, and the latter is a dg category,  the equality  $$d(g\circ f)=d(g)\circ f+(-1)^{|g|}g\circ d(f),\;\;\;\;\;\;\;\;\;\;(*)$$
holds for all homogeneous morphisms $f\in\text{HOM}_A(M,N)$ and $g\in\text{HOM}_A(N,P)$.  Then $Dg-A$ is also a dg category.

\end{proof}

\section{Dg algebras with enough idempotents versus small dg categories}

\label{dgalgebrascomingfromdgcategoriessection}

Let $A=(A,d)$ be a dg algebra with enough idempotents on which we fix a
distinguished family of orthogonal idempotents $(e_i)_{i\in I}$, which are
homogeneous of degree zero and such that $d(e_i)=0$, for all $i\in I$.
We can view $A$ as a small dg category as follows:

\begin{enumerate}
\item[.-] The set of objects is $\text{Ob}(A)=I$;
\item[.-] If $i,j\in A$, the set of morphisms of degree $n$ from $i$ to $j$ is
$A^n(i,j):=e_jA^ne_i$, for all $n\in\mathbb{Z}$;
\item[.-] The composition map
$A(j,k)\times A(i,j)=e_kAe_j\times e_jAe_i\longrightarrow e_kAe_i=A(i,k)$
is  the multiplication map.
\item[.-] The differential $d:A(i,j)=e_jAe_i\longrightarrow e_jAe_i=A(i,j)$
is  the differential of $A$ as a dg algebra, for all $i,j\in I$.
\end{enumerate}

It is routine to check that the data above make $A$ into a small dg category.
Conversely, let $\mathcal{A}$ be a small dg category. We can view $\mathcal{A}$
as a dg algebra with enough idempotents as follows:

\begin{enumerate}
\item[.-] The elements of $\mathcal{A}$ are those of
$\bigoplus_{A,B\in\text{Ob}(\mathcal{A})}\mathcal{A}(A,B)$, and we put
$$\mathcal{A}^n=\bigoplus_{A,B\in\text{Ob}(\mathcal{A})}\mathcal{A}^n(A,B)$$
for the $K$-module of elements of degree $n$ in $\mathcal{A}$, for all $n\in\mathbb{Z}$.
\item[.-] The multiplication  in $\mathcal{A}$ extends by
$K$-linearity the composition of morphisms in  $\mathcal{A}$.
\item[.-] The differential $d:\mathcal{A}\longrightarrow\mathcal{A}$
is the direct sum of the differentials
$d_{A,B}:\mathcal{A}(A,B)\longrightarrow\mathcal{A}(A,B)$, as $A,B$
vary on the set of objects of $\mathcal{A}$.
\end{enumerate}
It is routine to see that the data above make $\mathcal{A}$ into
a dg algebra with enough idempotents, where the identities
$e_A:=1_A$ ($A\in\text{Ob}(\mathcal{A})$) form a distinguished
family of orthogonal idempotents of degree zero. Note that we have
$\mathcal{A}(A,B)=e_B\mathcal{A}e_A$.

The processes explained above of passing from dg algebras with enough
idempotents to small dg categories, and viceversa, are clearly inverse
to each other. This allows us to pass freely from one language to the
other. Note, in particular, that the opposite dg algebra with enough
idempotent corresponds to the opposite dg category by this bijective correspondence.

To be consistent with our notation in the previous section,  we shall
denote by $Gr-K$ the category of graded $K$-modules with degree zero
morphisms and $GR-K$ the graded category with the same objects and where,
for each pair $(V,W)$ of objects, the graded $K$-module of morphisms is
$$\text{HOM}_K(V,W)=\bigoplus_{p\in\mathbb{Z}}\text{HOM}_K^p(V,W),$$
 where $\text{HOM}^p(V,W)=\text{Hom}_{Gr-K}(V,W[p])$ consists of those morphisms
of $K$-modules $f:V\longrightarrow W$ such that $f(V^n)\subseteq V^{n+p}$, for all
$n\in\mathbb{Z}$. Note that $GR-K$ is just the underlying graded category of the
dg category $Dg-K$.

Given a small dg category $\mathcal{A}$, a \emph{graded right
$\mathcal{A}$}-module was defined in \cite{K1} as a graded functor
$M:\mathcal{A}^{op}\longrightarrow GR-K$. The category $\mathcal{G}\mathcal{A}$
has as objects the graded right $\mathcal{A}$-modules and as morphisms their
natural transformations.   Note that, by definition, if  $f:M\longrightarrow N$ is a morphism
in $\mathcal{G}\mathcal{A}$ then  $f_A:M(A)\longrightarrow N(A)$ is a
morphism in $Gr-K$, for each $A\in\text{Ob}(\mathcal{A})$.
The category $\mathcal{G}\mathcal{A}$ comes with a graded functor $?[1]:\mathcal{G}\mathcal{A}\longrightarrow\mathcal{G}\mathcal{A}$  given on
objects by the rule $M[1](A)=M(A)[1]$.
If $a^o\in\mathcal{A}^{op}(A,B)=\mathcal{A}(B,A)$ is a homogeneous element, then
$M[1](a^o):M(A)[1]\longrightarrow M(B)[1]$ is the map $(-1)^{|a|}M(a^o ):M(A)\longrightarrow M(B)$.
We claim that with this definition we have a well-defined graded right
$\mathcal{A}$-module. Indeed, if  $b^o\in\mathcal{A}^{op}(B,C)=\mathcal{A}(C,B)$
is another homogeneous element, then we have

\begin{eqnarray*}
M[1](b^o\circ a^o)&=&(-1)^{|a| |b|}M[1]((a\circ b)^o)\\
&=&(-1)^{|a| |b|}(-1)^{|a|+|b|}M((a\circ b)^o)\\
&=&(-1)^{|a|+|b|}M[(-1)^{|a| |b|}(a\circ b)^o]\\
&=&(-1)^{|a|+|b|}M(b^o\circ a^o)
\end{eqnarray*}
while $$M[1](b^o)\circ M[1](a^o)=(-1)^{|a|+|b|}M(b^o)\circ M(a^o).$$

Therefore $M[1]$ is a well-defined graded right $\mathcal{A}$-module.
Note the discrepancy of the definition of $M[1]$ with the definition in
\cite{K1}.
The assignment $M\rightsquigarrow M[1]$ extends to an auto-equivalence of categories
$\mathcal{G}\mathcal{A}\longrightarrow\mathcal{G}\mathcal{A}$ which acts as
the identity on morphisms. Then the graded category $\text{Gra}\mathcal{A}$
was defined in \cite{K1} as the one having the same objects as
$\mathcal{G}\mathcal{A}$ and as graded $K$-module  of morphisms $\text{Hom}_{\text{Gra}\mathcal{A}}(M,N)=\bigoplus_n\text{Hom}_{\mathcal{G}\mathcal{A}}(M,N[n])$,
where the composition of homogeneous element is given as $g\circ f=g[p]\circ f$,
provided $|f|=p$.

A dg functor $M:\mathcal{A}^{op}\longrightarrow Dg-K$  is called a
\emph{right dg $\mathcal{A}$-module}. It becomes an object of $\text{Gra}\mathcal{A}$
when considering the composition
$\mathcal{A}^{op}\stackrel{M}{\longrightarrow} Dg-K\stackrel{forgetful}{\longrightarrow}GR-K$.
The category $\mathcal{C}_{dg}\mathcal{A}$ (see \cite{K2}), denoted $\text{Dif}\mathcal{A}$ in \cite{K1},
has as objects the right dg $\mathcal{A}$-modules with spaces of morphisms
$\text{Hom}_{\mathcal{C}_{dg}\mathcal{A}}(M,N)=\text{Hom}_{\text{Gra}\mathcal{A}}(M,N)$,
where the differential
$d:\text{Hom}_{\mathcal{C}_{dg}\mathcal{A}}(M,N)\longrightarrow\text{Hom}_{\mathcal{C}_{dg}\mathcal{A}}(M,N)$
acts as $d(f)=d_N\circ f-(-1)^{|f|}f\circ d_M$. Note that one extends $?[1]$
from $\text{Gra}\mathcal{A}$ to $?[1]:\mathcal{C}_{dg}\mathcal{A}
\longrightarrow\mathcal{C}_{dg}\mathcal{A}$, by defining the
differential as $d_{M[1]}=-d_M[1]$.

\begin{Theo1} \label{thm.dg modules versus dg functors}
Let $A=(A,d)$ be a graded algebra with enough idempotents, where
$(e_i)_{i\in I}$ is a fixed distinguished family of orthogonal idempotents,
all of them homogeneous of zero degree and annihilated by $d$. We also view $A$ as a dg category
with $\text{Ob}(A)=I$ as described above. Let $M,N$ be arbitrary objects
of $Dg-A$ and $f:M\longrightarrow N$ be a homogeneous morphism in this category.
The following assertions hold:
\begin{enumerate}
\item The assignments $i\rightsquigarrow\tilde{M}(i):=Me_i$, and
$e_iA^ne_j\longrightarrow\text{Hom}_{GR-K}^n(Me_i,Me_j)$,
$a^o\rightsquigarrow \tilde{M}(a^o)$, where $\tilde{M}(a^o)(x)=(-1)^{|a| |x|}xa$
for each $x\in Me_i$ homogeneous, define a dg functor
$\tilde{M}:A^{op}\longrightarrow\mathcal{C}_{dg}K$ and, hence,
an object of $\mathcal{C}_{dg}A$.
\item If $\tilde{f}=(\tilde{f}_i)_{i\in I}$, where
$\tilde{f}_i:=f_{| Me_i}:\tilde{M}(i)=Me_i\longrightarrow Ne_i=\tilde{N}(i)$,
for each $i\in I$,  then $\tilde{f}$ is a morphism $\tilde{M}\longrightarrow\tilde{N}$
of degree $|f|$ in $\mathcal{C}_{dg}A$.
\item The assignments $M\rightsquigarrow\tilde{M}$ and $f\rightsquigarrow\tilde{f}$
of the two previous assertions define an equivalence of dg categories
$Dg-A\stackrel{\cong}{\longrightarrow}\mathcal{C}_{dg}A$.
\end{enumerate}
\end{Theo1}

\begin{proof}
1) We first prove that $\tilde{M}$ is a graded functor between the underlying
graded categories of $A$ and $\mathcal{C}_{dg}K=Dg-K$, for which we just need
to check that $\tilde{M}(b^o\circ a^o)=\tilde{M}(b^o)\circ\tilde{M}(a^o)$
whenever $a^o\in A^{op}(i,j)=e_iAe_j$ and $b^o\in A^{op}(j,k)=e_jAe_k$ are
homogeneous elements, where $i,j,k\in I$. Note that both sides of the desired
equality are then $K$-linear maps
$\tilde{M}(i)=Me_i\longrightarrow Me_k=\tilde{M}(k)$.
When applying them to a homogeneous element $x\in Me_i$, we have:
\begin{eqnarray*}
\tilde{M}(b^o\circ a^o)(x)&=&(-1)^{|a| |b|}\tilde{M}((ab)^o)(x)\\&=&(-1)^{|a| |b|}(-1)^{|ab| |x|}x(ab)\\
&=&(-1)^{|a| |b|+|a| |x|+|b| |x|}x(ab)\\
&=&(-1)^{|a| |x|}(-1)^{|b| |xa|}(xa)b\\
&=&(-1)^{|a| |x|}\tilde{M}(b^o)(xa)\\
&=&\tilde{M}(b^o)[(-1)^{|a| |x|} xa]\\
&=&(\tilde{M}(b^o)\circ\tilde{M}(a^o))(x)
\end{eqnarray*}
The desired equality then holds due to the fact that $M$ is a right $A$-module.

In order to see that $\tilde{M}$ is a dg functor, we need to check that it
commutes with the differentials on Hom spaces  of $A^{op}$ and $\mathcal{C}_{dg}K$.
That is,  that if $a\in e_iAe_j=A^{op}(i,j)$ is a homogeneous element, then
$\tilde{M}(d(a)^o)=d_{\text{HOM}_K(\tilde{M}(i),\tilde{M}(j))}(\tilde{M}(a))$.
To check this equality, we evaluate the two maps on a homogeneous element
$x\in\tilde{M}(i)=Me_i$. We then have:
\begin{eqnarray*}
d_{\text{HOM}_K(\tilde{M}(i),\tilde{M}(j))}(\tilde{M}(a^o))(x)&=&
[d_{Me_j}\circ\tilde{M}(a^o)-(-1)^{|\tilde{M}(a^o)|}\tilde{M}(a^o)\circ d_{Me_i}](x)\\
&=&d_M((-1)^{|a| |x|}xa)-(-1)^{|a|}(-1)^{|a| |d_M(x)|}d_M(x)a\\
&=&(-1)^{|a| |x|}(d_M(xa)-d_M(x)a)\\
&=&(-1)^{|a| |x|}(-1)^{|x|}xd(a)\\
&=&(-1)^{(|a|+1) |x|}xd(a)\\
&=&(-1)^{|d(a)| |x|}xd(a)\\
&=&\tilde{M}(d(a)^o)(x),
\end{eqnarray*}
as desired.

\bigskip

2) In order to prove this assertion, we first show that $\widetilde{M[1]}$
is isomorphic to $\widetilde{M}[1]$. On objects, we have
$\widetilde{M[1]}(i)=M[1]e_i=Me_i[1]=\widetilde{M}[1](i)$, for each $i\in I$.
On the other hand, if $a^o\in A^{op}(i,j)=e_iAe_j$ and
$x\in\widetilde{M[1]}(e_i)=Me_i$ are homogeneous elements, then we have that
$\widetilde{M[1]}(a^o)(x)=(-1)^{|a| |x|_{M[1]}}xa$, where $|x|_{M[1]}$ denotes
the degree of $x$ as an element of $M[1]$. We know that $|x|_{M[1]}=|x|-1$,
where $|x|$ is the degree of $x$ as an element of $M$. Therefore we have
$\widetilde{M[1]}(a^o)(x)=(-1)^{|a| (|x|-1)}xa$. On the other side, by
definition of the shift in $Gra-K$, we have that
$\widetilde{M}[1](a^o)=(-1)^{|a|}\widetilde{M}(a^o)$. It follows that $$(\widetilde{M}[1])(a^o)(x)=(-1)^{|a|}\widetilde{M}(a^o)(x)=(-1)^{|a|}(-1)^{|a| |x|}xa.$$
As a consequence, we have that $\widetilde{M[1]}(a^o)(x)=(\widetilde{M}[1])(a^o)(x)$.

Let $a^o\in A^{op}(i,j)=e_iAe_j$ be homogeneous and assume that $|f|=n$,
that is, that $f:M\longrightarrow N[n]$ is a morphism in $Gr-A$.
We need to prove that the following diagram in $Gr-K$ is commutative:
$$\xymatrix{\widetilde{M}(i)\ar[r]^{\widetilde{f}_i}\ar[d]_{\widetilde{M}(a^\circ)}&
\widetilde{N}[n](i)\ar[d]^{(\widetilde{N}[n])(a^\circ)}\\
\widetilde{M}(j)\ar[r]^{\widetilde{f}_j}&
\widetilde{N}[n](j)}$$
Indeed, for each $x\in\widetilde{M}(i)=Me_i$ homogeneous, we have
$$(\widetilde{f}_j\circ\widetilde{M}(a^o))(x)=\widetilde{f}_j((-1)^{|a| |x|_M}xa)=(-1)^{|a| |x|_M}f(xa), $$
and, using the previous paragraph, we also have
$$(\widetilde{N}[n](a^o)\circ \widetilde{f}_i)(x)= (\widetilde{N[n]}(a^o)\circ
\widetilde{f}_i)(x)=\widetilde{N[n]}(a^o)(f(x))=(-1)^{|a| |f(x)_{N[n]}|}f(x)a.$$
Note that $|f(x)|_{N[n]}=|x|_M$ because $f:M\longrightarrow N[n]$ is a morphism of
degree zero. On the other hand, the product $f(x)a$ is considered in the graded
right $A$-module $N[n]$. The commutativity of the desired diagram follows
from that fact that $f:M\longrightarrow N[n]$ is a morphism of right $A$-modules.

\bigskip

3) Since the definition of  $(?)^{\widetilde{}}$ on morphisms is the `identity', i.e.
$\widetilde{f}_i:\widetilde{M}(i)=Me_i\longrightarrow\widetilde{N}(i)=Ne_i$  is just the
restriction of $f$ to $Me_i$, we readily see that we have a well-defined
$K$-linear functor between the underlying graded categories
$(?)^{\widetilde{}}:GR-A\longrightarrow Gra-A$. This functor is clearly graded
since $|\widetilde{f}|=|f|$ for each homogeneous morphism $f$ in $GR-A$.
Moreover the differential of the graded $K$-module $Me_i=\widetilde{M}(i)$ is
the same when coming from $M$ that when coming from $\widetilde{M}$. Again the
fact that the action of $(?)^{\widetilde{}}$ on morphisms is the identity
implies that $(?)^{\widetilde{}}$ commutes with the
differentials on Hom spaces. That is, it is actually a dg functor
$Dg-A\longrightarrow\mathcal{C}_{dg}A$.

On the other hand, the `identity' condition on the action on morphisms
immediately implies that $(?)^{\widetilde{}}$ is a faithful functor.
We shall now prove that $(?)^{\widetilde{}}$ is full.
If $\psi :\widetilde{M}\longrightarrow\widetilde{N}[n]$
is a morphism in $\mathcal{G}A$, then for each $a^o\in A^{op}(i,j)=e_iAe_j$,
we have that $\psi_j\circ\widetilde{M}(a^o)=\widetilde{N}[n](a^o)\circ\psi_i$.
When applying both members of this equality to an element $x\in\widetilde{M}(i)=Me_i$, we get that
$$(\psi_j\circ\widetilde{M}(a^o))(x)=\psi_j[(-1)^{|a| |x|_M}xa]= (-1)^{|a| |x|_M}\psi_j(xa)$$
while
$$(\widetilde{N}[n](a^o)\circ\psi_i)(x)=(\widetilde{N[n]}(a^o)\circ\psi_i)(x)=
(-1)^{|a| |\psi_i(x)|_{N[n]}}\psi_i(x)a. $$
Bearing in mind that $|\psi_i(x)|_{N[n]}=|x|_M$,
we conclude that $\psi_j(xa)=\psi_i(x)a$. This means that if we
define $f:M=\bigoplus_{i\in I} Me_i\longrightarrow\bigoplus_{i\in I} N[n]e_i=N[n]$ as the direct sum of the
$\psi_i:Me_i=\widetilde{M}(i)\longrightarrow\widetilde{N}[n](i)=N[n]e_i$, then
$f$ is a morphism of graded right $A$-modules such that $\widetilde{f}=\psi$.
We showed that the functor $(?)^{\widetilde{}}$ is also full.

It remains to check that $(?)^{\widetilde{}}$ is a dense functor.
Let $F$ be an object of $\mathcal{C}_{dg}A$ and consider the dg  $K$-module $M_F:=\bigoplus_{i\in I}F(i)$.
We endow $M_F$ with a structure of graded right $A$-module as follows.
Given $x\in F(i)$ and $a^o\in A^{op}(j,k)=e_jAe_k$, we define
$xa:=(-1)^{|a| |x|}\delta_{ij}F(a^o)(x)$, where $\delta_{ij}$ is the Kronecker symbol.
Then one extends this multiplication by $K$-linearity. In order to see
that this rule gives $M:=M_F$ the structure of a graded right $A$-module,
we just need to consider $x\in F(i)=Me_i$, $a\in e_iAe_j$ and $b\in e_jAe_k$
homogeneous elements and check that $x(ab)=(xa)b$. Indeed,
bearing in mind that $(ab)^o=(-1)^{|a| |b|}(b^o\circ a^o)$ when looking
at the elements of $A$ as morphism in the underlying graded $K$-category,
we have an equality
\begin{eqnarray*}
x(ab)&=&(-1)^{|ab| |x|}F((ab)^o)(x)\\
&=&(-1)^{|ab| |x|}(-1)^{|a| |b|}F(b^o\circ a^o)(x)\\
&=&(-1)^{|a| |x|+|b| |x|+|a| |b|}[F(b^o)\circ F(a^o)](x)\\
&=&(-1)^{|b| |xa|}(-1)^{|a| |x|}F(b^o)(F(a^o)(x))\\
&=&(-1)^{|b| |xa|}F(b^o)(xa)=(xa)b,
\end{eqnarray*}
which shows that $x(ab)=(xa)b$ as desired.

We claim that the diferential $d=\oplus d_{F(i)}:M=\oplus_{i\in I} F(i)\longrightarrow\oplus_{i\in I} F(i)=M$ satisfies Libniz rule, thus making $M$ into a right dg $A$-module. To see this, note that since $F:A^{op}\longrightarrow\mathcal{C}_{dg}K$ is a dg functor we have an equality $\delta (F(a^o))=F(d(a)^o)$, for any morphism $a^o\in A^{op}(i,j)=e_iAe_j$, where $d:e_iAe_j\longrightarrow e_iAe_j$ is the restriction of the differential of $A$ and $\delta :\text{HOM}_K(F(i),F(j))=\text{HOM}_K(Me_i,Me_j)\longrightarrow \text{HOM}_K(Me_i,Me_j)$ is the differential on $Hom$ spaces of $\mathcal{C}_{dg}K$. We then have that $F(d(a)^o)=d_{Me_j}\circ F(a^o)-(-1)^{|F(a^o)|}F(a^o )\circ d_{Me_i}$. Bearing in mind that $|F(a^o)|=|a|$, when making act both members of the last equality on a homogeneous element $x\in Me_i$, we have $$(-1)^{|d(a)| |x|}xd(a)=(-1)^{|a| |x|}d_{Me_j}(xa)-(-1)^{|a|}(-1)^{|d_{Me_i}(x)| |a|}d_{Me_i}(x)a. $$ Cancelling $(-1)^{|a| |x|}$ from this equality, we get that $(-1)^{|x|}xd(a)=d_M(xa)-d_M(x)a$, from which Leibniz rule immediately follows.

The fact that $\tilde{M_F}\cong F$
follows in a straightforward way, and hence $(?)^{\tilde{}}$ is a dense functor.
\end{proof}

\begin{Rem} \label{rem.projective versus representable}
Note that the equivalence of categories $(?)^{\tilde{}}:Dg-A\longrightarrow\mathcal{C}_{dg}A$ given
by Theorem~\ref{thm.dg modules versus dg functors} takes $e_iA$ to the representable dg $A$-module
$i^\wedge:A(?,i):A^{op}\longrightarrow\mathcal{C}_{dg}K=Dg-K$
(see Example \ref{ex.the regular dg bifunctor} and its proof).
\end{Rem}

\section{Right versus left dg modules} \label{section.right versus left}
From the definition of left dg $A$-module we get the following:

\begin{Lemma1} \label{lem.left versus right dg module}
If $(M,d_M)$ is a left dg $A$-module, then it is a right dg $A^{op}$-module with the multiplication map $M\otimes A^{op}\longrightarrow M$ defined as $(x,a^o)\rightsquigarrow xa^o:=(-1)^{|a| |x|}ax$, for all homogeneous elements $a\in A$ and $x\in M$. Conversely, if $(M,d_M)$ is a right dg $A^{op}$-module, then it is a left dg $A$-module, where the multiplication map $A\otimes M\longrightarrow M$ takes $(a,x)\rightsquigarrow ax:=(-1)^{|a| |x|}xa^o$, for $a$ and $x$ as above.
\end{Lemma1}
\begin{proof}
We just prove the first implication, the reverse one being then clear. Given $x\in M$ and $a,b\in A$ homogeneous elements, we have:

$$(xa^o)b^o=(-1)^{|xa^o| |b|}b(xa^o)=(-1)^{(|x|+|a|)|b|)}(-1)^{|x| |a|}b(ax)$$
and
$$x(a^ob^o)=(-1)^{|a| |b|}x(ba)^o=
(-1)^{|a| |b|}(-1)^{|x| |ba|}(ba)x=(-1)^{|a| |b|}(-1)^{|x|(|b|+|a|)}(ba)x.$$
Therefore we have $(xa^o)b^o=x(a^ob^o)$. Since up to here the differential has played no role, we have actually proved that any graded left $A$-module is a graded right $A^{op}$-module.

We next check that the differential of $M$ as a left $A$-module satisfies Leibniz
rule as a right $A^{op}$-module. We have that
$$d_M(xa^o)=(-1)^{|x| |a|}d_M(ax)=(-1)^{|x| |a|}[d(a)x+(-1)^{|a|}ad_M(x)],$$
while we have $$d_M(x)a^o+(-1)^{|x|}xd(a)^o=(-1)^{(|x|+1)|a|}ad_M(x)+(-1)^{|x|}(-1)^{|x|(|a|+1)}d(a)x,  $$
so that $d_M(xa^o)=d_M(x)a^o+(-1)^{|x|}xd(a)$, for all homogeneous elements
$x\in M$ and $a^o\in A^{op}$.
\end{proof}

As with graded right $A$-modules, one first defines the category $A-Gr$ of graded left
$A$-modules, where morphisms are the graded morphisms of zero degree. By the sign trick of the previous lemma this category should be canonically identified with $Gr-A^{op}$. We next need to define a shift functor  $?[1]:A-Gr\longrightarrow A-Gr$ which,  viewed as a functor $Gr-A^{op}\longrightarrow Gr-A^{op}$,   coincides with the shift for graded right modules (see the paragraph after the proof of Lemma \ref{lem.tensor product dg algebras}). This forces the definition of the multiplication map $A\otimes M[1]\longrightarrow M[1]$ ($(a,x)\rightsquigarrow a\cdot x$). Indeed we will have $a\cdot x=(-1)^{|a| |x|_{M[1]}}xa^o$. But the multiplication $xa^o$ is the same in $M[1]$ and $M$, due the the definition of $?[1]$ for graded right $A^{op}$-modules. Then in $M$ we have $xa^o=(-1)^{|a| |x|_M}ax$, where $ax$ is given by the multiplication $A\otimes M\longrightarrow M$. We then get that the multiplication map in $M[1]$ is given by $$a\cdot x=(-1)^{|a| |x|_{M[1]}}(-1)^{|a| |x|_M}ax=(-1)^{|a|}ax,$$ where $ax$ is the multiplication in $M$.

This readily gives
a graded $K$-category $A-GR$ whose objects are the objects of $A-Gr$ and where the space
of morphisms $HOM_{A^{op}}(M,N)$ between two objects $M$ and $N$ is graded in such a way that
the $n$-th homogeneous component is $HOM^n_{A^{op}}(M,N)=Hom_{A-Gr}(M,N[n])$.
Note that, viewing an element $f\in Hom_{A-Gr}(M,N[n])$, as a morphism
$f:M\longrightarrow N$ of degree $n$, we have
$$f(ax)=((-1)^{|a|})^naf(x)=(-1)^{|f| |a|}af(x),$$
for all homogeneous elements $a\in A$ and $x\in M$. This is due to the fact
that the multiplication map $A\otimes N[n]\longrightarrow N[n]$ acts as
$$a\cdot y=((-1)^{|a|})^nay=(-1)^{na}ay,$$ for all homogeneous elements $a\in A$
and $y\in N[n]$, where $ay$ is the product in $N$.

\begin{Rem} \label{rem.forgetful for graded left modules}
We have an obvious forgetful functor $A-Gr\longrightarrow A-Mod$ acting as the identity on objects and morphisms. However, we have such a functor for the category $A-GR$ only in case $A$ is evenly graded (i.e. $A^{2k+1}=0$, for all $k\in\mathbb{Z}$). In the general case one has a forgetful `pseudo-functor' $A-GR\longrightarrow A-\text{Mod}$. It acts as the identity on objects, but takes $f\rightsquigarrow\hat{f}$, where $\hat{f}(x)=(-1)^{|f| |x|}f(x)$, for all homogeneous elements $f\in\text{HOM}_{A^{op}}(M,N)$ and $x\in M$ (note that $\hat{f}$ is a morphism in $\text{Mod}-A$). This assignment satisfies the equality $\widehat{g\circ f}=(-1)^{|f| |g|}\hat{g}\circ\hat{f}$, for all homogeneous morphisms $f,g$ in $A-GR$. The ultimate reason for this disruption is that $A^{op}$ is not the opposite algebra of $A$ as an ungraded algebra.
\end{Rem}

We have essentially proved the following expected result.

\begin{Prop} \label{prop.left-versus-right}
Let $A$ be a dg algebra with enough idempotents. The following data give a
dg category $A-Dg$ which is equivalent to the dg category $Dg-A^{op}$:

\begin{enumerate}
\item[-] The objects of $A-Dg$ are the left dg $A$-modules (see Definition \ref{def.dg module});
\item[-] The morphisms in $A-Dg$ are defined as in the category $A-GR$;
\item[-] For each pair $(M,N)$ of objects, the differential
$d:HOM_{A^{op}}(M,N)\rightarrow HOM_{A^{op}}(M,N)$ on $Hom$ spaces is
defined by the rule $d(f):=d_N\circ f-(-1)^{|f|}f\circ d_M$.
\end{enumerate}
\end{Prop}
\begin{proof}

Once we identify the objects of $A-Dg$ with those of $Dg-A^{op}$, we
only need to identify the spaces of morphisms in $A-GR$ and $GR-A^{op}$
for the differential $d:HOM_{A^{op}}(M,N)\longrightarrow HOM_{A^{op}}(M,N)$
is just the restriction of that of the dg $K$-module $\text{HOM}_K(M,N)$. Indeed,
if $f:M\longrightarrow N$ is a homogeneous morphism in $A-GR$,
then
$$
f(xa^o)=(-1)^{|a| |x|}f(ax)=(-1)^{|a| |x|}(-1)^{|a| |f|}af(x)=(-1)^{|a| |f(x)|}af(x)=f(x)a^o.
$$

\end{proof}

As in the case of right modules, the shift functor $?[1]:A-Gr\longrightarrow A-Gr$ extends to a functor $?[1]:A-Dg\longrightarrow A-Dg$ such that $d (f[1])=-d(f)[1]$, for all homogeneous $f\in\text{HOM}_{A^{op}}(M,N)$, where $d$ is the differential on $Hom$ spaces.

\section{Dg bimodules}

\begin{Def} \label{def.dg bimodule}
Let $A$ and $B$ be dg algebras with enough idempotents.
A \emph{graded $A-B-$bimodule} is graded $K$-module $M$ together with the following data:
\begin{enumerate}
\item A morphism of graded $K$-vector spaces $\mu_{left}:A\otimes M\longrightarrow M$ ($a\otimes x\rightsquigarrow ax$) making $M$ into a graded left $A$-module,
\item and a morphism of graded $K$-vector spaces $\mu_{right}:M\otimes B\longrightarrow M$ ($x\otimes b\rightsquigarrow xb$) making $M$ into a graded right $B$-module,
\item such that $(ax)b=a(xb)$, for all  $(a,x,b)\in A\times M\times B$.
\end{enumerate}
  A \emph{differential graded (dg) $A-B-$bimodule} is a pair $(M,d_M)$ consisting of a graded $A-B-$bimodule $M$ and a morphism $d_M:M\longrightarrow M$ in $GR-K$ of degree $+1$, called the differential,  such that $d_M\circ d_M=0$ and  $$d_M(axb)=d_A(a)xb+(-1)^{|a|}ad_M(x)b+(-1)^{|a|+|x|}axd_B(b),$$ for all homogeneous elements $a\in A$, $x\in M$ and $b\in B$. This latter formula is called Leibniz rule (for the the dg bimodule).
	
	\end{Def}
	
As in the case of right or left dg modules, we successively consider the category $A-Gr-B$ of graded $A-B-$bimodules with morphisms of zero degree, the graded category $A-GR-B$, where $\text{HOM}_{A-B}(M,N):=\text{Hom}_{A-GR-B}(M,N)$ is the graded $K$-module with $n$-th homogeneous component $\text{HOM}^n_{A-B}(M,N)=\text{Hom}_{A-Gr-B}(M,N[n])$, for each $n\in\mathbb{Z}$, and where $g\circ f:=g[p]\circ f$ in case $f$ and $g$ are homogeneous elements of $\text{HOM}_{A-B}(M,N)$, with $|f|=p$. Finally, $A-Dg-B$ will denote the dg category whose objects are the dg $A-B-$bimodules with morphisms as in $A-GR-B$.

\begin{Prop} \label{prop.dg bimodules as dg modules}
Let $A$ and $B$ be dg algebras with enough idempotents. The following three terms `are' synonymous:
\begin{enumerate}
\item Dg $A-B-$bimodule;
\item Right dg $B\otimes A^{op}-$module;
\item Left dg $A\otimes B^{op}$-module.
\end{enumerate}
In particular, there are equivalences of dg categories $$Dg-(B\otimes A^{op})\cong A-Dg-B\cong (A\otimes B^{op})-Dg.$$
\end{Prop}
\begin{proof}
We define the map $\Phi :A\otimes B^{op}\longrightarrow (B\otimes A^{op})^{op}$ by the rule $\Phi (a\otimes b^o)=(-1)^{|a| |b|}(b\otimes a^o)^o$, for all homogeneous elements $a\in A$ and $b\in B$. We will prove that $\Phi$ is an isomorphism of dg algebras, which will imply that left dg $A\otimes B^{op}$ is synonymous of right dg $B\otimes A^{op}$-module using  Proposition \ref{prop.left-versus-right}. We clearly have that $\Phi$ a morphism (of zero degree) of graded $K$-modules. Moreover, if $a_1,a_2\in A$ and $b_1,b_2\in B$ are homogeneous elements, then we have equalities

\begin{eqnarray*}
\Phi [(a_1\otimes b_1^o)\cdot (a_2\otimes b_2^o)]&=&(-1)^{|b_1| |a_2|}\Phi (a_1a_2\otimes b_1^ob_2^o)\\
&=&(-1)^{|b_1| |a_2|+|b_1| |b_2|}\Phi (a_1a_2\otimes (b_2b_1)^o)\\
&=&(-1)^{|b_1| |a_2|+|b_1| |b_2|+(|a_1|+|a_2|)(|b_1|+|b_2|)}[b_2b_1\otimes (a_1a_2)^o]^o
\end{eqnarray*}

and

\begin{eqnarray*}
\Phi (a_1\otimes b_1^o)\lefteqn{\cdot\Phi (a_2\otimes b_2^o)=}\\
&=&(-1)^{|a_1| |b_1| +|a_2| |b_2|}(b_1\otimes a_1^o)^o\cdot (b_2\otimes a_2^o)^o\\
&=&(-1)^{|a_1| |b_1| +|a_2| |b_2|+(|b_1|+|a_1|)(|b_2|+|a_2|)}[(b_2\otimes a_2^o)(b_1\otimes a_1^o)]^o\\
&=&(-1)^{|a_1| |b_1| +|a_2| |b_2|+(|b_1|+|a_1|)(|b_2|+|a_2|)+|a_2| |b_1|}[b_2b_1\otimes a_2^oa_1^o]^o\\
&=&(-1)^{|a_1| |b_1| +|a_2| |b_2|+(|b_1|+|a_1|)(|b_2|+|a_2|)+|a_2| |b_1|+|a_1| |a_2|}[b_2b_1\otimes (a_1a_2)^o]^o.
\end{eqnarray*}
We compare the signs of the two expressions.
\begin{eqnarray*}
|a_2| |b_1|+\lefteqn{|b_1| |b_2|+|a_1| |b_1|+|a_1| |b_2|+|a_2| |b_1|+|a_2| |b_2|=}\\
&=&|b_1| |b_2|+|a_1| |b_1|+|a_1| |b_2|+|a_2| |b_2|\\
&=&|a_1| |b_1|+|a_2| |b_2|+|b_1| |b_2|+|b_1| |a_2|+|a_1| |b_2|+|a_1| |a_2|+|a_2| |b_1|+|a_1| |a_2|
\end{eqnarray*}
We conclude that $\Phi$ is a (clearly bijective) homomorphism of graded algebras.
In order to prove that it is actually an isomorphism of dg algebras, we need to
check that it is compatible with the differentials. Indeed, if $a\in A$ and $b\in B$
are homogeneous elements, then we have
\begin{eqnarray*}
[d\circ\Phi ](a\otimes b^o)&=&(-1)^{|a| |b|}d[(b\otimes a^o)^o]\\
&=&(-1)^{|a| |b|}[(d(b)\otimes a^o)^o+(-1)^{|b|}(b\otimes d(a)^o)^o]\\
&=&(-1)^{|a| |b|}(d(b)\otimes a^o)^o+(-1)^{(|a|+1) |b|}(b\otimes d(a)^o)^o\\
&=&(-1)^{(|a|+1) |b|}(b\otimes d(a)^o)^o+(-1)^{|a|+|a| (|b|+1)}(d(b)\otimes a^o)^o\\
&=&\Phi [d(a)\otimes b^o+(-1)^{|a|}a\otimes d(b)^o]\\
&=&[\Phi\circ d](a\otimes b^o)
\end{eqnarray*}
which shows that $\Phi\circ d=d\circ\Phi$ and, hence, that $\Phi$ is an isomorphism of dg algebras.

If $M$ is a dg $A-B-$bimodule, then we define a multiplication map
$M\otimes (B\otimes A^{op})\longrightarrow M$ which takes
$x\otimes b\otimes a^o\rightsquigarrow x(b\otimes a^o):=(-1)^{(|x|+|b|) |a|}axb$,
whenever $x\in M$, $b\in B$ and $a\in A$ are homogeneous elements. We claim that
this map endows $M$ with a structure of graded right $B\otimes A^{op}$-module.
For this we just need to check the equality
$x[(b_1\otimes a_1^o)(b_2\otimes a_2^o)]=[x(b_1\otimes a_1^o)](b_2\otimes a_2^o)$,
for all homogeneous elements $a_1,a_2\in A$, $b_1,b_2\in B$ and $x\in M$.
Indeed we have:
\begin{eqnarray*}
x[(b_1\otimes a_1^o)(b_2\otimes a_2^o)]&=&(-1)^{|a_1| |b_2|}x[b_1b_2\otimes a_1^oa_2^o]\\
&=&(-1)^{|a_1| |b_2|+|a_1| |a_2|}x[b_1b_2\otimes (a_2a_1)^o]\\
&=&(-1)^{|a_1| |b_2|+|a_1| |a_2|+(|x|+|b_1|+|b_2|)(|a_1|+|a_2|)}a_2a_1xb_1b_2\\
&=&(-1)^{(|x|+|b_1|) |a_1|+(|a_1|+|x|+|b_1|+|b_2|) |a_2|}a_2a_1xb_1b_2\\
&=&(-1)^{(|x|+|b_1|) |a_1|}a_1xb_1(b_2\otimes a_2^o)\\
&=&[x(b_1\otimes a_1^o)](b_2\otimes a_2^o)
\end{eqnarray*}
We conclude that the above mentioned multiplication map endows $M$ with a structure
of graded right $B\otimes A^{op}$-module. We finally check that the differential
$d_M:M\longrightarrow M$ satisfies Leibniz rule $d_M[x(b\otimes a^o)]=d_M(x)(b\otimes a^o)+(-1)^{|x|}xd(b\otimes a^o)$,
for all homogeneous elements $x\in M$, $b\in B$ and $a\in A$. Indeed we have an equality
\begin{eqnarray*}
d_M(x)\lefteqn{(b\otimes a^o)+(-1)^{|x|}xd(b\otimes a^o)=}\\
&=&(-1)^{(|x|+1+|b|) |a|}ad_M(x)b+(-1)^{|x|}x[d(b)\otimes a^o+(-1)^{|b|}b\otimes d(a)^o]\\
&=&(-1)^{(|x|+1+|b|) |a|}ad_M(x)b+(-1)^{|x|}(-1)^{(|x|+|b|+1) |a|}axd(b)\\
&&+(-1)^{|x|+|b|}(-1)^{(|x|+|b|)(|a|+1)}d(a)xb\\
&=&(-1)^{(|x|+|b|) |a|}[(-1)^{|a|}ad_M(x)b+(-1)^{|x|+|a|}axd(b)+(-1)^{2(|x|+|b|)}d(a)xb]\\
&=&(-1)^{(|x|+|b|) |a|}d(axb)\\
&=&d_M[x(b\otimes a^o)]
\end{eqnarray*}
where the last equation holds by the definition of right $B\otimes A^{op}$-module
structure on $M$. Then $(M,d_M)$ is a right dg $B\otimes A^{op}$-module.

Obviously, one can reverse the arguments step by step, so that if $X$ is a
right dg $B\otimes A^{op}$-module, then it is also a dg $A-B-$bimodule when
taken with multiplication $axb=(-1)^{(|x|+|b|) |a|}x(b\otimes a^o)$, for all homogeneous
elements $a\in A$, $x\in M$ and $b\in B$.
\end{proof}

\begin{Rem}
We emphasize the structures of right dg $B\otimes A^{op}$-module and left
dg $A\otimes B^{op}$-module coming from Proposition \ref{prop.dg bimodules as dg modules} and its proof.
If $M$ is a dg $A-B-$bimodule and $a\in A$, $b\in B$ and $x\in M$ are homogeneous
elements, then we have:
\begin{enumerate}
\item $x(b\otimes a^o)=(-1)^{(|x|+|b|) |a|}axb$;
\item $(a\otimes b^o)x=(-1)^{|b| |x|}axb$.
\end{enumerate}
Indeed the first equality appears in the proof of  Proposition \ref{prop.dg bimodules as dg modules} and then,
by Proposition \ref{prop.left-versus-right}, we have a structure of left dg
$(B\otimes A^{op})^{op}$-module on $M$. Using then the isomorphism
$\Phi :A\otimes B^{op}\longrightarrow (B\otimes A^{op})^{op}$ from the proof
of  , we make $M$ into a left dg $A\otimes B^{op}$-module.
The reader can check that this module structure is given by equality 2.
\end{Rem}

\begin{Example}
If $A$ is a dg algebra with enough idempotents, then it is a dg $A-A-$bimodule with its
canonical multiplication and differential. We will call it the \emph{regular dg bimodule}.
\end{Example}

\section{Homotopy category and derived category}

As in the case of small dg categories, given a dg algebra with enough idempotents $A$,
the 0-cycle category $Z^0(Dg-A)$, denoted by $\mathcal{C}(A)$ in the sequel,
has two structures to take into account. It is a bicomplete abelian category,
where the exact sequences are those sequences
$0\rightarrow L\longrightarrow M\longrightarrow N\rightarrow 0$ of morphisms in $\mathcal{C}(A)$
which are exact as sequences in $Gr-A$. But, even more relevant to us, it
has a Quillen exact structure where the conflations (=admissible short
exact sequences) are the short exact sequences which split in $Gr-A$ (see \cite{B} and \cite{K3} for the axioms and details about exact categories). It is
called the \emph{semi-split exact structure}.  We are now going to give an explicit description
of the projective (=injective) objects for this exact structure.

\begin{Lemma} \label{lem.canonical conflation}
Any conflation in $\mathcal{C}(A)$ is isomorphic to one whose
underlying exact sequence in $Gr-A$ is
$0\rightarrow L\stackrel{\begin{pmatrix}1 \\0 \end{pmatrix}}{\longrightarrow}
L\oplus N\stackrel{\begin{pmatrix}0 & 1 \end{pmatrix}}{\longrightarrow}N\rightarrow 0$,
where the differential of $L\oplus N$ is of the form
$\delta =\begin{pmatrix}d_L & s\\ 0 & d_N \end{pmatrix}$, for some morphism
$s:N\longrightarrow L$ of degree $1$ in $GR-A$ such that $d_L\circ s+s\circ d_N=0$.
\end{Lemma}

\begin{proof}
By the definition of conflations, the underlying exact sequence in $Gr-A$ of
such a conflation is always as indicated, where $L$ and $N$ are right dg $A$-modules. We initially put
$\delta =\begin{pmatrix}d_{11} & d_{12}\\ d_{21} & d_{22} \end{pmatrix}$
($\begin{pmatrix}x \\ y \end{pmatrix}\rightsquigarrow
\begin{pmatrix}d_{11}(x)+d_{12}(y) \\ d_{21}(x)+d_{22}(y)\end{pmatrix}$).
Since $L\stackrel{\begin{pmatrix}1 \\0 \end{pmatrix}}{\longrightarrow}L\oplus N$
should be an element of $Z^0(HOM_A(L,L\oplus N))$, we should have
$\delta \circ\begin{pmatrix}1 \\0 \end{pmatrix}-\begin{pmatrix}1 \\0 \end{pmatrix}\circ d_L=0$.
From this equality we get that $d_{21}=0$ and $d_{11}=d_L$. Using that
$L\oplus N\stackrel{\begin{pmatrix}0 & 1 \end{pmatrix}}{\longrightarrow}N$ is in $Z^0(HOM_A(L\oplus N,N))$
we then get $d_N\circ\begin{pmatrix}0 & 1 \end{pmatrix}-\begin{pmatrix}0 & 1 \end{pmatrix}\circ\delta$,
from which we get that $d_{22}=d_N$. Finally, from the equality $\delta\circ\delta =0$
we get that $d_L\circ s+s\circ d_N=0$, with $s=d_{12}$.
\end{proof}

\begin{Rem} \label{rem.cone}
If $f:M\longrightarrow N$ is a morphism in $\mathcal{C}(A)$,
then one can consider the split exact sequence
$0\rightarrow N\stackrel{\begin{pmatrix}1 \\0 \end{pmatrix}}{\longrightarrow}
N\oplus M[1]\stackrel{\begin{pmatrix}0 & 1 \end{pmatrix}}{\longrightarrow}M[1]\rightarrow 0$ in $Gr-A$.
Viewing $f$ as morphism of degree $+1$ from $M[1]$ to $N$, Lemma~\ref{lem.canonical conflation}
makes $N\oplus M[1]$
into a right dg $A$-module with differential
$\delta =\begin{pmatrix}d_N & f\\ 0 & d_{M[1]} \end{pmatrix}$.
This dg $A$-module is known as the \emph{cone of $f$} and will be
denoted by $C(f)$ in the sequel. Note that we have an associated
conflation $0\rightarrow N\longrightarrow C(f)\longrightarrow M[1]\rightarrow 0$ in $\mathcal{C}(A)$.
\end{Rem}

\begin{Prop} \label{prop.projective-injective in C(A)}
For a right dg $A$-module $P$, the following assertions are equivalent:

\begin{enumerate}
\item $P$ is projective with respect to the semi-split exact structure;
\item $P$ is injective with respect to the semi-split exact structure;
\item $P$ is isomorphic to a direct summand of a cone $C(1_M)$, for a right dg $A$-module $M$.
\end{enumerate}
Such a $P$ is acyclic. In particular $\mathcal{C}(A)$ is a Frobenius exact category with the semi-split structure.
\end{Prop}

\begin{proof}
The acyclic condition $C(1_M)$ is well-knosn. It is then enough to check that $C(1_M)$ is projective and injective, for
each right dg $A$-module $M$. Once this is proved, if $P$ is projective
(resp. injective) object, then the canonical conflation
$0\rightarrow P[-1]\longrightarrow C(1_{P[-1]})\rightarrow P\rightarrow 0$
(resp. $0\rightarrow P\longrightarrow C(1_{P})\rightarrow P[1]\rightarrow 0$)
must split in $\mathcal{C}(A)$, and the rest of the proof will be trivial.

By Lemma \ref{lem.canonical conflation}, any deflation (=admissible
epimorphism) in $\mathcal{C}(A)$ can be identified with
$\begin{pmatrix} 0 & 1\end{pmatrix}:L\oplus N\longrightarrow N$, where $L\oplus N$
is made into a right dg $A$-module with the differential
$\delta =\begin{pmatrix}d_L & s\\ 0 & d_N \end{pmatrix}$ described there.
If now $f:C(1_M)\longrightarrow N$ is any morphism in $\mathcal{C}(A)$,
then, viewed in $Gr-A$, it is a morphism
$f=\begin{pmatrix} u & v \end{pmatrix}:M\oplus M[1]\longrightarrow N$
such that
$\begin{pmatrix} u & v \end{pmatrix}\circ\begin{pmatrix} d_M & 1_{M[1]}\\ 0 & d_{M[1]} \end{pmatrix}
=d_N\circ\begin{pmatrix} u & v \end{pmatrix}$.
The second component of this equality gives that $u+v\circ d_{M[1]}=d_N\circ v$,
which we express as $u=\hat{d}(v)$, where $v$ is viewed as a morphism of degree
$-1$ in $GR-A$ and $\hat{d}$ is the internal differential
$\hat{d}:HOM_A(M,N)\longrightarrow HOM_A(M,N)$ in the dg category $Dg-A$.
We will now check that the morphism in $Gr-A$ given in matrix form as
$\alpha =\begin{pmatrix} s\circ v & 0\\ \hat{d}(v) & v \end{pmatrix}:M\oplus M[1]\longrightarrow L\oplus N$
is a morphism in $\mathcal{C}(A)$,
$\alpha :C(1_M)\longrightarrow (L\oplus N,\delta )$, such that
$\begin{pmatrix} 0 & 1 \end{pmatrix}\circ\alpha =\begin{pmatrix} u & v \end{pmatrix}=f$.
In order to see that $\alpha$ is a morphism in $\mathcal{C}(A)$,
we just need to check the equality
$$\begin{pmatrix} s\circ v & 0 \\ \hat{d}(v) & v\end{pmatrix}\circ\begin{pmatrix} d_M & 1_{M[1]} \\ 0 & d_{M[1]} \end{pmatrix}=\begin{pmatrix} d_L & s \\ 0 & d_N\end{pmatrix}\circ \begin{pmatrix} s\circ v & 0 \\ \hat{d}(v) & v\end{pmatrix}.$$
We check the equality entry by entry:
\begin{enumerate}
 \item[(11)] We need to check that $s\circ v\circ d_M=d_L\circ s\circ v+s\circ\hat{d}(v)$.
 But we have $s\circ\hat{d}(v)=s\circ (d_N\circ v-(-1)^{|v|}v\circ d_M)=s\circ d_N\circ v+s\circ v\circ d_M$
 and, using the fact that $s\circ d_N=-d_L\circ s$, the desired equality follows.

\item[(12)] The equality $s\circ v=s\circ v$ is clear.

\item[(21)] We need to check that $\hat{d}(v)\circ d_M=d_N\circ\hat{d}(v)$.
But this is an immediate consequence of the fact that $0=\hat{d}(\hat{d}(v))=
d_N\circ\hat{d}(v)-(-1)^{| \hat{d}(v)|}\hat{d}(v)\circ d_M=d_N\circ\hat{d}(v)-\hat{d}(v)\circ d_M$.

\item[(22)] We need to check that $\hat{d}(v)+v\circ d_{M[1]}=d_N\circ v$.
This is a direct consequence of the equality
$\hat{d}(v)=d_N\circ v-(-1)^{|v|}v\circ d_M=d_N\circ v+v\circ d_M$ and the fact that $d_{M[1]}=-d_M$.
\end{enumerate}
Once we know that $\alpha$ is a morphism in $\mathcal{C}(A)$, it is clear that
$\begin{pmatrix} 0 & 1 \end{pmatrix}\circ\alpha =\begin{pmatrix}\hat{d}(v) & v \end{pmatrix}=f$.
Therefore $C(1_M)$ is projective with respect to the semi-split exact
structure of $\mathcal{C}(A)$. It is also injective can be proved using a dual argument.
\end{proof}

\begin{Def} \label{def.contractible dg module}
A right dg $A$-module $P$ is called \emph{contractible} when it satisfies
any one of the equivalent conditions of Proposition~\ref{prop.projective-injective in C(A)}.
\end{Def}

\begin{Cor} \label{cor.homotopy category}
When $\mathcal{C}(A)$ is considered with its semi-split (Frobenius)
exact structure, its stable category $\underline{\mathcal{C}(A)}=:\mathcal{H}(A)$
is a triangulated category with arbitrary (set-indexed) coproducts, where the suspension
functor is induced by the shift functor $?[1]$ of $Dg-A$ and where the triangles are, up to isomorphism,  the images by the quotient functor $\mathcal{C}(A)\longrightarrow\mathcal{H}(A)$ of conflations in $\mathcal{C}(A)$. Moreover, $\mathcal{H}(A)$ is
equivalent to the $0$-homology category  $H^0(Dg-A)$.
\end{Cor}

\begin{proof}
It is a standard fact (see \cite[Section I.2]{H}) that the stable category $\underline{\mathcal{E}}$ of a Frobenius exact
category $\mathcal{E}$ is a triangulated category whose suspension functor is the cosyzygy functor and whose triangles are, up to isomorphism, the images by the projection functor $\mathcal{E}\longrightarrow\underline{\mathcal{E}}$ of conflations in $\mathcal{E}$.
But, for each object $M$ of $\mathcal{C}(A)$, we have a conflation
$0\rightarrow M\longrightarrow C(1_M)\longrightarrow M[1]\rightarrow 0$, where $C(1_M)$
is contractible. It follows that the shift functor $?[1]:Dg-A\longrightarrow Dg-A$ (or
$?[1]:\mathcal{C}(A)\longrightarrow\mathcal{C}(A)$) induces the suspension functor of
$\underline{\mathcal{C}(A)}=:\mathcal{H}(A)$.

For the last assertion, we need to prove that a morphism $f:M\longrightarrow N$ in $\mathcal{C}(A)$
factors through a contractible dg $A$-module if, and only if,  it is in the image of
the internal differential $d:HOM_A(M,N)\longrightarrow HOM_A(M,N)$. To avoid confusions,
we will denote by $\hat{d}$ this internal differential.  Indeed, the morphism $f$
factors through a contractible dg $A$-module if, and only if, it factors through
the canonical deflation (= admissible epimorphism)
$C(1_{N[-1]})\stackrel{\begin{pmatrix} 0 & 1 \end{pmatrix}}{\longrightarrow}N$.
This happens if, and only if, there is a morphism $\sigma :M\longrightarrow N$
of degree $-1$ in $GR-A$ such that
$\begin{pmatrix} \sigma \\ f\end{pmatrix}:M\longrightarrow C(1_{N[-1]}) \equiv N[-1]\oplus N$
is a morphism in $\mathcal{C}(A)$. This in turn is equivalent to the existence of such a $\sigma$
such that the matrix
equality $\begin{pmatrix}d_{N[-1]} & 1_N\\ 0 & d_N \end{pmatrix}\circ\begin{pmatrix} \sigma \\ f\end{pmatrix}=
\begin{pmatrix} \sigma \\ f\end{pmatrix}\circ d_M$ holds. That is, $f$ factors through
a projective object if, and only if, there is morphism $\sigma :M\longrightarrow N$ of degree
$-1$ in $GR-A$ such that $d_{N[-1]}\circ\sigma +f=\sigma\circ d_M$, which is equivalent to
saying that $f=\hat{d}(\sigma)$.
\end{proof}

\begin{Def} \label{def.homotopy category}
The category $\mathcal{H}(A)$ of last corollary is called the \emph{homotopy category of $A$}. A morphism $f:M\longrightarrow N$ in $\mathcal{C}(A)$ is called \emph{null-homotopic} when it is a $0$-boundary $f\in B^0(\text{HOM}_A(M,N))$ (i.e. $f=d_N\circ\sigma +\sigma\circ d_M$, for some $\sigma\in\text{HOM}_A^{-1}(M,N)$). This is equivalent to say that $f$ is mapped onto zero by the projection functor $\mathcal{C}(A)\longrightarrow\mathcal{H}(A)$.
\end{Def}

Note that if $f:M\longrightarrow N$ is a morphism in $\mathcal{C}(A)$ then we get induced morphisms of $K$-modules $Z^k(f):=f_{| Z^k(M)}:Z^k(M)\longrightarrow Z^k(N)$ and $B^k(f):= f_{| B^k(M)}:B^k(M)\longrightarrow B^k(N)$, for all $k\in\mathbb{Z}$. They give rise to functors $Z^k,B^k:\mathcal{C}(A)\longrightarrow\text{Mod}-K$, for all $k\in\mathbb{Z}$, and, gathering all together, to functors $Z^*,B^*:\mathcal{C}(A)\longrightarrow Gr-K$ given by $Z^*(M)=\oplus_{k\in\mathbb{Z}}Z^k(M)$ (resp. $B^*(M)=\oplus_{k\in\mathbb{Z}}B^k(M)$) and $Z^*(f)=\oplus_{k\in\mathbb{Z}}Z^k(f)$ (resp. $B^*(f)=\oplus_{k\in\mathbb{Z}}B^k(f)$). These functors are compatible with the inclusions $B^k(?)\hookrightarrow Z^k(?)$ and, hence, they give rise to functors $H^k:\mathcal{C}(A)\longrightarrow\text{Mod}-K$ ($k\in\mathbb{Z}$) and $H^*:\mathcal{C}(A)\longrightarrow Gr-K$, where $H^k(M)=Z^k(M)/B^k(M)$ and $H^*(M)=\oplus_{k\in\mathbb{Z}}H^k(M)$. We call $H^k$ the \emph{k-th homology functor} and, without mentioning the degree, we call $H^*$ the \emph{homology functor}.   If $f$ is null-homotopic, and hence $f=d_N\circ\sigma +\sigma\circ d_M$, for some morphism $\sigma:M\longrightarrow N[-1]$ in $Gr-A$, then $\text{Im}(Z^k(f))\subseteq B^k(N)$, for all $k\in\mathbb{N}$. This implies that the functor $H^k$ vanishes on null-homotopic morphisms, for all $k\in\mathbb{Z}$, which implies that we have a uniquely determined functor, still denoted and called the same,  $H^k:\mathcal{H}(A)\longrightarrow\text{Mod}-K$ such that the composition $\mathcal{C}(A)\stackrel{proj}{\longrightarrow}\mathcal{H}(A)\stackrel{H^k}{\longrightarrow}\text{Mod}-K$ is the $k$-th homology functor. We also get a corresponding  functor $H^*:\mathcal{H}(A)\longrightarrow\text{Mod}-A$.

\begin{Def} \label{def.quasi-isomorphism}
A \emph{quasi-isomorphism} of dg modules is a morphism $f:M\longrightarrow N$ in $\mathcal{C}(A)=Z^0(Dg-A)$ such that $H^*(f)$ is an isomorphism in $Gr-K$. This is equivalent to saying that its cone  $C(f)$ is an acyclic dg $A$-module (see Remark \ref{rem.cone}).
\end{Def}

As in the case of small dg categories and their dg modules, the class of quasi-isomorphisms is a
multiplicative system in $\mathcal{H}(A)$ compatible with the triangulation,
in the sense of Verdier (see  \cite[Section II.2]{Verdier}, where we refer the reader for the concepts and terminology concerning localization of triangulated categories used in this paper).

\begin{Def} \label{def.derived category}
The localization of $\mathcal{H}(A)$ with respect to the class of quasi-isomorphisms,
denoted by $\mathcal{D}(A)$, is called the \emph{derived category of $A$}. It is  
a triangulated category with arbitrary coproducts and the canonical functor
$q:\mathcal{H}(A)\longrightarrow\mathcal{D}(A)$ is a triangulated functor. The shift in $\mathcal{D}(A)$ is induced by that of $\mathcal{H}(A)$ and the triangles in $\mathcal{D}(A)$ are, up to isomorphism, the images by $q$ of triangles in $\mathcal{H}(A)$. 
\end{Def}

Note that, by the universal property of the localized category, since the functor $H^*:\mathcal{H}(A)\longrightarrow Gr-K$ takes quasi-isomorphisms to isomorphism, there is a uniquely determined functor, still denoted and called the same,  $H^*:\mathcal{D}(A)\longrightarrow Gr-K$ such that the composition $\mathcal{H}(A)\stackrel{q_A}{\longrightarrow}\mathcal{D}(A)\stackrel{H^*}{\longrightarrow} Gr-K$ is the homology functor.

\begin{Rem}
What we have done for $A$ can be done also for $A^{op}$ and for $B\otimes A^{op}$, where $B$ is another algebra with enough idempotents, obtaining the categories
$\mathcal{C}(A^{op})$, $\mathcal{H}(A^{op})$ and $\mathcal{D}(A^{op})$ (resp. $\mathcal{C}(B\otimes A^{op})$, $\mathcal{H}(B\otimes A^{op})$ and $\mathcal{D}(B\otimes A^{op})$). Due to the
equivalences of dg categories $A-Dg\cong Dg-A^{op}$ (see Proposition \ref{prop.left-versus-right}) and $A-Dg-B\cong Dg(B\otimes A^{op})$ (see Proposition \ref{prop.dg bimodules as dg modules}),
we will  look at $\mathcal{C}(A^{op})$ (resp. $\mathcal{H}(A^{op})$) and $\mathcal{C}(B\otimes A^{op})$ (resp. $\mathcal{H}(B\otimes A^{op})$)
as the $0$-cycle categories $Z^0(A-Dg)$  and $Z^0(A-Dg-B)$ (resp. $0$-homology categories $H^0(A-Dg)$ and $H^0(A-Dg-B)$).
In particular, the objects of $\mathcal{D}(A^{op})$ are considered to be left dg $A$-modules and those of $\mathcal{D}(B\otimes A^{op})$ as dg $A-B-$bimodules.
\end{Rem}

As in the case of the derived category of an abelian category (see \cite{Verdier}), we have:

\begin{Prop} \label{prop.exact sequences give triangles}
Let $A$ be a dg algebra with enough idempotents. The canonical composition  functor $\mathcal{C}(A)\stackrel{p}{\longrightarrow}\mathcal{H}(A)\stackrel{q}{\longrightarrow}\mathcal{D}(A)$ takes short exact sequences in  $\mathcal{C}(A)$ (for the abelian structure) to triangles in $\mathcal{D}(A)$. 
\end{Prop}
\begin{proof}
Let $0\rightarrow L\stackrel{u}{\longrightarrow}M\stackrel{v}{\longrightarrow}N\rightarrow 0$ be a short exact sequence in $\mathcal{C}(A)$ and fix an inflation (=admissible monomorphism with respect to the semi-split exact structure of $\mathcal{C}(A)$) $j:L\longrightarrow I$, where $I$ is a contractible dg $A$-module. If $X$ is the lower right corner of the pushout of $u$ and $j$, then we get the following commutative diagram whose rows are exact sequences:

$$\xymatrix{0 \ar[r] & L \ar[r]^{\begin{pmatrix}u\\ j \end{pmatrix}} \ar@{=}[d] & M\oplus I \ar[r] \ar[d]^{(1,0)} & X \ar[r] \ar[d]^{\epsilon}  & 0 \\ 0 \ar[r] & L \ar[r]^{u} & M \ar[r]^{v \hspace{0.3 cm}} & N \ar[r] & 0}.$$  It then follows that the right square of this diagram is bicartesian and, as a consequence, that $I=\text{Ker}[\begin{pmatrix} 1 & 0\end{pmatrix}]\cong\text{Ker}(\epsilon)$.  Since $I$ is acyclic we get that $\epsilon$ is a quasi-isomorphism and, hence, the three vertical arrows of last diagram are quasi-isomorphisms. Then the images of the 'rows' of this diagram by the canonical functor $q\circ p:\mathcal{C}(A)\longrightarrow\mathcal{D}(A)$ are isomorphic. But the upper row of the diagram is a conflation (see \cite[Proposition 2.12]{B}), whose image by $q\circ p$ 'is' then a triangle in $\mathcal{D}(A)$(see Definition \ref{def.derived category}). 
\end{proof}

Recall that if $\mathcal{D}$ is a triangulated category with (set-indexed) coproducts, then an object $C$ of $\mathcal{D}$ is called \emph{compact} when the functor $\text{Hom}_\mathcal{D}(C,?):\mathcal{D}\longrightarrow Ab$ preserves coproducts. The category $\mathcal{D}$ is said to be \emph{compactly generated} when there is a set $\mathcal{S}$ of compact objects such that $\bigcap_{n\in\mathbb{Z},S\in\mathcal{S}}\text{Ker}(\text{Hom}_\mathcal{D}(S,?[n])=0$, and $\mathcal{D}$ is said to be \emph{algebraic} when it is triangle equivalent to the stable category of some Frobenius exact category (see \cite[Section I.2]{H}).
As an immediate consequence of \cite[Theorem 4.3]{K1}, our
Theorem \ref{thm.dg modules versus dg functors} and its proof, we get:

\begin{Cor} \label{cor.Keller theorem}
For a triangulated category  $\mathcal{D}$, the following assertions are equivalent:
\begin{enumerate}
\item $\mathcal{D}$ is compactly generated and algebraic.
\item $\mathcal{D}$ is triangle equivalent to  $\mathcal{D}(A)$, for some
dg algebra with enough idempotents $A$.
\end{enumerate}
\end{Cor}

\section{Derived functors}

We call the attention of the reader on the following fact, that we shall freely use.

\begin{Rem} \label{rem.homological natural transformation}
If in Definition \ref{def.homological natural transformation} one has
$\mathcal{A}=Dg-A$ and $\mathcal{B}=Dg-B$, for some dg algebras with enough
idempotents $A$ and $B$, then the homological condition translate into the
fact that $\tau_M$ commutes with the differentials. That is,  that
$d_{G(M)}\circ\tau_M=\tau_M\circ d_{F(M)}$, for each right dg $A$-module $M$.
\end{Rem}

We start with the following  observation

\begin{Lemma} \label{lem.dg functor commutes with shift}
Let $A$ and $B$ be dg algebras with enough idempotents and $F:Dg-A\longrightarrow Dg-B$
be a dg functor. Then there is a natural isomorphism
$\rho_{F,?} :F\circ (?[1])\cong (?[1])\circ F$ which is natural on $F$. That is,
such that if $\tau :F\longrightarrow G$ is a natural transformation of dg functors,
then the following diagram in $Dg-B$ is commutative, for each right dg $A$-module $M$:
$$\xymatrix{
F(M[1])\ar[r]^{\xi_{F,M}}\ar[d]_{\tau_{M[1]}}&F(M)[1]\ar[d]^{\tau_{M}[1]}\\
G(M[1])\ar[r]^{\xi_{G,M}}&G(M)[1]
}$$
\end{Lemma}

\begin{proof}
We consider the morphism $1_M^-\in\text{HOM}_A^{-1}(M,M[1])=\text{Hom}_{Gr-A}(M,M)$
given by $1_M^-=1_M$. Then
$$F(1_M^-)\in\text{HOM}_B^{-1}(F(M),F(M[1]))=\text{HOM}_{Gr-B}(F(M),F(M[1])[-1]),$$ and hence
$F(1_M^-)[1]\in\text{Hom}_{Gr-B}(F(M)[1],F(M[1]))$.

Similarly, we have
$1_M^+\in\text{HOM}_A^{1}(M[1],M)=\text{HOM}_{Gr-A}(M[1],M[1])$ given by
$1_M^+=1_{M[1]}$, so that
$$F(1_M^+)\in\text{HOM}_B^{1}(F(M[1]),F(M))=\text{Hom}_{Gr-B}(F(M[1]),F(M)[1]).$$
We then get that $$1_{F(M[1])}=F(1_{M[1]})=F(1_M^-\circ 1_M^+)=F(1_M^-)\circ F(1_M^+).$$
By the definition of the composition of morphisms in $GR-B$, we then get that
$1_{F(M[1])}$ is equal to the composition $F(M[1])\stackrel{F(1_M^+)}{\longrightarrow}F(M)[1]\stackrel{F(1_M^-)[1]}{\longrightarrow}F(M[1])$.
On the other hand, we have $F(1_M^+)\circ (F(1_M^-)[1])=(F(1_M^+)[-1]\circ F(1_M^-))[1]$.
But due to the definition of the composition of morphisms in $GR-B$, we have
$$F(1_M^+)[-1]\circ F(1_M^-)=F(1_M^+)\circ F(1_M^-))=F(1_M^+\circ 1_M^-)=F(1_M)=1_{F(M)}.$$
We then have $$F(1_M^+)\circ (F(1_M^-)[1])=1_{F(M)}[1]=1_{F(M)[1]},$$
which shows that $F(1_M^+)$ and $(F(1_M^-)[1]$ are mutually inverse isomorphisms.

We define $\rho_{F,M}=F(1_M^+):F(M[1])\longrightarrow F(M)[1]$, for each
right dg $A$-module $M$. Note that if $\alpha :M\longrightarrow N$ is any
homogeneous morphism in $Dg-A$, then the compositions $M[1]\stackrel{\alpha[1]}{\longrightarrow}N[1]\stackrel{1_N^+}{\longrightarrow}N$ and
$M[1]\stackrel{1_M^+}{\longrightarrow}M\stackrel{\alpha}{\longrightarrow}N$ coincide in $GR-A$.
It follows that
\begin{eqnarray*}F(\alpha )[1]\circ\rho_{F,M}&=&F(\alpha )\circ F(1_M^+)=F(\alpha\circ 1_M^+)=
F(1_N^+\circ\alpha [1])\\
&=&F(1_N^+)\circ F(\alpha [1])=\rho_{F,N}\circ F(\alpha [1]),
\end{eqnarray*}
when we interpret $\alpha [1]$ as an element of $\text{HOM}_A^{-1}(M[1],N)$, using the
definition of the composition of morphisms in $GR-A$ and $GR-B$ and the functoriality
of $F$. It follows that $\rho =(\rho_{F,N})_{N\in Dg-A}$ defines a natural isomorphism
$F\circ (?[1])\cong (?[1])\circ F$.

It remains to check the commutativity of the diagram in the statement, whenever
$\tau :F\longrightarrow G$ is a natural transformation of dg functors.
But we have $\tau_M[1]\circ\rho_{F,M}=\tau_{M}[1]\circ F(1_M^+)=\tau_M\circ F(1_M^+)$,
when viewing $F(1_M^+)$ as an element of $\text{HOM}_B^{1}(M[1],M)$. The naturality of
$\tau$ then gives that
$\tau_M[1]\circ\rho_{F,M}=G(1_M^+)\circ\tau_{M[1]}=\rho_{G,M}\circ\tau_{M[1]}$, as desired.
\end{proof}

\begin{Prop} \label{prop.homotopical resolutions}
Let $A$ be a dg algebra with enough idempotents.
The canonical functor $q=q_A:\mathcal{H}(A)\longrightarrow\mathcal{D}(A)$ has
a left adjoint and a right adjoint, both of them triangulated and fully faithful.
\end{Prop}

\begin{proof}
Keller proved (see \cite[Theorems 3.1 and 3.2]{K1}) that, for each object
$M$ of $\mathcal{H}(A)$, we have quasi-isomorphisms $\pi=\pi_M:P_M\longrightarrow M$
and $\iota =\iota_M:M\longrightarrow I_M$ in $\mathcal{H}(A)$, where  $P_M$ and
$I_M$ are right dg  $A$-modules such that the functors
$\text{Hom}_{\mathcal{H}(A)}(P_M,?)$ and $\text{Hom}_{\mathcal{H}(A)}(?,I_M)$
vanish on acyclic complexes. By a standard argument, one sees that this
last property implies that the maps
$\text{Hom}_{\mathcal{H}(A)}(P_M,N)\longrightarrow\text{Hom}_{\mathcal{D}(A)}(P_M,N)$ and
$\text{Hom}_{\mathcal{H}(A)}(N,I_M)\longrightarrow\text{Hom}_{\mathcal{D}(A)}(N,I_M)$
defined by $q$ are both bijective, for any object $N$ of $\mathcal{H}(A)$.

We now define the left adjoint
$\Pi:\mathcal{D}(A)\longrightarrow\mathcal{H}(A)$ as follows.
For each right dg $A$-module $M$, we fix a quasi-isomorphism $\pi_M:P_M\longrightarrow M$
as above, and define $\Pi (M)=P_M$ on objects. If now $f:M\longrightarrow N$
is a morphism in $\mathcal{D}(A)$, then
$\pi_N^{-1}\circ f\circ\pi_M\in\text{Hom}_{\mathcal{D}(A)}(P_M,P_N)$. By the
last paragraph, we then get a unique morphism $\alpha :P_M\longrightarrow P_N$ in $\mathcal{H}(A)$
such that $q(\alpha )=\pi_N^{-1}\circ f\circ\pi_M$. We define $\Pi(f)=\alpha $.
It is routine to check that in this way we have defined a functor
$\Pi:\mathcal{D}(A)\longrightarrow\mathcal{H}(A)$. Moreover the map
$\text{Hom}_{\mathcal{H}(A)}(\Pi (M),N)\longrightarrow\text{Hom}_{\mathcal{D}(A)}(M,N)=
\text{Hom}_{\mathcal{D}(A)}(M,q(N))$
taking $\beta\rightsquigarrow q(\beta)\circ \pi_M^{-1}$
is bijective and natural on both arguments. Then $\Pi$ is left
adjoint to $q$. The co-unit of this adjunction is just
$\pi :\Pi\circ q\longrightarrow 1_{\mathcal{H}(A)}$, where $\pi_M:(\Pi\circ q)(M)=
P_M\longrightarrow M$ is the quasi-isomorphism fixed above. The unit
$\lambda :1_{\mathcal{D}(A)}\longrightarrow q\circ \Pi$ is given by
$\lambda_M=\pi_M^{-1}:M\longrightarrow (q\circ \Pi)(M)=P_M$. It follows
that $\lambda$ is a natural isomorphism, which implies that $\Pi$ is
fully faithful (see \cite[Proposition II.7.5]{HS}). On the other hand,
it is well-known that the left adjoint of a triangulated functor is also
triangulated (see \cite[Lemma 5.3.6]{Neeman}).

The existence of a right adjoint $\Upsilon:\mathcal{D}(A)\longrightarrow\mathcal{H}(A)$
acting on objects as $\Upsilon (M)=I_M$ is proved by an argument dual to the one
in the previous paragraphs.
\end{proof}

\begin{Def} \label{def.homotopical resolution}
A right dg $A$-module $P$ (resp. $I$) is called
\emph{homotopically projective (resp. injective)} if the functor
$\text{Hom}_{\mathcal{H}(A)}(P,?):\mathcal{H}(A)\longrightarrow\text{Mod}-K$ (resp.
$\text{Hom}_{\mathcal{H}(A)}(?,I)=\mathcal{H}(A)^{op}\longrightarrow\text{Mod}-K$)
vanishes on acyclic complexes. By the proof of
Proposition~\ref{prop.homotopical resolutions}, if $\Pi$ and
$\Upsilon$ are the left and right adjoints of
$q:\mathcal{H}(A)\longrightarrow\mathcal{D}(A)$, respectively,
then  the essential image $Im(\Psi)$ (resp. $Im(\Upsilon )$)
consists of homotopically projective (resp. injective) objects.
We will call $\Pi$ and $\Upsilon$ the \emph{homotopically projective
resolution functor} and \emph{homotopically injective resolution functor},
respectively. Given a right dg $A$-module $M$, a \emph{homotopically
projective resolution} (resp. \emph{homotopically injective resolution})
of $M$ will be a quasi-isomorphism $\pi :P\longrightarrow M$ (resp.
$\iota :M\longrightarrow I$), where $P$ is a homotopically projective
(resp. injective) right dg $A$-module.
\end{Def}

\begin{Rem} \label{rem.homotopically projective injective}
Note that we have $$H^k(\text{HOM}_A(M,N))=(H^0\circ (?[k])(\text{HOM}_A(M,N))=H^0(\text{HOM}_A(M,N[k]))=\text{Hom}_{\mathcal{H}(A)}(M,N[k]),$$ for all $k\in\mathbb{Z}$. Then saying
that $P$ (resp. $Y$) is a homotopically projective (resp. homotopically
injective) dg $A$-module is equivalent to saying that the dg $K$-module
$\text{HOM}_A(P,N)$ (resp. $\text{HOM}_A(N,Y)$) is acyclic whenever $N$
is an acyclic dg $A$-module.
\end{Rem}

\begin{Example} \label{ex.representables are homot.projective}
If $(e_i)_{i\in I}$ is a distinguished family of orthogonal idempotents of
$A$, then all right dg $A$-modules $e_iA$ are homotopically projective.
\end{Example}

\begin{proof}
It is a  consequence of Remark \ref{rem.projective versus representable}
and \cite[Theorem 3.1]{K1}.
\end{proof}

\medskip

Let us consider dg functors $F:\mathcal{A}\longrightarrow\mathcal{B}$ and
$G:\mathcal{B}\longrightarrow\mathcal{A}$ between dg categories. By
Examples \ref{ex.the regular dg bifunctor} and \ref{ex.tensor product of dg functors},
we then have dg functors
$\mathcal{B}(F(?),?):\mathcal{A}^{op}\otimes\mathcal{B}\longrightarrow Dg-K$ and
$\mathcal{B}(?,G(?)):\mathcal{A}^{op}\otimes\mathcal{B}\longrightarrow Dg-K$.

\begin{Def} \label{def.dg adjunction}
In the situation of last paragraph, we say that the pair $(F,G)$ is a \emph{dg adjunction}
or that \emph{$F$ is left dg adjoint to $G$} or that \emph{$G$ is right dg adjoint to $F$}
when there is a natural isomorphism of dg functors $\eta :\mathcal{B}(F(?),?)\stackrel{\cong}{\longrightarrow}\mathcal{A}(?,G(?))$. Due to
Lemma \ref{lem.dg bifunctor} and Definition \ref{def.homological natural transformation},
this means that, for each pair of objects $(A,B)\in\mathcal{A}\times\mathcal{B}$,
the map $\eta_{A,B}:\mathcal{B}(F(A),B)\longrightarrow\mathcal{A}(A,G(B))$ is an
isomorphism in $Gr-K$, natural on $A$ and $B$,  such that
$\eta_{A,B}(d_\mathcal{B}(\beta ))=d_\mathcal{A}(\eta_{A,B}(\beta ))$,
for each homogeneous element $\beta\in\mathcal{B}(F(A),B)$.
\end{Def}

\begin{Lemma} \label{lem.derived functor of dg functor}
Let $A$ and $B$ be dg algebras with enough idempotents.
If  $F:Dg-A\longrightarrow Dg-B$ (resp.  $F:(Dg-A)^{op}\longrightarrow Dg-B$)
is a dg functor, then the induced functor
$F=Z^0F:Z^0(Dg-A)\cong\mathcal{C}(A)\longrightarrow\mathcal{C}(B)=Z^0(Dg-B)$ (resp.
$F:=Z^0F:Z^0((Dg-A)^{op})\cong\mathcal{C}(A)^{op}\longrightarrow\mathcal{C}(B)=Z^0(Dg-B)$ )
preserves conflations. If moreover $F$ takes contractible dg modules to
contractible dg modules, then the induced functor
$F:=H^0F:H^0(Dg-A)\cong\mathcal{H}(A)\longrightarrow\mathcal{H}(B)=H^0(Dg-B)$ (resp.
$F:=H^0F:H^0((Dg-A)^{op})\cong\mathcal{H}(A)^{op}\longrightarrow\mathcal{H}(B)=H^0(Dg-B)$ )
is triangulated. When $F$ is part of a dg adjunction, it takes contractible
dg modules to contractible dg modules.
\end{Lemma}

\begin{proof}
We prove the covariant part of the lemma, the contravariant part being entirely
similar. Since $F$ is a dg functor the induced functor
$Z^0F:Z^0(Dg-A)=\mathcal{C}(A)\longrightarrow Z^0(Dg-B)=\mathcal{C}(B)$
preserves exact sequences which split in in the underlying graded categories.
That is, $Z^0F$ takes conflations to conflations. As a consequence, $H^0F$
takes triangles to triangles. The initial assertion then follows from
\cite[Lemma 2.27]{NS-Japan}, bearing in mind that, by
Lemma \ref{lem.dg functor commutes with shift}, we also have a natural
isomorphism $H^0F\circ (?[1])\cong (?[1])\circ H^0F$. To end the proof,
it will enough to show that if $(F,G)$ is a dg adjunction, then also
$(Z^0F,Z^0G)$ is an adjunction. Indeed, if this is proved and if $(G,F)$ is
a dg adjunction, then we will have that $(Z^0G,Z^0F)$ is an adjunction.
In any of the two situations,  \cite[Lemma 2.27]{NS-Japan} again gives that
$Z^0F$ preserves projective (=injective) objects with respect to the
semi-split exact structures, which amounts to saying that $F$ preserves
contractible dg modules.

Let $\eta :\text{HOM}_B(F(?),?)\stackrel{\cong}{\longrightarrow}\text{HOM}_A(?,G(?))$
be a graded natural isomorphism which commutes with the differentials.
Bearing in mind that,  for each $M\in Dg-A$ and $X\in Dg-B$, we have
$(Z^0F)(M)=F(M)$ and $(Z^0G)(X)=G(X)$, we then get an isomorphism of $K$-modules
$$\xymatrix{
\text{Hom}_{\mathcal{C}(B)}((Z^0F)(M),X)\ar@{=}[d]&&
\text{Hom}_{\mathcal{C}(A)}(M,(Z^0G)(X))\ar@{=}[d]\\
Z^0(\text{HOM}_B(F(M),X))\ar[rr]^{Z^0(\eta_{M,X})}&&
Z^0(\text{HOM}_A(M,G(X))),}
$$
which is natural on $M$ and $X$ since so is $\eta$.
\end{proof}

\medskip

The following functors will be very important in the sequel.

\begin{Def} \label{def.derived functors}
Let $A$ and $B$ be dg algebras with enough idempotents and let
$\Pi=\Pi_A:\mathcal{D}(A)\longrightarrow\mathcal{H}(A)$ and
$\Upsilon=\Upsilon_A :\mathcal{D}(A)\longrightarrow\mathcal{H}(A)$
be the homotopically projective and the homotopically
injective resolution functors, respectively.
\begin{enumerate}
\item  If $F:Dg-A\longrightarrow Dg-B$ is a dg functor which preserves
contractible dg modules and we also put
$H^0F=F:\mathcal{H}(A)\longrightarrow\mathcal{H}(B)$, then:
\begin{enumerate}
\item The composition $\mathbb{R}F:\mathcal{D}(A)\stackrel{\Upsilon}{\longrightarrow}\mathcal{H}(A)
    \stackrel{F}{\longrightarrow}\mathcal{H}(B)\stackrel{q_B}{\longrightarrow}\mathcal{D}(B)$
    is called the \emph{(total) right derived functor of $F$}.
\item The composition $\mathbb{L}F:\mathcal{D}(A)\stackrel{\Pi}{\longrightarrow}\mathcal{H}(A)
    \stackrel{F}{\longrightarrow}\mathcal{H}(B)\stackrel{q_B}{\longrightarrow}\mathcal{D}(B)$
    is called the \emph{(total) left derived functor of $F$}.
\end{enumerate}
\item If $F:(Dg-A)^{op}\longrightarrow Dg-B$ is a dg functor which preserves
contractible dg modules, then  the composition $\mathbb{R}F:\mathcal{D}(A)^{op}\stackrel{\Pi^o}{\longrightarrow}\mathcal{H}(A)^{op}
\stackrel{F}{\longrightarrow}\mathcal{H}(B)\stackrel{q_B}{\longrightarrow}\mathcal{D}(B)$
is called the \emph{(total) right derived functor of $F$}.
\end{enumerate}
\end{Def}

\begin{Rem} \label{rem.not left derived of contravariant}
If $F:(Dg-A)^{op}\longrightarrow Dg-B$ is a dg functor as in last definition,
we can interpret it also as dg functor $F^o:Dg-A\longrightarrow (Dg-B)^{op}$.
We then have that $\mathbb{R}F=(\mathbb{L}F^o)^o$, where $\mathbb{L}F^o$ is
the composition $\mathcal{D}(A)\stackrel{\Pi}{\longrightarrow}\mathcal{H}(A)
\stackrel{F^o}{\longrightarrow}\mathcal{H}(B)^{op}
\stackrel{q_B^o}{\longrightarrow}\mathcal{D}(B)^{op}$.
This is the reason for which we have not talked about left derived
functors of contravariant dg functors.
\end{Rem}

\begin{Rem} \label{rem.derived functor of exact dg functor}
If in any of the situations of last definition, the dg functor $F$ also
preserves acyclic dg modules, then the induced functor on the homotopy
categories $H^0F$ preserves quasi-isomorphisms. Then one gets a well-defined
unique triangulated functor $F:\mathcal{D}(A) \longrightarrow\mathcal{D}(B)$
(resp. $F:\mathcal{D}(A)^{op}\longrightarrow\mathcal{D}(B)$) such that
$q_B\circ H^0F\cong F\circ q_A$ (resp. $q_B\circ H^0F\cong F\circ q_A^o$).
It immediately follows that there are natural isomorphisms $\mathbb{L}F\cong F\cong\mathbb{R}F$.
\end{Rem}

If $F,G:Dg-A\longrightarrow Dg-B$ are dg functors and $\tau :F\longrightarrow G$
is a homological natural transformation, then, for each right dg $A$-module $M$,
we have that $\tau_M:F(M)\longrightarrow G(M)$ belongs to
$Z^0(\text{HOM}_B(F(M),G(M))=\text{Hom}_{\mathcal{C}(B)}(F(M),G(M))$.
We then get a morphism, still denoted the same,
$\tau_M:(H^0F)(M)=F(M)\longrightarrow G(M)=(H^0G)(M)$ in $\mathcal{H}(B)$. It is
seen in a straightforward way that, when $M$ varies, the $\tau_M$ give a natural
transformation $H^0F\longrightarrow H^0G$.
As a consequence we get  induced natural transformations
$$q(\tau_{\Pi_A (?)}):q_B\circ H^0F\circ\Pi=\mathbb{L}F\longrightarrow\mathbb{L}G=q_B\circ H^0F\circ\Pi$$
and $$q(\tau_{\Upsilon_A (?)}):q_B\circ H^0F\circ\Upsilon=
\mathbb{R}F\longrightarrow\mathbb{R}G=q_B\circ H^0F\circ\Upsilon.$$
An analogous fact holds
for when $F$ and $G$ are dg functors $(Dg-A)^{op}\longrightarrow Dg-B$.

\begin{Prop} \label{prop.dg transformation yields triang transformation}
Let $A$ and $B$ be dg algebras with enough idempotents. The following assertions hold:
\begin{enumerate}
\item If $F, G:Dg-A\longrightarrow Dg-B$ are  dg functors which take contractible
dg modules to contractible dg modules and $\tau :F\longrightarrow G$ is a homological
natural transformation of dg functors, then:
\begin{enumerate}
\item $q_B(\tau_{\Pi (?)}):\mathbb{L}F\longrightarrow\mathbb{L}G$ is a natural
transformation of triangulated functors $\mathcal{D}(A)\longrightarrow\mathcal{D}(B)$.
\item $q_B(\tau_{\Upsilon (?)}):\mathbb{R}F\longrightarrow\mathbb{R}G$ is a natural
transformation of triangulated functors $\mathcal{D}(A)\longrightarrow\mathcal{D}(B)$.
\end{enumerate}
\item If $F,G:(Dg-A)^{op}\longrightarrow Dg-B$ are dg functors which take contractible
dg modules to contractible dg modules and $\tau :F\longrightarrow G$ is a homological
natural transformation of dg functors, then $q_B(\tau_{\Pi_A^o ()}):\mathbb{R}F\longrightarrow\mathbb{R}G$
is natural transformation of triangulated functors $\mathcal{D}(A)^{op}\longrightarrow\mathcal{D}(B)$.
\end{enumerate}

Abusing of notation, all these natural transformations of triangulated
functors will be still denoted by $\tau$.  Moreover, if  in assertion
1.a (resp. 1.b, resp. 2),  $M$ is a right dg $A$-module such that $\tau_{\Pi_A(M)}$
(resp. $\tau_{\Upsilon_A(M)}$, resp. $\tau_{\Pi_A^o(M)}$) is a quasi-isomorphism
(e.g. an isomorphism in $\mathcal{H}(B)$ or $Dg-B$), then the evaluation of
$\tau :\mathbb{L}F\longrightarrow\mathbb{L}G$ (resp.
$\tau :\mathbb{R}F\longrightarrow\mathbb{R}G$ in 1.b and 2.b) at $M$ is an
isomorphism in $\mathcal{D}(B)$.
\end{Prop}

\begin{proof}
The proof in the three cases resemble each other very much.  We just prove 1.b.
The paragraph preceding this proposition shows that we have an induced
natural transformation
$q_B(\tau_{\Pi_A ()}):\mathbb{L}F=q_B\circ H^0F\circ\Pi_A\longrightarrow q_B\circ H^0G\circ\Pi_A=\mathbb{L}G$.
All we need to prove is that it is a natural transformation of triangulated
functors, for which it is enough to check that the induced natural transformation
$\tau:H^0F\longrightarrow H^0G$ is a natural transformation of triangulated functors.
Indeed, since $\Pi_A:\mathcal{D}(A)\longrightarrow\mathcal{H}(A)$ and
$q_B:\mathcal{H}(B)\longrightarrow\mathcal{D}(B)$ are triangulated functors,
it will follow that
$q_B(\tau_{\Pi (?)}):\mathbb{L}F=q_B\circ H^0F\circ\Pi_A\longrightarrow q_B\circ H^0G\circ\Pi_A=\mathbb{L}G$
is natural transformation of triangulated functors, as it is desired.

If $L\stackrel{\alpha}{\longrightarrow}M\stackrel{\beta}{\longrightarrow}N
\stackrel{\gamma}{\longrightarrow}L[1]$ (*) is a triangle in $\mathcal{H}(A)$,
then we may assume that  it comes from a conflation
$0\rightarrow L\stackrel{\alpha}{\longrightarrow}M\stackrel{\beta}{\longrightarrow}N\rightarrow 0$
in $\mathcal{C}(A)$, where $M$ is the cone of some morphism $\gamma [-1]:N[-1]\longrightarrow L$.
If now $\xi_F:F\circ (?[1])\cong (?[1])\circ F$ and  $\xi_G:G\circ (?[1])\cong (?[1])\circ G$
are the natural isomorphisms of Lemma \ref{lem.dg functor commutes with shift},
then the image of the triangle (*) by $H^0F$ is $$F(L)\stackrel{F(\alpha)}{\longrightarrow}F(M)\stackrel{F(\beta)}{\longrightarrow}F(N)
\stackrel{\xi_{F,L}\circ F(\gamma)}{\longrightarrow} $$ and the corresponding is true
when replacing $F$ by $G$. Due to the mentioned Lemma \ref{lem.dg functor commutes with shift},
we then have a commutative diagram in $\mathcal{H}(B)$ whose rows are triangles:
$$\xymatrix{
F(L)\ar[r]^{F(\alpha)}\ar[d]^{\tau_L}&F(M)\ar[r]^{F(\beta)}\ar[d]^{\tau_M}&
F(N)\ar[r]^{\xi_{F,L}\circ F(\gamma)}\ar[d]^{\tau_N}&F(L)[1]\ar[d]^{\tau_L[1]}\\
G(L)\ar[r]^{G(\alpha)}&G(M)\ar[r]^{G(\beta)}&G(N)\ar[r]^{\xi_{G,L}\circ G(\gamma)}&G(L)[1]
}$$
which shows that $q_B(\tau_{\Pi (?)}):H^0F\longrightarrow H^0G$ is a
natural transformation of  triangulated functors.

The last  statement is a direct consequence of the definition of the
triangulated transformation $\tau$ since the functor
$q_?:\mathcal{H}(?)\longrightarrow\mathcal{D}(?)$
takes quasi-isomorphisms to isomorphisms.
\end{proof}

\begin{Prop} \label{prop.dg adjunction - triangulated adjunction}
Let $A$ and $B$ be dg algebras with enough idempotents. The following assertions hold:
\begin{enumerate}
\item If $(F:Dg-A\longrightarrow Dg-B,G:Dg-B\longrightarrow Dg-A)$ is a dg adjunction of
dg functors, then $(\mathbb{L}F:\mathcal{D}(A)\longrightarrow\mathcal{D}(B),\mathbb{R}G:\mathcal{D}(B)
\longrightarrow\mathcal{D}(A))$ is an adjunction of triangulated functors.
\item If $(F^o:Dg-A\longrightarrow (Dg-B)^{op},G:(Dg-B)^{op}\longrightarrow Dg-A)$ is
a dg adjunction of dg functors, then $((\mathbb{R}F)^o:\mathcal{D}(A)\longrightarrow\mathcal{D}(B)^{op},\mathbb{R}G:\mathcal{D}(B)^{op}
\longrightarrow\mathcal{D}(A))$ is an adjunction of triangulated functors.
\end{enumerate}
\end{Prop}

\begin{proof}
In the proof of Lemma \ref{lem.derived functor of dg functor} we have seen that,
in the situation of assertion 1, one has that $(Z^0F:\mathcal{C}(A)
\longrightarrow\mathcal{C}(B),Z^0G:\mathcal{C}(B)
\longrightarrow\mathcal{C}(A))$ is an adjoint pair.
A similar argument  proves that, in the situation of assertion 2, one has that
$(Z^0(F^o)=(Z^0F)^o:\mathcal{C}(A)\longrightarrow\mathcal{C}(B)^{op},Z^0G:\mathcal{C}(B)^{op}
\longrightarrow\mathcal{C}(A))$ is an adjoint pair.

With the obvious adaptation, \cite[Lemma 2.27 and Proposition 2.28]{NS-Japan} and
their proofs show that assertion 1 holds. As for assertion 2, note that \cite[Lemma 2.27]{NS-Japan}
also shows that $(H^0(F^o)=(H^0F)^o:\mathcal{H}(A)\longrightarrow\mathcal{H}(B)^{op},H^0G:\mathcal{H}(B)^{op}
\longrightarrow\mathcal{H}(A))$ is an adjoint pair of triangulated functors.
Moreover, the adjoint pair $(\Pi_B:\mathcal{D}(B)\longrightarrow\mathcal{H}(B),q_B:\mathcal{H}(B)
\longrightarrow\mathcal{D}(B))$ implies that the pair $(q_B^o:\mathcal{H}(B)^{op}\longrightarrow\mathcal{D}(B)^{op},\Pi_B^o:\mathcal{D}(B)^{op}
\longrightarrow\mathcal{H}(B)^{op})$ is also an adjoint pair. It then follows that the composition
$(\mathbb{R}F)^o:\mathcal{D}(A)\stackrel{\Pi_A}{\longrightarrow}\mathcal{H}(A)
\stackrel{(H^0F)^o}{\longrightarrow}\mathcal{H}(B)^{op}
\stackrel{q_B^o}{\longrightarrow}\mathcal{D}(B)^{op}$ (see
Remark \ref{rem.not left derived of contravariant}) is left adjoint to the
composition $\mathcal{D}(A)\stackrel{q_A}{\longleftarrow}\mathcal{H}(A)
\stackrel{H^0G}{\longleftarrow}\mathcal{H}(B)^{op}
\stackrel{\Pi_B^o}{\longleftarrow}\mathcal{D}(B)^{op}$, which is precisely $\mathbb{R}G$.
\end{proof}

\begin{Prop} \label{prop.derived functor of composition}
Let $A$, $B$ and $C$ be dg algebras with enough idempotents and denote by
$\Pi_?:\mathcal{D}(?)\longrightarrow\mathcal{H}(?)$ and
$\Upsilon_?:\mathcal{D}(?)\longrightarrow\mathcal{H}(?)$ the
homotopically projective and homotopically injective resolution
functors, for $?=A,B,C$. Suppose that all the dg functors appearing
below preserve contractible dg modules. The following assertions hold:
\begin{enumerate}
\item Let $G:Dg-A\longrightarrow Dg-B$ and $F:Dg-B\longrightarrow Dg-C$
be dg functors. Then:
\begin{enumerate}
\item  There is a canonical natural transformation of triangulated
functors $\rho :\mathbb{R}(F\circ G)\longrightarrow\mathbb{R}F\circ\mathbb{R}G$.
When $M$ is a right dg $A$-module such that $G(\Upsilon_A(M))$ is homotopically
injective, then $\rho_M$ is an isomorphism.
\item There is a canonical natural transformation of triangulated functors
$\sigma :\mathbb{L}F\circ\mathbb{L}G\longrightarrow\mathbb{L}(F\circ G)$.
When $M$ is a right dg $A$-module such that $G (\Pi_A (M))$ is homotopically
projective, then $\sigma_M$ is an isomorphism.
\end{enumerate}
\item If $G:(Dg-A)^{op}\longrightarrow Dg-B$ and
$F:(Dg-B)^{op}\longrightarrow Dg-C$ are dg functors, then there is a
canonical natural transformation
$\tau:\mathbb{L}(F\circ G^o)\longrightarrow\mathbb{R}F\circ (\mathbb{R}G)^o$
of triangulated functors $\mathcal{D}(A)\longrightarrow\mathcal{D}(C)$.
When $M$ is a right dg $A$-module such that $G(\Pi_A(M))$ is homotopically
projective, then $\tau_M$ is an isomorphism.
\item If $G:(Dg-A)^{op}\longrightarrow Dg-B$ and $F:Dg-B\longrightarrow Dg-C$
are dg functors, then there is a canonical natural transformation of
triangulated functors $\omega:\mathbb{L}F\circ\mathbb{R}G\longrightarrow \mathbb{R}(F\circ G)$.
When $M$ is a right dg $A$-module such that $G(\Pi_A(M))$ is homotopically projective,
$\omega_M$ is an isomorphism.
\item If $G:Dg-A\longrightarrow Dg-B$ and $F:(Dg-B)^{op}\longrightarrow Dg-C$
are dg functors, then there is a canonical natural transformation of triangulated
functors $\theta:\mathbb{R}(F\circ G^o)\longrightarrow  \mathbb{R}F\circ(\mathbb{L}G)^o$.
When $M$ is a right dg $A$-module such that $G(\Pi_A(M))$ is homotopically
projective, $\theta_M$ is an isomorphism.
\end{enumerate}
\end{Prop}

\begin{proof}
The arguments for the proofs are all very much alike and rely on the explicit
definition of right and left derived functors in each case. We just provide the
proof of assertions 1.a and 2, leaving the rest as an exercise to the reader.

1.a) We consider the unit
$\lambda:1_{\mathcal{H}(B)}\longrightarrow\Upsilon_B\circ q_B$ of the
adjunction $(q_B,\Upsilon_B)$. We then get a canonical natural transformation
of triangulated functors
\begin{center}
$\rho:=(q_C\circ F)(\lambda_{(G\circ\Upsilon) (?)}):\mathbb{R}(F\circ G)=
q_C\circ F\circ G\circ\Upsilon_A=q_C\circ F\circ 1_{\mathcal{H}(B)}
\circ G\circ\Upsilon_A\longrightarrow q_C\circ F\circ
\Upsilon_B\circ q_B\circ G\circ\Upsilon_A=\mathbb{R}F\circ\mathbb{R}G$,
\end{center}
where $F=H^0F:\mathcal{H}(B)\longrightarrow\mathcal{H}(C)$ and
$G=H^0G:\mathcal{H}(A)\longrightarrow\mathcal{H}(B)$.
If now $G(\Upsilon_A)(M)$ is homotopically injective, then
$\lambda_{(G\circ\Upsilon_A)(M)}:(G\circ\Upsilon_A)(M)
\stackrel{\cong}{\longrightarrow}(\Upsilon_B\circ q_B\circ G\circ\Upsilon_A)(M)$
is an isomorphism, which implies that $\rho_M=(q_C\circ F)(\lambda_{(G\circ\Upsilon_A)(M)})$
is also an isomorphism.

2) The adjunction $(q_B^o:\mathcal{H}(B)^{op}\longrightarrow
\mathcal{D}(B)^{op},\Pi_B^o:\mathcal{D}(B)^{op}\longrightarrow
\mathcal{H}(B)^{op})$ yields a unit
$\mu:1_{\mathcal{H}(B)^{op}}\longrightarrow\Pi_B^o\circ q_B^o$.
Then we get a natural transformation of triangulated functors

\begin{center}
$\sigma :=(q_C\circ F)(\mu_{(G^o\circ\Pi_A)(?)}):\mathbb{L}(F\circ G^o)=
q_C\circ F\circ G^o\circ\Pi_A=
q_C\circ F\circ 1_{\mathcal{H}(B)^{op}}\circ G^o\circ\Pi_A\longrightarrow q_C\circ
F\circ \Pi_B^o\circ q_B^o\circ G^o\circ\Pi_A=\mathbb{R}F\circ(\mathbb{R}G)^o$.
\end{center}
If now $G(\Pi_A(M))$ is homotopically projective, then
$\mu_{(G^o\circ\Pi_A)(M)}:(G^o\circ\Pi_A)(M)\longrightarrow
(\Pi_B^o\circ q_B^o\circ G^o\circ\Pi_A)(M)$ is an isomorphism in
$\mathcal{H}(B)^{op}$, which implies that
$\sigma_M=(q_C\circ F)(\mu_{(G^o\circ\Pi_A)(M)})$ is an isomorphism.
\end{proof}

\medskip

Suppose now that $A,B,C$ are dg algebras with enough idempotents and that
$F:(Dg-A)\otimes (Dg-C)\longrightarrow Dg-B$ is a dg functor. We then have
induced functors
$$Z^0F:Z^0((Dg-A)\otimes (Dg-C))\longrightarrow Z^0(Dg-B)=\mathcal{C}(B)$$
and  $$H^0F:H^0((Dg-A)\otimes (Dg-C))\longrightarrow H^0(Dg-B)=\mathcal{H}(B).$$
On the other hand, the objects of $\mathcal{C}(A)\otimes\mathcal{C}(B)$
are those of $Z^0((Dg-A)\otimes (Dg-C))$ (i.e. those of $(Dg-A)\otimes (Dg-C)$).
But if $f$ is morphism in $\mathcal{C}(A)$ and $g$ is a morphism in $\mathcal{C}(C)$,
then, viewed as a morphism in $(Dg-A)\otimes (Dg-C)$, we have that
$f\otimes g$ is a $0$-cycle and, hence, a morphism of $Z^0((Dg-A)\otimes (Dg-C))$.
Indeed we have
$$d(f\otimes g)=d(f)\otimes g+(-1)^{|f|}f\otimes d(g)=0,$$
because $f$ and $g$ are morphisms in $Z^0(Dg-A)=\mathcal{C}(A)$ and
$Z^0(Dg-C)=\mathcal{C}(C)$, respectively. The assignments
$(M,X)\rightsquigarrow (M,X)$ and $f\otimes g\rightsquigarrow f\otimes g$ give a functor
$j:\mathcal{C}(A)\otimes\mathcal{C}(C)\longrightarrow Z^0((Dg-A)\otimes (Dg-C))$
and a composition
$$\mathcal{C}(A)\otimes\mathcal{C}(C)\stackrel{j}{\longrightarrow}Z^0((Dg-A)\otimes (Dg-C))\stackrel{Z^0F}{\longrightarrow}\mathcal{C}(B).$$
Abusing the notation, we still denote by $Z^0F$ this composition functor.
Considering now $f$ and $g$ as above, suppose that either $f$ or $g$ is
null-homotopic. We claim that $j(f\otimes g)=f\otimes g$ is a $0$-boundary
of $(Dg-A)\otimes (Dg-C)$. Indeed if, say, $g=d(g')$ then
$$d(f\otimes g')=d(f)\otimes g'+(-1)^{|f|}f\otimes d(g')=f\otimes g$$
since
$d(f)=0$ and $|f|=0$. A similar argument applies if we assume $f=d(f')$.
This means that we have an induced functor
$\underline{j}:\mathcal{H}(A)\otimes\mathcal{H}(C)\longrightarrow H^0((Dg-A)\otimes (Dg-C))$
and a corresponding composition
$$\mathcal{H}(A)\otimes\mathcal{H}(C)\stackrel{\underline{j}}{\longrightarrow}
H^0((Dg-A)\otimes (Dg-C))\stackrel{H^0F}{\longrightarrow}\mathcal{H}(B),$$
which we shall still denote by $H^0F$.

A procedure similar to the one depicted in the previous paragraph can be
undertaken with a dg functor $F:(Dg-A)^{op}\otimes (Dg-C)\longrightarrow Dg-B$,
getting then functors
$Z^0F:\mathcal{C}(A)^{op}\otimes\mathcal{C}(C)\longrightarrow\mathcal{C}(B)$ and
$H^0F:\mathcal{H}(A)^{op}\otimes\mathcal{H}(C)\longrightarrow\mathcal{H}(B)$.

\begin{Prop} \label{prop.Z0 and H0 for dg bifunctor}
Let $A,B,C$ be dg algebras with enough idempotents and let
$F:(Dg-A)\otimes (Dg-C)\longrightarrow Dg-B$ (resp.
$F:(Dg-A)^{op}\otimes (Dg-C)\longrightarrow Dg-B$) be a dg functor. The following assertions hold:
\begin{enumerate}
\item The functor $F=Z^0F:\mathcal{C}(A)\otimes\mathcal{C}(C)\longrightarrow\mathcal{C}(B)$ (resp.
$Z^0F:\mathcal{C}(A)^{op}\otimes\mathcal{C}(C)\longrightarrow\mathcal{C}(B)$) preserves conflations
on each variable.
\item If $F(P,X)$ and $F(M,Q)$ are contractible dg $B$-modules whenever $P$ and $Q$ are a
contractible dg $A$-module and a contractible dg $C$-module, respectively, then the functor
$F=H^0F:\mathcal{H}(A)\otimes\mathcal{H}(C)\longrightarrow\mathcal{H}(B)$ (resp.
$F=H^0F:\mathcal{H}(A)^{op}\otimes\mathcal{H}(C)\longrightarrow\mathcal{H}(B)$)
is triangulated on both variables.
\end{enumerate}
\end{Prop}

\begin{proof}
Assertion 1 is a direct consequence of Lemma \ref{lem.derived functor of dg functor}
bearing in mind Lemma \ref{lem.dg bifunctor}. Moreover, if $M$ is fixed, then the
dg functor $F_M=F(M,?):(Dg-C)\longrightarrow (Dg-B)$ takes contractible dg modules to
contractible dg modules, which implies by Lemma \ref{lem.derived functor of dg functor} that
$H^0F_M:\mathcal{H}(C)\longrightarrow\mathcal{H}(B)$ is a triangulated functor. But we
clearly have $H^0F_M=H^0F(M,?)$, which says that $F=H^0F$ is triangulated on the second
variable. A symmetric argument proves that it is triangulated on the first variable.
\end{proof}

\medskip

Our next goal is to see that, when a dg functor is part of a dg `bifunctor' and certain
conditions are satisfied, also its derived functor is part of a bifunctor which is
triangulated on both variables.

\begin{Def} \label{def.derived functor of a dg bifunctor}
Let $A,B,C$ be dg algebras with enough idempotents.
\begin{enumerate}
 \item If  $F:(Dg-A)\otimes (Dg-C)\longrightarrow Dg-B$ is a dg functor which preserves
 contractible dg modules in each variable, then we put
\begin{enumerate}
\item $\mathbb{L}F:\mathcal{D}(A)\otimes\mathcal{D}(C)
    \stackrel{\Pi_A\otimes\Pi_C}{\longrightarrow}\mathcal{H}(A)\otimes\mathcal{H}(C)
    \stackrel{H^0F}{\longrightarrow}\mathcal{H}(B)\stackrel{q_B}{\longrightarrow}\mathcal{D}(B)$.
\item $\mathbb{R}F:\mathcal{D}(A)\otimes\mathcal{D}(C)
    \stackrel{\Upsilon_A\otimes\Upsilon_C}{\longrightarrow}\mathcal{H}(A)\otimes\mathcal{H}(C)
    \stackrel{H^0F}{\longrightarrow}\mathcal{H}(B)\stackrel{q_B}{\longrightarrow}\mathcal{D}(B)$
\end{enumerate}
\item If $F:(Dg-A)^{op}\otimes (Dg-C)\longrightarrow Dg-B$ be a dg functor which preserves
contractible dg modules on each variable, then we put $\mathbb{R}F:\mathcal{D}(A)^{op}\otimes\mathcal{D}(C)
\stackrel{\Pi_A^o\otimes\Upsilon_C}{\longrightarrow}\mathcal{H}(A)\otimes\mathcal{H}(C)
\stackrel{H^0F}{\longrightarrow}\mathcal{H}(B)\stackrel{q_B}{\longrightarrow}\mathcal{D}(B)$.
\end{enumerate}
By their definition all these functors are triangulated in each variable.
\end{Def}

For each dg functor $F:(Dg-A)\otimes (Dg-C)\longrightarrow Dg-B$ as in last definition,
fixing an object $M$ in $Dg-A$ and $X$ in $Dg-C$, we get dg functors
$F_M=F(M,?):Dg-C\longrightarrow Dg-B$ and $F^X=F(?,X):Dg-A\longrightarrow Dg-B$.
It is natural to ask whether we have natural isomorphisms $\mathbb{L}F(M,?)\cong\mathbb{L}F_M$ and
$\mathbb{L}F(?,X)\cong\mathbb{L}F^X$, and similarly for the right derived versions.
For this, we have the following criterion:

\begin{Prop} \label{prop.triangulated part of bitriangulated}
Let $A,B,C$ be dg algebras with enough idempotents. The following assertions hold:

\begin{enumerate}
\item Let $F:(Dg-A)\otimes (Dg-C)\longrightarrow Dg-B$  be a dg functor. Then
\begin{enumerate}
\item if $F(?,Q):Dg-A\longrightarrow Dg-B$ preserves acyclic dg modules whenever $Q$ is
homotopically projective (resp. homotopically injective), then there is a natural
isomorphism of triangulated functors
$\mathbb{L}F(M,?)\cong \mathbb{L}F_M:\mathcal{D}(C)\longrightarrow\mathcal{D}(B)$
(resp. $\mathbb{R}F(M,?)\cong \mathbb{R}F_M:\mathcal{D}(C)\longrightarrow\mathcal{D}(B)$),
for each right dg $A$-module $M$,
\item if $F(P,?):Dg-A\longrightarrow Dg-B$ preserves acyclic dg modules whenever $P$ is
homotopically projective (resp. homotopically injective), then there is a natural
isomorphism of triangulated functors
$\mathbb{L}F(?,X)\cong\mathbb{L}F^X:\mathcal{D}(A)\longrightarrow\mathcal{D}(B)$,
for each right dg $C$-module $X$.
\end{enumerate}
\item Let $F:(Dg-A)^{op}\otimes (Dg-C)\longrightarrow Dg-B$  be a dg functor. Then
\begin{enumerate}
\item if $F(?,Q):(Dg-A)^{op}\longrightarrow Dg-B$ preserves acyclic dg modules
whenever $Q$ is homotopically injective, then there is a natural isomorphism of
triangulated functors $\mathbb{R}F(M,?)\cong\mathbb{R}F_M:\mathcal{D}(C)\longrightarrow\mathcal{D}(B)$,
for each right dg $A$-module $M$.
\item if $F(P,?):Dg-C\longrightarrow Dg-B$ preserves acyclic dg modules whenever $P$ is
homotopically projective, then there is a natural isomorphism of triangulated functors
$\mathbb{R}F(?,X)\cong\mathbb{R}F^X:\mathcal{D}(A)^{op}\longrightarrow\mathcal{D}(B)$
for each right dg $C$-module $X$.
\end{enumerate}
\end{enumerate}
\end{Prop}

\begin{proof}
We will prove  1.b and 2.a, and leave to the reader the other ones whose
proof follows entirely similar patterns. For 1.b, note that the action of
$\mathbb{L}F(?,X)$ and $\mathbb{L}F^X$ on objects is given by
$$\mathbb{L}F(?,X)(M)=\mathbb{L}F(M,X)=F(\Pi_A(M),\Pi_C(X))$$ and
$$\mathbb{L}F^X(M)=F^X(\Pi_A(M))=F(\Pi_A(M),X).$$
Moreover if $f:M\longrightarrow N$
is a morphism in $\mathcal{D}(A)$, then
$$\mathbb{L}F(?,X)(f)=\mathbb{L}F(f,1_X)=F(\Pi_A (f),\Pi_C(1_X))=F(\Pi_A(f),1_{\Pi_C(X)}).$$
It follows that we can identify $\mathbb{L}F(?,X)=\mathbb{L}F^{\Pi_C(X)}$, where
$F^{\Pi_C(X)}=F(?,\Pi_C(X)):Dg-A\longrightarrow Dg-B$ is the `left part' of $F$
when the fixed second variable is $\Pi_C(X)$. We fix now the homotopically
projective resolution map $\pi :\Pi_C(X)\longrightarrow X$. Note that $\pi$ is a
morphism in $\mathcal{H}(C)$, and hence the image of a morphism in $\mathcal{C}(C)$
by the canonical functor $\mathcal{C}(C)\longrightarrow\mathcal{H}(C)$. Fixing a
lift, we can think of $\pi$ as a morphism in $\mathcal{C}(C)$. We claim that
$\pi_*:F^{\Pi_C(X)}=F(?,\Pi_C(X))\longrightarrow F(?,X)=F^X$ is a homological
natural transformation of dg functors $Dg-A\longrightarrow Dg-B$. Indeed note that,
for a fixed $M$ in $Dg-A$, we have $(\pi_*)_M:F(M,\Pi_C(X))\longrightarrow F(M,X)$
is the morphism $(\pi_*)_M=F(1_M,\pi)$.  If now $f:M\longrightarrow N$ is a homogeneous
morphism in $Dg-A$, then
$$F^{\Pi_C(X)}(f)=F(f,1_{\Pi_C(X)})=F(M,\Pi_C(X))\longrightarrow
F(N,\Pi_C(X))$$ while $$F^X(f)=F(f,1_X):F(M,X)\longrightarrow F(N,X).$$
We then have an equality
\begin{eqnarray*}
F^X(f)\circ (\pi_*)_M&=&F(f,1_X)\circ F(1_M,\pi)\\
&=&F(f,\pi)\\
&=&(-1)^{|f| |\pi|}F(1_N,s)\circ F(f,1_{\Pi_C(X)})\\
&=&(\pi_*)_N\circ F^{\Pi_C(X)}(f),
\end{eqnarray*}
using Lemma \ref{lem.dg bifunctor} and the fact that $|\pi|=0$. It follows that
$\pi_*$ is a natural transformation of $K$-linear graded functors. In order to see
that it is homological it remains to check that
$(\pi_*)_M=F(1_M\otimes\pi)\in Z^0(\text{HOM}_B(F(M,\Pi_C(X))),F(M,X))$,
for all $M$ in $Dg-A$. But this is clear since $F(1_M\otimes\pi)$ is the image
of $1_M\otimes\pi$ by the functor $Z^0F:\mathcal{C}(A)\otimes\mathcal{C}(C)\longrightarrow\mathcal{C}(B)$
(see Proposition \ref{prop.Z0 and H0 for dg bifunctor}).

Once we know that $\pi_*:F^{\Pi_C(X)}\longrightarrow F^X$ is a homological natural
transformation of dg functors, Proposition \ref{prop.dg transformation yields triang transformation}
says that we have an induced natural transformation
$\pi_*:\mathbb{L}F^{\Pi_C(X)}\longrightarrow\mathbb{L}F^X$ of triangulated functors
$\mathcal{D}(A)\longrightarrow\mathcal{D}(B)$. But, when evaluating at
an object $M$ of $\mathcal{D}(A)$, we have that
$$(\pi_*)_M:\mathbb{L}F^{\Pi_C(X)}(M)=F(\Pi_A(M),\Pi_C(X))\longrightarrow F(\Pi_A(M),X)=\mathbb{L}F^X(M)$$
is an isomorphism in $\mathcal{D}(B)$. Indeed, by hypothesis
$F(\Pi_A(M),?):Dg-C\longrightarrow Dg-B$ preserves acyclic dg modules,
which implies that the induced triangulated functor $F(\Pi_A(M),?):H^0(Dg-C)=\mathcal{H}(C)\longrightarrow\mathcal{H}(B)=H^0(Dg-B)$
preserves quasi-isomorphisms. It follows that
$(\pi_*)_M=F(1_{\Pi_A(M)},\pi):F(\Pi_A(M),\Pi_C(X))\longrightarrow F(\Pi_A(M),X)$
is a quasi-isomorphism in $\mathcal{H}(B)$, which implies that (after applying
$q_B:\mathcal{H}(B)\longrightarrow\mathcal{D}(B)$) it becomes an isomorphism in
$\mathcal{D}(B)$. Then $(\pi_*):\mathbb{L}F^{\Pi_C(X)}\longrightarrow\mathbb{L}F^X$
is a natural transformation of triangulated functor which is pointwise an isomorphism.
Therefore it is a natural isomorphism. In particular
$\mathbb{L}F(?,X)=\mathbb{L}F^{\Pi_C(X)}$ is naturally isomorphic
to $\mathbb{L}F^X$ as triangulated functors.

\bigskip

The proof of 2.a follows an entirely similar pattern. We outline the argument,
leaving the details to the reader. We have that
$\mathbb{R}F(M,?)(X)=F(\Pi_A(M),\Upsilon_A (X))$ and
$\mathbb{R}F_M(X)=F_M(\Upsilon (X))=F(M,\Upsilon (X))$. We can then identify
$\mathbb{R}F(M,?)=\mathbb{R}F_{\Pi_A(M)}$, where $F_{\Pi_A(M)}=F(\Pi_A(M),?):Dg-C\longrightarrow Dg-B$.
If now $\pi :\Pi (M)\longrightarrow M$ is homotopically projective resolution,
which we view as a morphism in $\mathcal{C}(A)$, then
$\pi^* :F_M=F(M,?)\longrightarrow F_{\Pi_A(M)}=F(\Pi_A(M),?)$ is a homological
natural transformation of dg functors $Dg-C\longrightarrow Dg-B$. The associated
natural transformation of triangulated functor $\pi^*:\mathbb{R}F_M\longrightarrow\mathbb{R}F_{\Pi_A(M)}$,
when evaluated at an object $X$ of $\mathcal{D}(C)$, is
$(\pi^*)_X=F(\pi, 1_{\Upsilon (X)}):\mathbb{R}F_M(X)=F(M,\Upsilon (X))\longrightarrow F(\Pi_A(M),\Upsilon (X))$.
This is a quasi-isomorphism of dg $B$-modules, and hence an isomorphism in
$\mathcal{D}(B)$, for each $X\in\mathcal{D}(C)$. It follows that $\pi^*$ is
a natural isomorphism
$\mathbb{R}F_M\stackrel{\cong}{\longrightarrow}\mathbb{R}F_{\Pi_A(M)}=\mathbb{R}F(M,?)$.
\end{proof}

\section{The classical dg bifunctors} \label{section.contravariant HOM}

All throughout this section, we fix  dg algebras with enough idempotents $A$,
$B$ and $C$ and fix a distinguished family of orthogonal idempotents $(e_i)_{i\in I}$ in $B$,
$(\epsilon_j)_{j\in J}$ in $A$ and $(\nu_k)_{k\in\mathcal{K}}$ in $C$,
all of  which are homogeneous of degree zero and killed by the differential.
If $M$ a dg  $C-A-$bimodule and $X$ is a dg $B-A-$bimodule,
the space of morphisms $\text{HOM}_A(M,X)$ in $Dg-A$ has a canonical
structure of non-unitary graded $B-C-$bimodule, where the multiplication
map is identified by the rule  $(bfc)(m)=bf(cm)$, for all homogeneous
elements $b\in B$, $f\in\text{HOM}_A(M,X)$, $c\in C$ and $m\in M$.
To avoid the non-unitary problem, we consider the largest unitary graded
$B-C-$sub-bimodule
$$\overline{HOM}_A(M,X)=B\text{HOM}_A(M,X)C$$
of $\text{HOM}_A(M,X)$. Note that, expressed in terms of the distinguished
families of orthogonal idempotents,  $\overline{HOM}_A(M,X)$ consists
of those $f\in\text{HOM}_A(M,X)$ such that
$\text{Im(f)}\subseteq\oplus_{i\in I'}e_iX$, for some finite subset
$I'\subseteq I$, and $f(\nu_kM)=0$, for all but finitely many
$k\in\mathcal{K}$. We can say much more about the just defined graded $B-C-$bimodule

\begin{Lemma} \label{lem.B-module structure of M-upstar}
The differential $d:\text{HOM}_A(M,X)\longrightarrow\text{HOM}_A(M,X)$ satisfies Leibniz equality $$d(bfc)=d_B(b)fc+(-1)^{|b|}bd(f)c+(-1)^{|b|+|f|}bfd_C(c),$$ for all homogeneous elements $b\in B$, $f\in\text{HOM}_A(M,X)$ and $c\in C$. Moreover, it satisfies that $d(\overline{HOM}_A(M,X))\subseteq\overline{HOM}_A(M,X)$, so that, endowed  with the restricted differential, $\overline{HOM}_A(M,X)$
becomes a (unitary!) dg $B-C-$bimodule.
\end{Lemma}

\begin{proof}
We let act both members of the desired Leiniz equality on a homogeneous
element $m\in M$. We then have
\begin{eqnarray*}
d(bfc)(m)&=&[d_X\circ (bfc)-(-1)^{|bfc|}(bfc)\circ d_M](m)\\
&=&d_X(bf(cm))-(-1)^{|b|+|f|+|c|}bf(cd_M(m))\\
&=&d_B(b)f(cm)+(-1)^{|b|}bd_X(f(cm))-(-1)^{|b|+|f|+|c|}bf(cd_M(m))\\
&=&d_B(b)f(cm)+(-1)^{|b|}b(d_X\circ f-(-1)^{|f|}f\circ d_M)(cm)\\
&&+(-1)^{|b|+|f|}bf(d_M(cm)-(-1)^{|c|}cd_M(m))\\
&=&d_B(b)f(cm)+(-1)^{|b|}bd(f)(cm)+(-1)^{|b|+|f|}bf(d_C(c)m)\\
&=&[d_B(b)fc+(-1)^{|b|}bd(f)c+(-1)^{|b|+|f|}bfd_C(c)](m).
\end{eqnarray*}

To prove the last statement, take a
homogeneous element $f\in\overline{HOM}_A(M,X)$. We have
$d(f)=d_X\circ f-(-1)^{|f|}f\circ d_M$, which implies that
$\text{Im}(d(f))\subseteq d_X(Im(f))+Im(f)$. But if $I'\subseteq I$
is any finite subset such that $\text{Im}(f)\subseteq\oplus_{i\in I'}e_iX$,
then $\text{Im}(d(f))\subseteq\oplus_{i\in I'}e_iX$. This is because
$d_X(e_iX)\subseteq e_iX$ since $d_B(e_i)=0$ for all $i\in I'$.
By analogous reason, we have that $d_M(\nu_kM)\subseteq\nu_kM$.
This implies that if $\mathcal{K}'\subseteq\mathcal{K}$ is any
finite subset such that $f(\nu_kM)=0$, for all
$k\in\mathcal{K}\setminus\mathcal{K}'$, then we also have
$d(f)(\nu_kM)=0$, for all $k\in\mathcal{K}\setminus\mathcal{K}'$.
As a consequence, we get that $d(f)\in\overline{HOM}_A(M,X)$.
 \end{proof}

\medskip


We can now prove

\begin{Prop} \label{prop.bifunctor}
The assignment $(M,X)\rightsquigarrow\overline{HOM}_A(M,X)$
is the definition on objects of a dg functor
$\overline{HOM}_A(?,?):(C-Dg-A)^{op}\otimes (B-Dg-A)\longrightarrow B-Dg-C$.
\end{Prop}

\begin{proof}
We have obvious restriction of scalars functors
$\rho :C-Dg-A\longrightarrow Dg-A$ and $\rho':B-Dg-A\longrightarrow Dg-A$,
both of which are clearly dg functors since when $M$ and $N$ are dg
$C-A-$bimodules,  the differential $d:\text{HOM}_R(M,N)\longrightarrow\text{HOM}_R(M,N)$ is `the same'
when taking $R=A\otimes C^{op}$ or when taking $R=A$. On the other hand,
by Example \ref{ex.the regular dg bifunctor}, we have a canonical
dg functor $\text{HOM}_A(?,?):(Dg-A)^{op}\otimes (Dg-A)\longrightarrow Dg-K$.
We then get an induced dg functor
$$\text{HOM}_A(?,?):(C-Dg-A)^{op}\otimes (B-Dg-A)
\stackrel{\rho^o\otimes\rho'}{\longrightarrow}(Dg-A)^{op}\otimes (Dg-A)
\stackrel{HOM_A(?,?)}{\longrightarrow}Dg-K.$$
Recall that  if  $\alpha :N\longrightarrow M$ and  $\varphi :X\longrightarrow Y$
are homogeneous morphisms in $C-Dg-A$ and $B-Dg-A$, respectively, then
$\text{HOM}_A(\alpha^o\otimes\varphi ):\text{HOM}_A(M,X)\longrightarrow\text{HOM}_A(N,Y)$
is the map defined by $\text{HOM}_A(\alpha^o\otimes\varphi )(f)=
(-1)^{(|\varphi|+|f|)|\alpha|}\varphi\circ f\circ\alpha$.

By definition, we have that $\overline{HOM}_A(M,X)$ is a dg
$B-C-$subbimodule of $\text{HOM}_A(M,X)$ with the restricted differential.
In order to have an induced graded  functor
$\overline{HOM}_A(?,?):(C-Dg-A)^{op}\otimes (B-Dg- A)\longrightarrow B-Dg-C$,
using Lemma \ref{lem.dg bifunctor}, it is enough to prove the following two conditions:
\begin{enumerate}
\item[a)] For $\varphi$ as above and each dg $C-A-$bimodule $M$, the map
$\varphi_*:=\text{HOM}_A(1_M^o,\varphi):\text{HOM}_A(M,X)\longrightarrow\text{HOM}_A(M,Y)$
is a morphisms of non-unitary $B-C-$bimodules such that
$\varphi_*(\overline{HOM}_A(M,X))\subseteq\overline{HOM}_A(M,Y)$.
\item[b)] For $\alpha$ as above and each dg $B-A-$bimodule $X$, the map
$\alpha^*:=\text{HOM}_A(\alpha^0,1_X):\text{HOM}_A(M,X)\longrightarrow\text{HOM}_A(N,X)$
is a morphism of non-unitary graded $B-C-$bimodules such  that
$\alpha^*(\overline{HOM}_A(M,X))\subseteq\overline{HOM}_A(N,X)$.
\end{enumerate}
For condition a), we first check that $\varphi_*$ is a morphism
(of degree $|\varphi|$) of nonunitary graded left $B$-modules.
According to the comments after Lemma \ref{lem.left versus right dg module},
although applied to non-unitary graded left $B$-modules,  we need to check that
$\varphi_*(bf)=(-1)^{|\varphi||b|}b\varphi_* (f)$ or, equivalently, that
$\varphi\circ (bf)=(-1)^{|\varphi||b|}b(\varphi\circ f)$, for any homogeneous
element $f\in\text{HOM}_A(M,X)$. By  letting act both members of the desired equality
on a homogeneous element $m\in M$, we get that
$$[\varphi\circ (bf)](m)=\varphi (bf(m))=(-1)^{|\varphi | |b|}b\varphi (f(m))=
[(-1)^{|\varphi||b|}b(\varphi\circ f)](m),$$
bearing in mind that $\varphi$ is a morphism of graded $B-A-$bimodules and,
hence, also a morphism of graded left $B$-modules. On the other hand,
if $c\in C$ is a homogeneous element, we have that
$\varphi_* (fc)=\varphi\circ (fc)$ while $\varphi_*(f)c=(\varphi\circ f)c$.
Both maps take $m\rightsquigarrow (\varphi\circ f)(cm)$, for each $m\in M$.
Then $\varphi_*$ is a morphism of nonunitary graded $B-C-$bimodules.
Moreover, if $f\in\overline{HOM}_A(M,X)$ and we fix finite subsets $I'\subset I$
and $F\subset\mathcal{K}$ such that
$\text{Im}(f)\subseteq\oplus_{i\in I'}e_iX$ and $f(\nu_kM)=0$, for all $k\in\mathcal{K}\setminus F$,
then we have that
$$\text{Im}(\varphi_*(f))=\text{Im}(\varphi\circ f)\subseteq
\varphi (\oplus_{i\in I'}e_iX)\subseteq\oplus_{i\in I'}e_iY,$$
because $\varphi$ is in particular a morphism left $B$-modules, and that
$\varphi_*(f)(\nu_kM)=(\varphi\circ f)(\nu_kM)=0$, for all $k\in\mathcal{K}\setminus F$.
Therefore we have $\varphi_*(f)\in\overline{HOM}_A(M,Y)$.

For condition b), we first prove that $\alpha^*$ is morphism of non-unitary
graded left $B$-modules, which amounts to prove that
$\alpha^*(bf)=(-1)^{|a| |\alpha|}b\alpha^*(f)$, for any homogeneous element $f\in\text{HOM}_A(M,X)$.
On one hand, we have $\alpha^*(bf)=(-1)^{| \alpha| (|b|+|f|)}(bf)\circ\alpha$ while
$(-1)^{|\alpha | |b|}b\alpha^*(f)=(-1)^{|\alpha | |b|}(-1)^{|\alpha| |f|}b(f\circ\alpha)$.
Since $(bf)\circ\alpha =b(f\circ\alpha )$ due to the definition of the multiplication map $B\otimes\text{HOM}_A(N,X)\longrightarrow\text{HOM}_A(N,X)$,
we conclude that $\alpha^*$ is a morphism of degree $|\alpha |$ of graded left $B$-modules.

If $c\in C$ is a homogeneous element, then
$\alpha^*(fc)=(-1)^{|\alpha |(|c|+|f|)}(fc)\circ\alpha$ and
$\alpha^*(f)c=(-1)^{|\alpha| |f|}(f\circ\alpha )c$. When we let
these morphisms act on a homogeneous element $x\in N$, we get
$$[\alpha^*(fc)](x)=(-1)^{|\alpha |(|c|+|f|)}(fc)(\alpha (x))=(-1)^{|\alpha |(|c|+|f|)}f(c\alpha (x))$$
while $$[\alpha^*(f)c](x)=(-1)^{|\alpha| |f|}f(\alpha (cx))=(-1)^{|\alpha| |f|}(-1)^{|\alpha | |c|}f(c\alpha (x)), $$
bearing in mind that $\alpha$ is a morphism of graded left modules
(see the comments after Lemma \ref{lem.left versus right dg module},
applied to non-unitary left $B$-modules). It then follows that $\alpha^*$
is a homogeneous morphism of graded $B-C-$bimodules. On the other hand,
if one considers finite subset $I'\subset I$ and $F\subset\mathcal{K}$
as for condition a), then
$\text{Im}(\alpha^*(f))=\text{Im}(f\circ\alpha )\subseteq\text{Im}(f)\subseteq\oplus_{i\in I'}e_iX$
while $\alpha^*(f)(\nu_kN)=f(\alpha (\nu_kN))\subseteq f(\nu_kM)=0$, for all $k\in\mathcal{K}\setminus F$,
bearing in mind that $\alpha$ is in particular a morphism of left $C$-modules.

To finish the proof, we just need to check that the induced map
\begin{center}
$\overline{HOM}_A(?,?):\text{Hom}_{(C-Dg-A)^{op}\otimes (B-Dg-A)}[(N,X),(M,Y)]=
\text{HOM}_{C-A}(M,N)\otimes\text{HOM}_{B-A}(X,Y)
\longrightarrow\text{HOM}_{B-C}(\overline{HOM}_A(M,X),\overline{HOM}_A(N,Y))$
\end{center}
commutes with the differentials. But what we have done above shows that  the map
$\overline{HOM}_A(?,?)$ is induced by the map
\begin{center}
$\text{HOM}_A(?,?):\text{Hom}_{(C-Dg-A)^{op}\otimes (B-Dg-A)}[(N,X),(M,Y)]=\text{HOM}_{C-A}(M,N)\otimes\text{HOM}_{B-A}(X,Y)\longrightarrow
\text{HOM}_{K}(\text{HOM}_A(M,X),\text{HOM}_A(N,Y))$ (*)
\end{center}
given by the dg functor
$$\text{HOM}_A(?,?):(C-Dg-A)^{op}\otimes (B-Dg-A)
\stackrel{\rho^o\otimes\rho'}{\longrightarrow}(Dg-A)^{op}\otimes (Dg-A)
\stackrel{HOM_A(?,?)}{\longrightarrow}Dg-K.$$
The differential in $\overline{HOM}_A(P,Z)$ is the restriction of that of
$\text{HOM}_A(P,Z)$, for each $P\in C-Dg-A$ and $Z\in B-Dg-A$, and
the differential on $\text{HOM}_{B-C}(\text{HOM}_A(M,X),\text{HOM}_A(N,Y))$
is the restriction of the differential in $\text{HOM}_{K}(\text{HOM}_A(M,X),\text{HOM}_A(N,Y))$.
Therefore $\overline{HOM}_A(?,?)$ commutes with the differentials due to the
fact that the map (*)  commutes with the differentials.
\end{proof}

\medskip

We want to  emphasize  a sort of `dual' situation. Suppose now that $X$
is again a dg $B-A-$bimodule and that $W$ is a dg $B-C-$bimodule. Then the
graded $K$-module $\text{HOM}_{B^{op}}(W,X)$ consisting of the morphisms
$W\longrightarrow X$ in $B-Dg$ should have a structure of non-unitary dg
$C-A-$bimodule. Indeed, we can think of $W$ and $X$ as a dg $C^{op}\otimes B^{op}$-bimodule
and a dg $A^{op}-B^{op}-$ bimodule, respectively. Then the first paragraph
of this section says that $\text{HOM}_{B^{op}}(W,X)$ has a structure of
non-unitary dg $A^{op}-C^{op}$-bimodule, which is equivalent to saying
that it has a structure of non-unitary graded $C-A-$bimodule. Taking then $\overline{HOM}_{B^{op}}(W,X)=C\text{HOM}_{B^{op}}(W,X)A$, we get a (now unitary)
dg $C-A-$bimodule. Our following result makes explicit this structure.

\begin{Cor} \label{cor.right A-structure on HOMB(U,X)}
In the situation of last paragraph, the following assertions hold:

\begin{enumerate}
\item The structure of graded $C-A-$bimodule on $\overline{HOM}_{B^{op}}(W,X)$
is given by the rule $(cfa)(w)=(-1)^{(|c|+|a|)|w|+|c| |f|}f(wc)a$, for all
homogeneous elements $c\in C$, $f\in \overline{HOM}_{B^{op}}(W,X)$, $a\in A$
and $w\in W$.
\item  $\overline{HOM}_{B^{op}}(W,X)$ consists of the $f\in \text{HOM}_{B^{op}}(W,X)$
such that $\text{Im}(f)\subset\oplus_{j\in J'}X\epsilon_j$, for some finite
subset $J'\subseteq J$, and $f(W\nu_k)=0$ for all but finitely many $k\in\mathcal{K}$.
\item The assignment $(W,X)\rightsquigarrow\overline{HOM}_{B^{op}}(W,X)$ is
the definition on objects of a dg functor $(B-Dg-C)^{op}\otimes (B-Dg-A)\longrightarrow C-Dg-A$.
\end{enumerate}
\end{Cor}

\begin{proof}
1) Interpreting $\text{HOM}_{B^{op}}(W,X)$ as a non-unitary dg
$A^{op}-C^{op}-$bimodule, the first paragraph of this section tells us
that this structure is given by the rule $(a^ofc^o)(w)=a^of(c^ow)$.
But, by the identification of modules over a dg algebras as dg modules
on the other side over the opposite dg algebra, we get that
\begin{eqnarray*}
(a^ofc^o)(w)&=&a^of(c^ow)=(-1)^{|c| |w|}a^of(wc)=(-1)^{|c| |w|}(-1)^{|f(wc)| |a|}f(wc)a\\
&=&(-1)^{|c| |w|+|a|(|f|+|w|+|c|)}f(wc)a
\end{eqnarray*}
But, by analogous reason, we have an equality

\begin{center}
$(a^ofc^o)(w)=(-1)^{(|a|+|f|)|c|}[c(a^of)](w)=(-1)^{(|a|+|f|)|c|}(-1)^{|a| |f|}(cfa)(w)=
(-1)^{|a| |c|+|f| |c|+|a| |f|}(cfa)(w)$.
\end{center}
Comparing these two expressions and cancelling signs appearing in both expressions, we get that
$(-1)^{(|c|+|a|)|w|}f(wc)a=(-1)^{|f| |c|}(cfa)(w)$, which gives the equality of assertion 1.

2) Considering the distinguished families of orthogonal idempotents
$(\epsilon_j^o)_{j\in J}$ and $(\nu_k^0)_{k\in\mathcal{K}}$ in $A^{op}$ and $C^{op}$,
respectively, we know that $\overline{HOM}_{B^{op}}(W,X)=A^{op}\text{HOM}_{B^{op}}(W,X)C^{op}$
consists of those $f\in\text{HOM}_{B^{op}}(W,X)$ such that
$\text{Im}(f)\oplus_{j\in J'}\epsilon_j^oX$, for some finite subset
$J'\subseteq J$, and $f(\nu_k^oW)=0$, for all but finitely many $k\in\mathcal{K}$.
Bearing in mind that $\epsilon_j^oX=X\epsilon_j$ and $\nu_k^oW=W\nu_k$, for all
$j\in J$ and $k\in\mathcal{K}$, the assertion follows.

3) is a direct consequence of Proposition \ref{prop.bifunctor}.
\end{proof}

\vspace*{0.3cm}

Let $X$ be again a dg $B-A-$bimodule and let $U$ be a dg $C-B-$bimodule.
Then the dg $K$-module  $U\otimes X:=U\otimes_KX$  has a canonical structure of dg
$C-A-$bimodule by defining $c(u\otimes x)a=(cu)\otimes (xa)$. Clearly, this
multiplication makes $U\otimes X$ into a graded $C-A-$bimodule. Moreover if
$u\in U$, $x\in X$, $c\in C$ and $a\in A$ are homogeneous elements, then we have
\begin{eqnarray*}
d[c(u\otimes x)a]&=&d(cu\otimes xa)\\
&=&d_U(cu)\otimes xa+(-1)^{|cu|}cu\otimes d_X(xa)\\
&=&d_U(cu)\otimes xa+(-1)^{|c|+|u|}cu\otimes (d_X(x)a+(-1)^{|x|}xd_A(a))\\
&=&d_U(cu)\otimes xa+(-1)^{|c|+|u|}cu\otimes d(x)a+(-1)^{|c|+|u|+|x|}cu\otimes xd_A(a)\\
&=&(d_C(c)u+(-1)^{|c|}cd_U(u))\otimes xa+(-1)^{|c|+|u|}cu\otimes d_X(x)a\\
&&+(-1)^{|c|+|u|+|x|}cu\otimes xd_A(a)\\
&=&d_C(c)u\otimes xa+(-1)^{|c|}cd_U(u)\otimes xa+(-1)^{|c|+|u|}cu\otimes d_X(x)a\\
&&+(-1)^{|c|+|u|+|x|}cu\otimes xd_A(a)\\
&=&d_C(c)u\otimes xa +(-1)^{|c|}c[d_U(u)\otimes x+(-1)^{|u|}u\otimes d_X(x)]a\\
&&+ (-1)^{|c|+|u| +|x|}cu\otimes xd_A(a)\\
&=&d_C(c)(u\otimes x)a+(-1)^{|c|}cd_{U\otimes X}(u\otimes x)a+(-1)^{|c|+|u\otimes x|}c(u\otimes x)d_A(a)
\end{eqnarray*}
%
so that the differential of $U\otimes X$ satisfies Leibniz rule as a $C-A-$bimodule.

The $K$-submodule $N$ of $U\otimes X$ generated by all differences $ub\otimes x-u\otimes bx$,
where $u\in U$, $x\in X$ and $b\in B$ are homogeneous elements, is a graded $C-A-$subbimodule
of $U\otimes X$. We will show  that $d(N)\subseteq N$, which will imply that we get an induced
graded map of degree $+1$,
$$d:U\otimes_BX:=\frac{U\otimes X}{N}\longrightarrow \frac{U\otimes X}{N}=U\otimes_BX,$$
making $U\otimes_BX$ into a dg $C-A-$bimodule.
Indeed, we leave to the reader checking the following equality, for all
homogeneous elements $u\in U$, $x\in X$ and $b\in B$:
\begin{eqnarray*}
d(ub\otimes x-u\otimes bx)
&=&d_U(u)b\otimes x-d_U(u)\otimes bx+(-1)^{|u|}(ud_B(b)\otimes x-u\otimes d_B(b)x)\\
&&+(-1)^{|u|+|v|}(ub\otimes d_X(x)-u\otimes bd_X(x)).
\end{eqnarray*}
This shows that  $d(ub\otimes x-u\otimes bx)\in N$ and, hence, that $d(N)\subseteq N$ as desired.

\begin{Prop} \label{prop.tensor functor}
Let $A$, $B$ and $C$ be dg algebras with enough idempotents. The assignment
$(U,X)\rightsquigarrow U\otimes_BX$ is the definition on objects of a dg functor
$$?\otimes_B?:\left(C-Dg-B\right)\otimes \left(B-Dg-A\right)\longrightarrow C-Dg-A.$$
\end{Prop}

\begin{proof}
For simplicity, we denote by $T$ the dg functor that we want to define, so that
$T(U,X)=U\otimes_AX$. If now $\alpha :U\longrightarrow V$ and $\varphi :X\longrightarrow Y$

are homogeneous morphisms in $C-Dg-B$ and $B-Dg-A$, respectively, we define
$T(\alpha\otimes\varphi):U\otimes_BX\longrightarrow V\otimes_BY$ by the rule
$T(\alpha\otimes\varphi )(u\otimes x)=(-1)^{| \varphi | |u|}\alpha (u)\otimes\varphi (x)$,
for all homogeneous elements $u\in U$ and $x\in X$. We first prove that
$T(\alpha\otimes\varphi )$ is well-defined. If $b\in B$ is any homogeneous element,
then we have
\begin{eqnarray*}
T(\alpha\otimes\varphi )(ub\otimes x)&=&(-1)^{| \varphi| (|u|+|b|)}\alpha (ub)\otimes\varphi (x)\\
&=&(-1)^{| \varphi| (|u|+|b|)}\alpha (u)b\otimes\varphi (x)\\
&=&(-1)^{| \varphi| |u|}(-1)^{| \varphi| |b|}\alpha (u)\otimes b\varphi (x)\\
&=&(-1)^{| \varphi| |u|}\alpha (u)\otimes\varphi (bx)\\
&=&T(\alpha\otimes\varphi ) (u\otimes bx)
\end{eqnarray*}
using the facts that $\alpha$ is a morphism of graded right
$B$-modules and $\varphi$ is a morphism of graded left $B$-modules.
Therefore $T(\alpha\otimes\varphi )$ is a well-defined morphism in $GR-K$,
and we clearly have $|T(\alpha\otimes\varphi )|=| \alpha|+| \varphi|$.
It is very easy to see that $T(\alpha\otimes\varphi)$ is morphism in $GR-A$.
On the other hand, if $c\in C$ is a homogeneous element, then we have
\begin{eqnarray*}
T(\alpha\otimes\varphi)[c(u\otimes x)]&=&T(\alpha\otimes\varphi)(cu\otimes x)\\
&=&(-1)^{| \varphi| (|c|+|u|)}\alpha (cu)\otimes\varphi (x)\\
&=&(-1)^{| \varphi| (|c|+|u|)}(-1)^{|\alpha | |c|}c\alpha (u)\otimes\varphi (x)\\
&=&(-1)^{(|\alpha|+|\varphi|) |c|}(-1)^{| \varphi| |u|}c(\alpha (u)\otimes\varphi (x))\\
&=&(-1)^{|T(\alpha\otimes\varphi )| |c|}[cT(\alpha\otimes\varphi)](u\otimes x),
\end{eqnarray*}
bearing in mind that $\alpha$ is a morphism of graded left $C$-modules.
%
It follows that
$T(\alpha\otimes\varphi)[c(u\otimes x)]=
(-1)^{|T(\alpha\otimes\varphi )| |c|}[cT(\alpha\otimes\varphi)](u\otimes x)$,
so that $T(\alpha\otimes\varphi)$ is also a morphism in $C-GR$,
and hence a morphism in $C-GR-A$ (see the comments after Lemma \ref{lem.left versus right dg module}).

We now check conditions 2.a-2.c of Lemma \ref{lem.dg bifunctor}:
\bigskip

\emph{Property 2.c:} Note that we have $T(\alpha\otimes 1_Z)=\alpha\otimes 1_Z,  $
for each dg $B-A-$bimodule $Z$, while
$T(1_W\otimes\varphi )(w\otimes x)=(-1)^{|\varphi| |w|}w\otimes\varphi (x)$,
for each dg $C-B-$module $W$ and all homogeneous elements $w\in W$ and $x\in X$.
We then have
\begin{eqnarray*}
[T(\alpha\otimes 1_Y)\circ T(1_M\otimes\varphi)](u\otimes x)&=&
(-1)^{|\varphi | |u|}T(\alpha\otimes 1_Y)(u\otimes\varphi (x))\\
&=&
(-1)^{|\varphi | |u|}\alpha (u)\otimes\varphi (x)\\
&=&T(\alpha\otimes\varphi )(u\otimes x)\\
&=&(-1)^{|\varphi ||u|}\alpha (u)\otimes\varphi (x)\\
&=&(-1)^{|\varphi ||\alpha|}T(1_V\otimes\varphi) (\alpha (u)\otimes x)\\
&=&(-1)^{|\varphi ||\alpha|}[T(1_V\otimes\varphi)\circ T(\alpha\otimes 1_X)](u\otimes x)
\end{eqnarray*}
%
for all homogeneous elements $u\in U$ and $x\in X$.
Therefore condition 2.c of the mentioned lemma is satisfied.

\bigskip

\emph{Property 2.a} If $U$ is a fixed dg $C-B-$bimodule and we consider
the assignments $T_U:B-Dg-A\longrightarrow C-Dg-A$ given by
$T_U(X)=U\otimes_BX$ on objects and $T_U(\varphi )=T(1_U\otimes\varphi )$
on morphisms, we need to check that $T_U$ is a dg functor. We have
$T_U(1_X)=T(1_U\otimes 1_X):u\otimes x\rightsquigarrow (-1)^{|1_X| |u|}u\otimes x=u\otimes x$,
so that $T_U(1_X)=1_{T_U(X)}$. Moreover, if
$\varphi :X\longrightarrow Y$ and $\psi :Y\longrightarrow Z$
are homogeneous morphisms in $B-Dg-A$, then we have
\begin{eqnarray*}
T_U(\psi\circ\varphi)(u\otimes x)&=&T(1_U\otimes (\psi\circ\varphi))(u\otimes x)\\
&=&(-1)^{(| \varphi|+| \psi|)|u|}u\otimes(\psi\circ\varphi)(x)\\
&=&(-1)^{|\varphi | |u|}(-1)^{|\psi | |u|}u\otimes (\psi\circ\varphi)(x)\\
&=&(-1)^{|\varphi | |u|}T(1_U\otimes\psi)(u\otimes\varphi (x))\\
&=&[T(1_U\otimes\psi)\circ T(1_U\otimes\varphi)](u\otimes x)\\
&=&[T_U(\psi)\circ T_U(\varphi)](u\otimes x)
\end{eqnarray*}
%
It then follows that $T_U$ is a graded functor. We need to see that
it commutes with the differentials, which means that the diagram
$$\xymatrix{\text{HOM}_{B-A}(X,Y)\ar[r]^d\ar[d]^{T_U}&\text{HOM}_{B-A}(X,Y)\ar[d]^{T_U}\\
\text{HOM}_{C-A}(U\otimes_BX,U\otimes_BY)\ar[r]^{\delta}&
\text{HOM}_{C-A}(U\otimes_BX,U\otimes_BY)
}$$ commutes, where $d$ and $\delta$ are the differentials on Hom spaces of
$B-Dg-A$ and $C-Dg-A$, respectively. We fix any homogeneous element
$\varphi\in\text{HOM}_{B-A}(X,Y)$ and shall prove that
$(\delta\circ T_U)(\varphi )=(T_U\circ d)(\varphi)$.
Letting act the two members of the desired equality on $u\otimes x$,
where $u\in U$ and $x\in X$ are homogeneous elements, we get:
\begin{eqnarray*}
[(\delta\circ T_U)\lefteqn{(\varphi )](u\otimes x)=}\\
&=&[d_{U\otimes_BY}\circ T_U(\varphi )-(-1)^{|\varphi|}T_U(\varphi)\circ d_{U\otimes_BX}](u\otimes x)\\
&=&(-1)^{|\varphi | |u|}d_{U\otimes_BY}(u\otimes\varphi (x))-
(-1)^{|\varphi|}T_U(\varphi)(d_U(u)\otimes x+(-1)^{|u|}u\otimes d_X(x))\\
&=&(-1)^{|\varphi | |u|}[d_U(u)\otimes\varphi (x)+(-1)^{|u|}u\otimes d_Y(\varphi (x))]\\
&&-(-1)^{|\varphi|}[(-1)^{|\varphi |(|u|+1)}d_U(u)\otimes\varphi (x)+
(-1)^{|u|}(-1)^{|\varphi| |u|}u\otimes\varphi (d_X(x))]\\
&=&(-1)^{(|\varphi |+1)|u|}u\otimes d_Y(\varphi (x))-
(-1)^{(|\varphi |+1)|u|}(-1)^{|\varphi|}u\otimes\varphi (d_X(x))\\
&=&(-1)^{(|\varphi |+1)|u|}u\otimes d(\varphi )(x)\\
&=&T(1_U\otimes d(\varphi ))(u\otimes x)\\
&=&[(T_U\circ d)(\varphi )](u\otimes x).
\end{eqnarray*}

\bigskip

\emph{Condition 2.b} Let us fix a dg $B-A-$bimodule $X$ and consider the
assignments $T^X=?\otimes_BX: (C-Dg-B)\longrightarrow C-Dg-A$ given on objects
by $U\rightsquigarrow U\otimes_BX$ and on morphisms by
$\alpha\rightsquigarrow T(\alpha\otimes 1_X)=\alpha\otimes 1_X$.
It is straightforward to see that
$T^X(\beta\circ\alpha)=T^X(\beta)\circ T^X(\alpha)$, whenever $\alpha$ and $\beta$
are composable morphisms in $C-Dg-B$, and that $T^X(1_U)=1_{T^X(U)}$, so that
we have a graded functor   $C-Dg-B\longrightarrow C-Dg-A$. It remains to see that
$T^X$ commutes with the differentials.
For that, we fix arbitrary  dg $C-B$-bimodules $U$ and $V$ and denote by
$d:HOM_B(U,V)\longrightarrow HOM_B(U,V)$  and
$\delta :HOM_A(U\otimes_BX,V\otimes_BX)\longrightarrow HOM_A(U\otimes_BX,V\otimes_BX)$
the respective differentials on $Hom$ spaces. We need to check
that $(?\otimes_BX)(d(\alpha ))=\delta [(?\otimes_BX)(\alpha )]$.
That is, we need to check that $d (\alpha )\otimes 1_X=\delta (\alpha\otimes 1_X)$,
for any homogeneous element $\alpha\in HOM_B(U,V)$. But if $u\in U$ and
$x\in X$ are homogeneous elements, then we have  an equality
\begin{eqnarray*}
\delta (\alpha\otimes 1_X)(u\otimes x)&=&[d_{V\otimes_BX}\circ (\alpha\otimes 1_X)-
(-1)^{|\alpha\otimes 1_X|}(\alpha\otimes 1_X)\circ d_{U\otimes_BX}](u\otimes x)\\
&=&d_{V\otimes_BX}(\alpha (u)\otimes x)-
(-1)^{|\alpha|}(\alpha\otimes 1_X)[d_U(u)\otimes x+(-1)^{|u|}u\otimes d_X(x)]\\
&=&d_V(\alpha (u))\otimes x+(-1)^{|\alpha (u)|}\alpha (u)\otimes d_X(x)-
(-1)^{|\alpha|}\alpha (d_U(u))\otimes x\\
&&-(-1)^{|\alpha|+|u|}\alpha (u)\otimes d_X(x)\\
&=&d_V(\alpha (u))\otimes x-(-1)^{|\alpha|}\alpha (d_U(u))\otimes x\\
&=&(d_V\circ\alpha -(-1)^{|\alpha|}\alpha\circ d_U)(u)\otimes x\\
&=&d(\alpha )(u)\otimes x\\
&=&(d (\alpha )\otimes 1_X)(u\otimes x).
\end{eqnarray*}
\end{proof}

\section{The classical dg adjunctions} \label{section.tensor product and adjoint}

In this section we  show that the classical tensor-Hom adjunction and the
adjunction between contravariant Hom functors for module categories over rings with
unit can be extended to the dg setting.

\begin{Theorem} \label{thm.adjunction tensor-HOM}
Let $A$, $B$ and $C$ be dg algebras with enough idempotents and let $X$
be a dg $B-A-$bimodules. The pair
$$(?\otimes_BX:C-Dg-B\longrightarrow C-Dg-A, \overline{HOM}_A(X,?):C-Dg-A\longrightarrow C-Dg-B)$$
is a dg adjunction.  As a consequence, we have an adjunction
$$(?\otimes_B^\mathbb{L}X:\mathcal{D}(B\otimes C^{op})
\longrightarrow\mathcal{D}(A\otimes C^{op}),\mathbb{R}\text{Hom}_A(X,?):\mathcal{D}(A\otimes C^{op})\longrightarrow
\mathcal{D}(B\otimes C^{op}))$$
of triangulated functors, where $?\otimes_B^\mathbb{L}:=\mathbb{L}(?\otimes_BX)$ and $\mathbb{R}\text{Hom}_A(X,?):=\mathbb{R}(\overline{HOM}_A(X,?))$.
\end{Theorem}
\begin{proof}
As usual, we fix distinguished families of orthogonal idempotents
$(\epsilon_j)_{j\in J}$, $(e_i)_{i\in I}$ and $(\nu_k)_{k\in\mathcal{K}}$ in $A$, $B$ and $C$,
respectively.
Whenever $U$ and $M$ are objects in $C-Dg-B$ and $C-Dg-A$, respectively, we define
$$\eta_{U,M}:\text{HOM}_{C-A}(U\otimes_BX,M)\longrightarrow\text{HOM}_{C-B}(U,\overline{HOM}_A(X,M))$$
by the rule $[\eta_{U,M}(f)(u)](x)=f(u\otimes x)$, for all homogeneous elements
$f\in\text{HOM}_{C-A}(U\otimes_BX,M)$, $u\in U$ and $x\in X$. We need to check that
$\eta$ is well-defined. We start by checking that $\eta (f)(u)\in\overline{HOM}_A(X,M)$.
If $a\in A$ is any homogeneous element, then we have
$$[\eta (f)(u)](xa)=f(u\otimes xa)=f(u\otimes x)a=[\eta (f)(u)](x)a,$$ so that $\eta (f)(u)$
is  a homogeneous element of $\text{HOM}_A(X,M)$. Moreover if $i\in I$
is such that $ue_i=0$, then we have that $\eta (f)(u)(e_iX)=f(u\otimes e_iX)=0$,
because $f(u\otimes e_ix)=f(ue_i\otimes x)=0$. It then follows that $\eta (f)(u)$
vanishes on all but finitely many $e_iX$. On the other hand,  we know that
there is a finite subset $F\subset\mathcal{K}$ such that $\nu_ku=0$, for all
$k\in\mathcal{K}\setminus F$. It then follows that
$$[(\nu_k\eta (f))(u)](x)=\nu_kf(u\otimes x)=f(\nu_ku\otimes x)=0,$$ for all
$k\in\mathcal{K}\setminus F$. It follows that $\text{Im}(\eta (f)(u))\subseteq\sum_{k\in F}\nu_kM$.
Therefore we get that $\eta (f)(u)\in\overline{HOM}_A(X,M)$.

We next check that  $\eta (f):U\longrightarrow\overline{HOM}_A(X,M)$ is a morphism
of dg $C-B$-bimodules. We need to check that $$\eta (f)(cub)=(-1)^{|c| | \eta (f)|}c(\eta (f)(u))b=(-1)^{|c| |f|}c(\eta (f)(u))b,$$ for all homogeneous elements $c\in C$, $u\in U$ and $b\in B$ (see  the comments after Lemma \ref{lem.left versus right dg module}).
Indeed  we have
\begin{eqnarray*}
[\eta (f)(cub)](x)&=&f(cub\otimes x)\\
&=&(-1)^{|c| |f|}cf(u\otimes bx)\\
&=&(-1)^{|c| |f|}c[\eta (f)(u)](bx)\\
&=&(-1)^{|c| |f|}[c(\eta (f)(u))b](x),
\end{eqnarray*}
due to the definition of the structure of  $C-B-$bimodule on $\overline{HOM}_A(X,M)$.
It follows that $\eta =\eta_{U,M}$ is well-defined.

For the naturality of $\eta$, recall that $\text{HOM}_A(?,M):(C-Dg-A)^{op}\longrightarrow Dg-K$
takes a homogeneous morphism $\alpha :N\longrightarrow N'$ in $C-Dg-A$
to the map
$$\alpha^*=\text{HOM}_A(\alpha^o,M):\text{HOM}_A(N',M)\longrightarrow\text{HOM}_A(N,M)$$
given by $\alpha^*(\beta )=(-1)^{|\alpha | |\beta|}\beta\circ\alpha$
(see the proof of Proposition \ref{prop.bifunctor}). A similar fact is true for
$\text{HOM}_{C-B}(?,W):(C-Dg-B)^{op}\longrightarrow Dg-K$, for any  dg $C-B-$bimodule $W$.
With this in mind, if $\varphi :U\longrightarrow V$ is a homogeneous morphism in $C-Dg-B$,
then we have that
$$[\varphi^*\circ\eta_{V,M}](g)=(-1)^{|\eta (g)| |\varphi |}\eta_{V,M}(g)\circ\varphi =
(-1)^{|g| |\varphi |}\eta_{V,M}(g)\circ\varphi,$$
for each homogeneous element $g\in\text{HOM}_{C-A}(V\otimes_BX,M)$. On the
other hand, we have
$$[\eta_{U,M}\circ (\varphi\otimes 1_X)^*](g)=(-1)^{|g| |\varphi |}\eta_{U,M}(g\circ (\varphi\otimes 1_X)).$$
Taking homogeneous elements $u\in U$ and $x\in X$, we then have
\begin{eqnarray*}
[(\varphi^*\circ\eta_{V,M})(g)](u)(x)
&=&(-1)^{|g| |\varphi |}[\eta_{V,M}(g)\circ\varphi ](u)(x)=(-1)^{|g| |\varphi |}g(\varphi (u)\otimes x)\\
&=&(-1)^{|g| |\varphi |}[g\circ (\varphi\otimes 1_X)](u\otimes x)\\
&=&(-1)^{|g| |\varphi |}\eta_{U,M}(g\circ (\varphi\otimes 1_X))(u)(x)\\
&=&[(\eta_{U,M}\circ (\varphi\otimes 1_X)^*)(g)](u)(x).
\end{eqnarray*}
This shows that  $\varphi^*\circ\eta_{V,M}=\eta_{U,M}\circ (\varphi\otimes 1_X)^*$,
which proves the naturality of $\eta$ on the variable $U$. The naturality on the
variable $M$ is shown as in the classical (ungraded) context.

It remains to prove that $\eta$ commutes with the differentials. For this, we denote by
$d:\text{HOM}_{C-A}(U\otimes_BX,M)\longrightarrow\text{HOM}_{C-A}(U\otimes_BX,M)$ and
$\delta :\text{HOM}_{C-B}(U,\overline{HOM}_A(X,M))\longrightarrow\text{HOM}_{C-B}(U,\overline{HOM}_A(X,M))$
the respective differentials on Hom spaces in the dg categories $C-Dg-A$ and $C-Dg-B$,
respectively. We need to prove that $\delta (\eta (f))=\eta (d(f))$, for each
homogeneous element $f\in \text{HOM}_{C-A}(U\otimes_BX,M)$. If $u\in U$ and $x\in X$
are  homogeneous elements, then we have:
\begin{eqnarray*}
[\delta (\eta (f))](u)(x)&=&[d_{\overline{HOM}_A(X,M)}\circ\eta (f)-(-1)^{|\eta (f)|}\eta (f)\circ d_U](u)(x)\\
&=&[d_{\overline{HOM}_A(X,M)}(\eta (f)(u))-(-1)^{|f|}\eta (f)(d_U(u))](x)\\
&=&[d_M\circ\eta (f)(u)-(-1)^{|u|+|f|}\eta (f)(u)\circ d_X-(-1)^{|f|}\eta (f)(d_U(u))](x)\\
&=&d_M(f(u\otimes x))-(-1)^{|u|+|f|}f(u\otimes d_X(x))-(-1)^{|f|}f(d_U(u)\otimes x)\\
&=&d_M(f(u\otimes x))-(-1)^{|f|}f(d_U(u)\otimes x+(-1)^{|u|}u\otimes d_X(x))\\
&=&d_M(f(u\otimes x))-(-1)^{|f|}f(d_{U\otimes_BX}(u\otimes x))\\
&=&(d_M\circ f-f\circ d_{U\otimes_BX})(u\otimes x)\\
&=&d(f) (u\otimes x)\\
&=&\eta (d(f))(u)(x).
\end{eqnarray*}
Therefore we have $\delta (\eta (f))=\eta (d(f))$, as desired.
\end{proof}

\medskip

We now consider a particular case of the last adjunction.
Let $\iota :B\longrightarrow A$ be  a homomorphism of dg algebras with
enough idempotents. All throughout the rest of the  paper, we assume that such a
homomorphism makes $A$ into a (unitary!)   $B-B-$bimodule (equivalently,
 that  $A=\iota (B)A\iota (B)$). This means that if $(e_i)_{i\in I}$
is any distinguished family of orthogonal idempotents of $B$ then, after
deleting those $\iota (e_i)$ which are zero, the family $(\iota (e_i){i\in I})$
is a distinguished family of orthogonal idempotents of $A$. Note that we have
an obvious \emph{restriction of scalars functor}
$\iota_*:C-Dg-A\longrightarrow C-Dg-B$, which is clearly a dg functor that preserves  acyclic and contractible
dg modules.  In particular, we have $\mathbb{R}\iota_*=\iota_*$ (see
Remark \ref{rem.derived functor of exact dg functor}).
We can apply  the last proposition  to the bimodule $X={}_BA_A$.
But note the following:

\begin{Lemma} \label{lem.iso restriction scalars}
In the situation of preceding paragraph, consider the dg functor
$$\overline{HOM}_A(A,?):C-Dg-A\longrightarrow C-Dg-B.$$ There is a
natural isomorphism of dg functors   $\iota_*\cong\overline{HOM}_A(A,?)$.
As a consequence, there is a natural isomorphism of triangulated functors
$$\mathbb{R}\iota_*=\iota_*\cong\mathbb{R}\text{Hom}_A(A,?):\mathcal{D}(A\otimes C^{op})\stackrel{\cong}{\longrightarrow}\mathcal{D}(B\otimes C^{op}).$$
\end{Lemma}

\begin{proof}
Recall that if $M$ is a dg $C-A-$bimodule, then  $\overline{HOM}_A(A,M)$
consists of the morphisms $f:A\longrightarrow M$ in $Dg-A$ such that $f(e_iA)=0$,
for all but finitely many $i\in I$, and $\text{Im}(f)\subseteq\bigoplus_{k\in F}\nu_kM$,

for some finite subset $F\subset\mathcal{K}$. If $m\in M$ is any homogeneous element
and we consider the homogeneous morphism $\lambda_m:A\longrightarrow M$ in
$GR-A$ given by $\lambda_m(a)=ma$, then $\lambda_m\in\overline{HOM}_A(A,M)$.
Indeed since $me_i=0$, for all but finitely many $i\in I$, we get that also
$\lambda_m(e_iA)=me_iA=0$, for all but finitely many $i\in I$. On the other
hand, since there is a finite subset $F\subset\mathcal{K}$ such that
$\nu_km=0$, for all $k\in\mathcal{K}\setminus F$, we get that
$\text{Im}(\lambda_m)=mA\subseteq\oplus_{k\in F}\nu_kM$.

The induced map $\lambda_M :M\longrightarrow \overline{HOM}_A(A,M)$
is clearly a morphism in $C-GR-B$. Defining
$\Psi :\overline{HOM}_A(A,M)\longrightarrow M$ by the rule $\Psi (f)=\sum_{i}f(e_i)$,
we get an inverse for $\lambda$ in $C-GR-B$. Then, for each $M$ in $C-Dg-A$,
we have a morphism of degree zero
$\lambda_M:\iota_*(M)\longrightarrow\overline{HOM}_A(A,M)$ in $C-Dg-B$.
To check that, when $M$ varies, we get a bijective natural transformation
of dg functors $\lambda :\iota_*\longrightarrow\overline{HOM}_A(A,?)$
is easy and left to the reader.
In order to see that we have a natural isomorphism of dg functors we
just need to check that $\lambda$ is a homological natural transformation,
which amounts to check that if $d:\text{HOM}_{C-B}(M,\overline{HOM}_A(A,M))\longrightarrow\text{HOM}_{C-B}(M,\overline{HOM}_A(A,M))$
is the differential on Hom spaces of the dg category $C-Dg-B$,
then $d(\lambda_M)=0$. To see this, for each homogeneous element $m\in M$, we have:
\begin{eqnarray*}
[d(\lambda_M)](m)&=&[d_{\overline{HOM}_A(A,M)}\circ\lambda_M-(-1)^{|\lambda_M|}\lambda_M\circ d_M](m)\\
&=&d_{\overline{HOM}_A(A,M)}(\lambda_M(m))-\lambda_M(d_M(m))\\
&=&d_M\circ\lambda_M(m)-(-1)^{|m|}\lambda_M(m)\circ d_A-\lambda_M(d_M(m)).
\end{eqnarray*}
When applying both members of this equality to a homogeneous element $a\in A$, we get that $$\{[d(\lambda_M)](m)\}(a)=d_M(ma)-(-1)^{|m|}md_A(a)-d_M(m)a =0.$$

The corresponding natural isomorphism of triangulated functors
follows from Proposition \ref{prop.dg transformation yields triang transformation}.
\end{proof}

\medskip

This justifies the following terminology:

\begin{Def}
If $\iota :B\longrightarrow A$ is a homomorphism of dg algebras
with enough idempotents as above, then the dg functor
$?\otimes_BA:C-Dg-B\longrightarrow C-Dg-A$ is called the
\emph{extension of scalars functor} associated to $\iota$.
It is denoted by $\iota^*:C-Dg-B\longrightarrow C-Dg-A$.
\end{Def}

As an immediate consequence of Theorem \ref{thm.adjunction tensor-HOM} and Lemma \ref{lem.iso restriction scalars}, we get:

\begin{Cor} \label{cor.adjunction extension-restriction scalars}
Let $\iota :B\longrightarrow A$ a homomorphism of dg algebras as above.
The pair $(\iota^*:C-Dg-B\longrightarrow C-Dg-A,\iota_*:C-Dg-A\longrightarrow C-Dg-B)$
is a dg adjunction. Therefore we have an adjoint pair of triangulated functors
$(\mathbb{L}\iota^*:\mathcal{D}(B\otimes C^{op})\longrightarrow\mathcal{D}(A\otimes C^{op}),\mathbb{R}\iota_*=
\iota_*:\mathcal{D}(A\otimes C^{op})\longrightarrow\mathcal{D}(B\otimes C^{op}))$.
\end{Cor}

We move now to study a less known adjunction.

\begin{Theorem} \label{thm.adjunction X-reflexive}
Let $A$ and $B$ be dg algebras with enough idempotents
and let $X$ be a dg $B-A-$bimodule. The pair
$$(\overline{HOM}_{B^{op}}(?,X)^o:B-Dg-C\longrightarrow (C-Dg-A)^{op},
\overline{HOM}_A(?,X):(C-Dg-A)^{op}\longrightarrow B-Dg-C)$$ is a dg adjunction.
In particular, the pair $$(\mathbb{R}\text{Hom}_{B^{op}}(?,X)^o:\mathcal{D}(C\otimes B^{op})\longrightarrow
\mathcal{D}(A\otimes C^{op})^{op},\mathbb{R}\text{Hom}_{A}(?,X):\mathcal{D}(A\otimes C^{op})^{op}\longrightarrow
\mathcal{D}(C\otimes B^{op}))$$ is an adjoint pair of triangulated functors, where
$\mathbb{R}\text{Hom}(?,X):=\mathbb{R}(\overline{HOM}(?,X))$ in both cases.
\end{Theorem}

\begin{proof}
All throughout the proof, we fix distinguished families of orthogonal
idempotents $(e_i)_{i\in I}$,  $(\epsilon_j)_{j\in J}$ and $(\nu_k)_{k\in\mathcal{K}}$ in $B$, $A$ and $C$, respectively. Let $U$ be a dg $B-C-$bimodule and $M$ be a dg $C-A-$bimodule. By the initial paragraph of Section \ref{section.contravariant HOM} and by Corollary \ref{cor.right A-structure on HOMB(U,X)}, we have:

\begin{enumerate}
\item[a)] $\overline{HOM}_{B^{op}}(U,X)$ consists of the $f\in\text{HOM}_{B^{op}}(U,X)$ such that $\text{Im}(f)\subseteq\oplus_{j\in J'}X\epsilon_j$, for some finite subset $J'\subset J$, and $f(U\nu_k)=0$ for all but finitely many $k\in\mathcal{K}$.

\item[b)] $\overline{HOM}_{A}(M,X)$ consists of the  $g\in\text{HOM}_{A}(M,X)$ such that $\text{Im}(g)\subseteq\oplus_{i\in I'}e_iX$, for some finite subset $I'\subset I$, and $g(\nu_kM)=0$ for all but finitely many $k\in\mathcal{K}$.
\end{enumerate}

We define a $K$-linear map
$$\xymatrix{
\text{Hom}_{(C-Dg-A)^{op}}(\overline{HOM}_{B^{op}}(U,X),M)\ar@{=}[d]&&
\text{Hom}_{B-Dg-C}(U,\overline{HOM}_A(M,X))\ar@{=}[d]\\
\text{HOM}_{C-A}(M,\overline{HOM}_{B^{op}}(U,X))\ar[rr]^{\xi=\xi_{U,M}}&&
\text{HOM}_{B-C}(U,\overline{HOM}_A(M,X))}$$
by the rule $[\xi (f)(u)](m)=(-1)^{|u| |m|}f(m)(u)$, for all homogeneous elements
$u\in U$, $m\in M$, and
$f\in\text{HOM}_{C-A}(M,\overline{HOM}_{B^{op}}(U,X))$. We first check
that if $f$ and $u$ are fixed, then the assignment
$m\rightsquigarrow [\xi (f)(u)](m)=(-1)^{|u| |m|}f(m)(u)$
gives a homogeneous morphism $M\longrightarrow X$ in $GR-A$.
Indeed we have
$$[\xi (f)(u)](ma)=(-1)^{|u| (|m|+|a|)}f(ma)(u)=
(-1)^{|u| (|m|+|a|)}[f(m)a](u).$$
But, by the structure of  right
dg $A$-module on $\overline{HOM}_{B^{op}}(U,X)$ (see
Corollary \ref{cor.right A-structure on HOMB(U,X)}),  we see
$[f(m)a](u)=(-1)^{|a| |u|}f(m)(u)a$, hence
$[\xi (f)(u)](ma)=(-1)^{|u| |m|}f(m)(u)a$. On the other hand,
we get $[\xi (f)(u)](m)a=(-1)^{|u| |m|}f(m)(u)a$.
This shows that $\xi (f)(u)$  is a homogeneous morphism $M\longrightarrow X$ in $GR-A$.

In order to check that $\xi$ is well-defined, we also need to see that the just
defined morphism  $\xi (f)(u)$ is really in $\overline{HOM}_A(M,X)$ (see point b) above).
We have $u=\sum_{i\in F_u}e_iu$,
for some finite subset $F_u\subseteq I$, and then
$[f(m)](u)=\sum_{i\in F_u}e_if(m)(u)\in\bigoplus_{i\in F_u}e_iX$, for all
$m\in M$, using the fact that $f(m)$ is a morphism of graded left $B$-modules. That is, we have that $\text{Im}(\xi (f)(u))\subseteq\oplus_{i\in F_u}e_iX$. On the other hand, we have a finite subset $\mathcal{K}'\subset\mathcal{K}$ such that $u\nu_k=0$, for all $k\in\mathcal{K}\setminus\mathcal{K}'$. Bearing in mind the explicit definition of the $C-A-$bimodule structure on $\overline{HOM}_{B^{op}}(U,X)$ (see Corollary \ref{cor.right A-structure on HOMB(U,X)}) and the fact that $f$ is a morphism of $C-A-$bimodules, we get that
\begin{eqnarray*}
[\xi (f)(u)](\nu_km)&=&(-1)^{|u| |\nu_km|}f(\nu_km)(u)\\
&=&(-1)^{|u| |m|}(\nu_kf(m))(u)\\
&=&(-1)^{|u| |m|} (-1)^{|\nu_k| |u|+|\nu_k| |f(m)|}f(m)(u\nu_k),
\end{eqnarray*}
which implies that $[\xi (f)(u)](\nu_km)=0$, for all $k\in\mathcal{K}\setminus\mathcal{K}'$, independently of $m$. It follows that $\xi (f)(u)](\nu_kM)=0$ for almost all $k\in\mathcal{K}$, so that $\xi (f)(u)\in\overline{HOM}_A(M,X)$.

As a final step to check that $\xi$ is well-defined, we will see that $\xi (f)$
is really a morphism   $U\longrightarrow\overline{HOM}_A(M,X)$ in $B-GR-C$.
That is, we need to check the equalities $\xi (f)(bu)=(-1)^{|f| |b|}b\xi (f)(u)$
(see  the comments after Lemma \ref{lem.left versus right dg module}) and $\xi (f)(uc)=\xi(f)(u)c$, for all homogeneous elements $b\in B$, $u\in U$ and $c\in C$.
We apply both members of the first desired equality to a homogeneous element $m\in M$ and get:
\begin{eqnarray*}
[\xi (f)(bu)](m)&=&(-1)^{(|b|+|u|) |m|}f(m)(bu)\\
&=&(-1)^{(|b|+|u|) |m|}(-1)^{|f(m)| |b|}bf(m)(u)\\
&=&(-1)^{(|b|+|u|) |m|}(-1)^{((||f|+m|) |b|}bf(m)(u)\\
&=&(-1)^{|u| |m|+|f| |b|}bf(m)(u)\\
&=&(-1)^{|f| |b|}(-1)^{|u| |m|}bf(m)(u)\\
&=&(-1)^{|f| |b|}b[\xi (f)(u)(m)]\\
&=&(-1)^{|f| |b|}[b\xi (f)(u)](m),
\end{eqnarray*}
using that $f(m)$ is a morphism in $B-GR$ of degree $|f|+|m|$. On the other hand, we have
$$[\xi (f)(uc)](m)=(-1)^{|uc| |m|}f(m)(uc)=(-1)^{|u| |m|+|c| |m|}f(m)(uc)$$ while, using Corollary \ref{cor.right A-structure on HOMB(U,X)} and the comments after Lemma \ref{lem.left versus right dg module}, we also get
\begin{eqnarray*}
[\xi (f)(u)c](m)&=&\xi (f)(u)(cm)\\
&=&(-1)^{|u| |cm|}f(cm)(u)\\
&=&(-1)^{|u| |cm|}(-1)^{|c| |f|}[cf(m)](u)\\
&=&(-1)^{|u| |cm|+|c| |f|}(-1)^{|c| |u|+|c| |f(m)|}f(m)(uc)\\
&=&(-1)^{|u| |m|+|c| |m|}f(m)(uc).
\end{eqnarray*}
We now prove the naturality of $\xi$ on both variables.
Let $\alpha :M\longrightarrow N$ be a homogeneous morphism
in $C-Dg-A$. With the obvious meaning of the vertical arrows,
we need to prove that the following diagram is commutative:
$$
\xymatrix{
HOM_{C-A}(M,\overline{HOM}_{B^{op}}(U,X))\ar[r]^{\xi_{U,M}}&HOM_{B-C}(U,\overline{HOM}_A(M,X))\\
HOM_{C-A}(N,\overline{HOM}_{B^{op}}(U,X))\ar[r]^{\xi_{U,N}}\ar[u]_{\alpha^*}&HOM_{B-C}(U,\overline{HOM}_A(N,X))\ar[u]^{\overline{HOM}_A(\alpha,X)_*}
}
$$

For this, we take any homogeneous morphism
$g:N\longrightarrow\overline{HOM}_{B^{op}}(U,X)$ in $C-Dg-A$.
We then have that $(\xi_{U,M}\circ\alpha^*)(g)=(-1)^{|\alpha| |g| }\xi_{U,M}(g\circ\alpha )$
and therefore
\begin{eqnarray*}
 [(\overline{HOM}_A(\alpha ,X)_*\circ\xi_{U,N}) (g)(u)](m)&=&[(\overline{HOM}_A(\alpha ,X)\circ\xi (g))(u)](m)\\
 &=&(-1)^{|\alpha| |\xi (g)(u)|}[\xi (g)(u)\circ\alpha](m)\\
 &=&(-1)^{|\alpha| (|g|+|u|)}\xi (g)(u)(\alpha (m))\\
 &=&(-1)^{|\alpha| (|g|+|u|)}(-1)^{|u| |\alpha (m)|}g(\alpha (m))(u)\\
 &=&(-1)^{|\alpha| (|g|+|u|)}(-1)^{|u| (|\alpha |+|m|)}(g\circ\alpha )(m)(u)\\
 &=&(-1)^{|\alpha| |g|+|u| |m|}(g\circ\alpha )(m)(u)\\
 &=&(-1)^{|\alpha| |g| }\xi_{U,M}(g\circ\alpha )(u)(m)\\
 &=&[(\xi_{U,M}\circ\alpha^*)(g)(u)](m),
\end{eqnarray*}
which proves the naturality of $\xi$ on the variable $M$.

Let now $\varphi :U\longrightarrow V$ be a homogeneous morphism in $B-Dg-C$.
For the naturality of $\xi$ on the `variable' $U$, we need to check that
the following diagram is commutative:
$$\xymatrix{
HOM_{C-A}(M,\ol{HOM}_{B^{op}}(U,X))\ar[r]^{\xi_{U,M}}&
HOM_{B-C}(U,\ol{HOM}_{A}(M,X))\\
HOM_{C-A}(M,\ol{HOM}_{B^{op}}(V,X))\ar[r]^{\xi_{V,M}}\ar[u]_{{\ol{HOM}_{B^{op}}(\varphi,X)}_*}&
HOM_{B-C}(V,\ol{HOM}_{A}(M,X))\ar[u]^{\varphi^*}
}$$
Let  $f:M\longrightarrow\overline{HOM}_{B^{op}}(V,X)$ be a
homogeneous morphism in $C-Dg-A$. Then we have that
$[\xi_{U,M}\circ\overline{HOM}_{B^{op}}(\varphi ,X)_*](f)=
\xi_{U,M}(\overline{HOM}_{B^{op}}(\varphi ,X)\circ f)$, while
$(\varphi^*\circ\xi_{V,M})(f)=(-1)^{|f| |\varphi |}\xi_{V,M}(f)\circ\varphi$.
If now $u\in U$ and $m\in M$ are homogeneous elements, then we have equalities:
\begin{eqnarray*}
[\xi_{U,M}(\overline{HOM}_{B^{op}}(\varphi ,X)\circ f)(u)](m)
&=&(-1)^{|u| |m|}[(\overline{HOM}_{B^{op}}(\varphi ,X)\circ f)(m)](u)\\
&=&(-1)^{|u| |m|}(-1)^{| \varphi| |f(m)|}[f(m)\circ\varphi](u)\\
&=&(-1)^{|u| |m|+|\varphi | |f|+|\varphi | |m|}f(m)(\varphi (u))\\
&=&(-1)^{|f| |\varphi |}(-1)^{|\varphi (u)| |m|}f(m)(\varphi (u))\\
&=&(-1)^{|f| |\varphi |}[\xi_{V,M}(f)(\varphi (u))](m)\\
&=&(-1)^{|f| |\varphi |}[(\xi_{V,M}(f)\circ\varphi )(u)](m)
\end{eqnarray*}
which proves the naturality of  $\xi$  on the `variable' $U$.

Next we check that $\xi$ commutes with the differentials.
That is, we need to prove that, for each
 dg $C-A-$bimodule $M$ and each  dg $B-C-$bimodule $U$, the following diagram is commutative.
$$\xymatrix{
\text{HOM}_{C-A}(M,\overline{HOM}_{B^{op}}(U,X))\ar[r]^d\ar[d]^{\xi_{U,M}}&
\text{HOM}_{C-A}(M,\overline{HOM}_{B^{op}}(U,X))\ar[d]^{\xi_{U,M}}\\
\text{HOM}_{B-C}(U,\overline{HOM}_{A}(M,X))\ar[r]^{\delta}&
\text{HOM}_{B-C}(U,\overline{HOM}_{A}(M,X)),
}$$
%
%
%
%
Here $d$ and $\delta$ are the differentials on Hom spaces in $C-Dg-A$ and $B-Dg-C$,
respectively. Put $\xi =\xi_{U,M}$ for simplicity. Then,  for all homogeneous elements
$f\in\text{HOM}_{C-A}(M,\overline{HOM}_{B^{op}}(U,X))$, $u\in U$ and $m\in M$, we have
\begin{eqnarray*}
[(\xi\lefteqn{\circ d)(f)(u)](m)=}\\
&=&[\xi (d(f))(u)][m]\\
&=&(-1)^{|u| |m|}d(f)(m)(u)\\
&=&(-1)^{|u| |m|}[(d_{\overline{HOM}(U,X)}\circ f-(-1)^{|f|}f\circ d_M)(m)](u)\\
&=&(-1)^{|u| |m|}[d_X\circ f(m)-(-1)^{|f(m)|}f(m)\circ d_U-(-1)^{|f|}f(d_M(m))](u)\\
&=&(-1)^{|u| |m|}[d_{\overline{HOM}(U,X)}(f(m))-(-1)^{|f|}f(d_M(m))](u)\\
&=&(-1)^{|u| |m|}[d_X(f(m)(u))-(-1)^{|f|}f(d_M(m))(u)-(-1)^{|f|+|m|}f(m)(d_U(u))]\\
&=&(-1)^{|u| |m|}d_X(f(m)(u))-(-1)^{|f|+|u|}(-1)^{|u| (|m|+1)}f(d_M(m))(u)\\
&&-(-1)^{|f|}(-1)^{(|u|+1)|m|}f(m)(d_U(u))\\
&=&[d_X\circ\xi(f)(u)-(-1)^{|\xi(f)(u)|}\xi(f)(u)\circ d_M-(-1)^{|f|}\xi(f)(d_U(u))](m)\\
&=&[d_{\overline{HOM}(M,X)}(\xi (f)(u))-(-1)^{|f|}\xi(f)(d_U(u))](m)\\
&=&[(\delta\circ\xi)(f)(u)][m].
\end{eqnarray*}
where we used that
$(\delta\circ\xi)(f)=
d_{\overline{HOM}(M,X)}\circ\xi (f)-(-1)^{|\xi (f)|}\xi(f)\circ d_U$ in the last equality.
%

Finally, in order to prove the bijective condition of $\xi$,
note that exchanging the roles of $U$ and $X$ and of $A$ and $B^{op}$,
one has a well-defined $K$-linear map of degree zero $\xi'=\xi'_{U,M}:\text{HOM}_{B-C}(U,\overline{HOM}_A(M,X))
\longrightarrow\text{HOM}_{C-A}(M,\overline{HOM}_{B^{op}}(U,X))$,
given by the rule $[\xi' (g)(m)](u)=(-1)^{|u| |m|}g(u)(m)$.
Clearly $\xi'_{M,U}$ is inverse to $\xi_{M,U}$.
\end{proof}

\medskip

We will now show that, under appropriate assumptions, the
derived functors of covariant and contravariant $\overline{HOM}$
are part of a bifunctor which is triangulated on both variables.
We need the following auxiliary result.

\begin{Lemma} \label{lem.preservation of acyclic}
Let $A$, $B$ and $C$ be dg algebras with enough idempotents and let
$P$ and $Q$ be dg $C-A-$bimodules such that  $P$ is homotopically projective  and  $Q$ is   homotopically injective as right
dg $A$-modules. Then the following assertions hold:
\begin{enumerate}
\item The functor $\overline{HOM}_A(P,?):B-Dg-A\longrightarrow B-Dg-C$ preserves acyclic dg bimodules.
\item The functor $\overline{HOM}_A(?,Q):(B-Dg-A)^{op}\longrightarrow C-Dg-B$ preserves acyclic dg bimodules.
\end{enumerate}
\end{Lemma}

\begin{proof}
The proofs of the two assertions are rather similar.
We only prove 1. Let $X$ be an acyclic dg $B-A-$bimodule.
We know that the non-unitary  dg $B-C$-bimodule $\text{HOM}_A(P,X)$
is acyclic since $P_A$ is homotopically projective. By definition,
we have that $\overline{HOM}_A(P,X)=B\text{HOM}_A(P,X)C$ and the
differential on this (unitary)  dg $B-C-$bimodule is the restriction
of the differential of $\text{HOM}_A(P,X)$. Let now
$f\in Z^n(\overline{HOM}_A(P,X))$ be any $n$-cycle. By the
acyclicity of $\text{HOM}_A(P,X)$, we have a $g\in\text{HOM}_A(P,X)^{n-1}$
such that $f=d(g)$, where $d$ is the differential of $\text{HOM}_A(P,X)$.
If $(e_i)_{i\in I}$ and $(\nu_k)_{k\in\mathcal{K}}$ are  distinguished families of orthogonal idempotents of $B$ and $C$, respectively,
then there are finite subset $I'\subset I$ and $\mathcal{K}'\subset\mathcal{K}$ such that $f=\sum_{i\in I',k\in\mathcal{K}'}e_if\nu_k$.
Taking $g'=\sum_{i\in I',k\in\mathcal{K}'}e_ig\nu_k$ and using Leibniz rule for the  non-unitary  dg $B-C$-bimodule $\text{HOM}_A(P,X)$ (see Lemma \ref{lem.B-module structure of M-upstar}), we get an element $g'\in\overline{HOM}_A(P,X)^{n-1}$
such that $$d(g')=d(\sum_{i\in I',k\in\mathcal{K}'}e_ig\nu_k)=\sum_{i\in I',k\in\mathcal{K}'}e_id(g)\nu_k=\sum_{i\in I',k\in\mathcal{K}'}e_if\nu_k=f.$$
\end{proof}

We say that an algebra with enough idempotents $A$ is \emph{$K$-projective} (resp. \emph{K-flat}) when it is projective (resp. flat) as a $K$-module.

\begin{Cor} \label{cor.RHom as part of bifunctor}
Let $A$, $B$ and $C$ be dg algebras with enough idempotents. The dg functor
$\overline{HOM}_A(?,?):(C-Dg-A)^{op}\otimes (B-Dg-A)\longrightarrow B-Dg-C$
preserves contractibility on each variable.
If $\mathbb{R}\text{HOM}_A(?,?):=\mathbb{R} (\overline{HOM}(?,?)):\mathcal{D}(A\otimes C^{op})^{op}\otimes\mathcal{D}(A\otimes B^{op})\longrightarrow\mathcal{D}(C\otimes B^{op})$ is the associated
bi-triangulated functor (see Definition \ref{def.derived functor of a dg bifunctor}),
then the following assertions hold:
\begin{enumerate}
\item If $C$ is $K$-projective or $X$ is a homotopically injective dg $B-A-$bimodule, then there is a natural isomorphism of triangulated functors
$$\mathbb{R}\text{HOM}_A(?,X)\cong\mathbb{R}\text{Hom}_A(?,X):=
    \mathbb{R}(\overline{HOM}_A(?,X)):\mathcal{D}(A\otimes C^{op})^{op}\longrightarrow\mathcal{D}(C\otimes B^{op}).$$

\item If either $B$ is $K$-flat or $M$ is a homotopically projective dg $C-A$-bimodule,
then there is a natural isomorphism of triangulated functors $$\mathbb{R}\text{HOM}_A(M,?)\cong\mathbb{R}\text{Hom}_A(M,?):=\mathbb{R}(\overline{HOM}_A(M,?)):\mathcal{D}(A\otimes B^{op})\longrightarrow\mathcal{D}(C\otimes B^{op}).$$
\end{enumerate}
In particular, if $C$ is $K$-projective (e.g. if $C=K$) one has an isomorphism $\mathbb{R}\text{Hom}_A(M,?)(X)\cong\mathbb{R}\text{Hom}_A(?,X)(M)$ in $\mathcal{D}(C\otimes B^{op})$, for all dg $B-A-$bimodules $X$ and all homotopically projective dg $C-A-$bimodules $M$.
\end{Cor}
\begin{proof}
By Theorems \ref{thm.adjunction tensor-HOM} and \ref{thm.adjunction X-reflexive},
we know that, for fixed $M$ and $X$ in $C-Dg-A$ and $B-Dg-A$, the dg functors
$\overline{HOM}_A(M,?):B-Dg-A\longrightarrow B-Dg-C$ and
$\overline{HOM}_A(?,X):(C-Dg-A)^{op}\longrightarrow B-Dg-C$ are part of a
dg adjunction. By Lemma \ref{lem.derived functor of dg functor}, both of them preserve
contractible dg modules, which shows the first statement of the corollary.

The last statement is a direct consequence of  assertions 1 and 2. Note that when $C$ is $K$-projective (resp. $B$ is $K$-flat), the restriction of scalars functor  $C-Dg-A\longrightarrow Dg-A$ (resp.  $B-Dg-A\longrightarrow Dg-A$) preserves homotopically projective (resp. homotopically injective) dg modules (this is well-known in the context of dg modules over small dg categories, but the reader can easily adapt the proof of \cite[Lemma 3.6]{NS-Japan} to get a direct proof by her/himself). Then assertion 1, when $C$ is $K$-projective,
is a direct consequence of
Proposition \ref{prop.triangulated part of bitriangulated}(2.b) and
Lemma \ref{lem.preservation of acyclic}. Similarly, assertion 2 for  $K$-flat $B$  follows from
 from Proposition \ref{prop.triangulated part of bitriangulated}(2.a) and Lemma \ref{lem.preservation of acyclic}.

To check what remains of assertions 1 and 2, we just prove what remains of assertion 2 since the argument for assertion 1 is entirely dual.   Recall from the proof of
Proposition \ref{prop.triangulated part of bitriangulated} that we
have a natural isomorphism of triangulated functors $\mathcal{D}(A\otimes B^{op})\longrightarrow\mathcal{D}(C\otimes B^{op})$
$$\mathbb{R}\text{HOM}_{A}(M,?)\cong\mathbb{R}(\overline{HOM}_{A}(\Pi_{C-A}(M),?))=
\mathbb{R}\text{Hom}_A(\Pi_{C-A}(M),?).$$
Then $\mathbb{R}\text{HOM}_A(M,?)$ is the composition
$$\mathcal{D}(A\otimes B^{op})
\stackrel{\Upsilon_{B-A}}{\longrightarrow}\mathcal{H}(A\otimes B^{op})
\stackrel{\overline{HOM}_A(\Pi_{C-A}(M),?)}{\longrightarrow}\mathcal{H}(C\otimes B^{op})
\stackrel{q_{C\otimes B^{op}}}{\longrightarrow}\mathcal{D}(C\otimes B^{op}),$$
where $\Upsilon_{B-A}$ (resp. $\Pi_{C-A}$) is the homotopically injective (resp. homotopically projective) resolution functor for dg $B-A-$bimodules (resp. dg $C-A-$bimodules).
But, when $M$ is homotopically projective, the homotopically projective
resolution $\pi :\Pi_{C-A}(M)\longrightarrow M$ is an isomorphism in
$\mathcal{H}(A\otimes C^{op})$. Considering the bi-triangulated functor
$$\overline{HOM}_A(?,?):\mathcal{H}(A\otimes C^{op})^{op}\otimes\mathcal{H}(A\otimes B^{op})
\longrightarrow\mathcal{H}(C\otimes B^{op})$$ (see Proposition \ref{prop.Z0 and H0 for dg bifunctor}), we deduce that $\pi$ induces
a natural isomorphism of triangulated functors
$$\pi^*:\overline{HOM}_{A}(M,?)
\stackrel{\cong}{\longrightarrow}\overline{HOM}_A(\Pi_{C-A}(M),?):\mathcal{H}(A\otimes B^{op})
\longrightarrow\mathcal{H}(C\otimes B^{op}).$$
Then $\mathbb{R}\text{HOM}_A(M,?)$
is naturally isomorphic to the composition
$$\mathcal{D}(A\otimes B^{op})\stackrel{\Upsilon_{B-A}}{\longrightarrow}\mathcal{H}(A\otimes B^{op})
\stackrel{\overline{HOM}_A(M,?)}{\longrightarrow}\mathcal{H}(C\otimes B^{op})
\stackrel{q_{C\otimes B^{op}}}{\longrightarrow}\mathcal{D}(C\otimes B^{op}),$$
which is precisely $\mathbb{R}\text{Hom}_A(M,?):\mathcal{D}(A\otimes B^{op})
\longrightarrow\mathcal{D}(C\otimes B^{op})$.
\end{proof}

\section{Dualities for perfect complexes} \label{section.dualities}

In this final section, we shall consider the adjunction of Theorem \ref{thm.adjunction X-reflexive} when $C=K$. That is, we consider the adjunction $$(\mathbb{R}\text{Hom}_{B^{op}}(?,X)^o:\mathcal{D}(B^{op})\longrightarrow\mathcal{D}(A)^{op},\mathbb{R}\text{Hom}_{A}(?,X):\mathcal{D}(A)^{op}\longrightarrow\mathcal{D}(B^{op})).$$

\begin{Rem} \label{rem.unit and counit}
We denote the unit of this adjunction by
$$\lambda :1_{\mathcal{D}(B^{op})}
\longrightarrow\mathbb{R}\text{Hom}_{A}(?,X)\circ\mathbb{R}\text{Hom}_{B^{op}}(?,X)^o.$$
Note that the counit
$\rho^o:\mathbb{R}\text{Hom}_{B^{op}}(?,X)^o\circ\mathbb{R}\text{Hom}_{A}(?,X)
\longrightarrow 1_{\mathcal{D}(A)^{op}}$, when evaluated at any right dg $A$-module,
is a morphisms in $\mathcal{D}(A)^{op}$. We then change this perspective,
and see it as natural transformation $\rho :1_{\mathcal{D}(A)}
\longrightarrow \mathbb{R}\text{Hom}_{B^{op}}(?,X)\circ\mathbb{R}\text{Hom}_{A}(?,X)^o$.
\end{Rem}

Recall that if
$(F:\mathcal{C}\longrightarrow\mathcal{D},G:\mathcal{D}\longrightarrow\mathcal{C})$
is an adjoint pair of arbitrary categories $\mathcal{C}$ and
$\mathcal{D}$, an object $C\in\text{Ob}(\mathcal{C})$
(resp. $D\in\text{Ob}(\mathcal{D})$) is called \emph{reflexive} (resp \emph{coreflexive})
with respect to the given adjunction when the evaluation of the unit at
$C$ (resp. the evaluation of the counit at $D$) is an isomorphism.
The following fact is well-known:

\begin{Lemma} \label{lem.equivalence reflexive-coreflexive}
In the situation of the previous paragraph, the functors $F$ and $G$
define by restriction mutually quasi-inverse equivalences of categories
between the full subcategories of reflexive and coreflexive objects.
\end{Lemma}

Coming back to the situation of Theorem \ref{thm.adjunction X-reflexive}, with $C=K$,
the following definition comes then naturally.

\begin{Def} \label{def.X-(co)reflexive objects}
A left dg $B$-module $U$ will be called \emph{homologically $X$-reflexive}
when the unit map
$$\lambda_U:U\longrightarrow\mathbb{R}\text{Hom}_{A}(\mathbb{R}\text{Hom}_{B^{op}}(U,X),X)$$
is an isomorphism. A right dg $A$-module $M$ is called
\emph{homologically $X$-coreflexive} when the counit map
$$\rho_M:M\longrightarrow\mathbb{R}\text{Hom}_{B^{op}}(\mathbb{R}\text{Hom}_{A}(M,X),X)$$
is an isomorphism. Fixing again distinguished families of orthogonal
idempotents $(e_i)_{i\in I}$ and $(\epsilon_j)_{j\in J}$ in $B$ and $A$,
respectively, we shall say that the dg $B-A-$bimodule $X$ is
\emph{left (resp. right) homologically faithfully balanced} (see \cite{NS-GTT})
when each $Be_i$ (resp. $\epsilon_jA$) is homologically $X$-reflexive
(resp. homologically $X$-coreflexive).
We will say that $X$ is \emph{homologically faithfully balanced}
when it is left and right homologically faithfully balanced.
\end{Def}

 In the sequel we denote by $\text{per}(A)$ (resp. $\text{per}(A^{op})$) the (thick) subcategory of $\mathcal{D}(A)$ (resp. $\mathcal{D}(A^{op}$)) formed by the compact objects. It will be called the \emph{perfect right (resp. left) derived category of $A$}.

Recall that if $\mathcal{C}$ and $\mathcal{D}$ are triangulated categories, then a \emph{triangulated duality} or a \emph{duality of triangulated categories} $\mathcal{C}\stackrel{\cong^o}{\longrightarrow}\mathcal{D}$   is  an equivalence of triangulated categories categories $\mathcal{C}^{op}\stackrel{\cong}{\longrightarrow}\mathcal{D}$.
As the following proposition shows, the definitions \ref{def.X-(co)reflexive objects}  are independent
of the considered distinguished families of idempotents.

\begin{Prop} \label{prop.RHom dualities}
Let $A$ and $B$ be dg algebras with enough idempotents, on which we
fix distinguished families of orthogonal idempotents
$(\epsilon_j)_{j\in J}$ and $(e_i)_{i\in I}$, respectively,
and let $X$ be a dg $B-A-$bimodule. The following assertions hold:

\begin{enumerate}
\item  $X$ is left homologically faithfully balanced if, and only
if, all objects of $\text{per}(B^{op})$ are homologically $X$-reflexive.
Then $\mathbb{R}\text{Hom}_{B^{op}}(?,X)$ and $\mathbb{R}\text{Hom}_{A}(?,X)$
define  quasi-inverse dualities of triangulated categories $\text{per}(B^{op})
\stackrel{\cong^o}{\longleftrightarrow}\text{thick}_{\mathcal{D}(A)}(e_iX\text{: }i\in I)$.
\item $X$ is right homologically faithfully balanced if, and only if,
all objects of $\text{per}(A)$ are homologically $X$-coreflexive.
Then $\mathbb{R}\text{Hom}_{A}(?,X)$ and $\mathbb{R}\text{Hom}_{B^{op}}(?,X)$
define  quasi-inverse dualities of triangulated categories $\text{per}(A)
\stackrel{\cong^o}{\longleftrightarrow}\text{thick}_{\mathcal{D}(B^{op})}(X\epsilon_j\text{: }i\in I)$.
\item If $A=B$ then the regular dg bimodule $X=A$ is homologically
faithfully balanced. In particular $\mathbb{R}\text{Hom}_{A}(?,A)$ and
$\mathbb{R}\text{Hom}_{A^{op}}(?,A)$ define quasi-inverse dualities
$\text{per}(A)\stackrel{\cong^o}{\longrightarrow}\text{per}(A^{op})$.
\end{enumerate}
\end{Prop}

\begin{proof}
The classes
$\text{Ref}(X)=\{U\in\mathcal{D}(B^{op})\text{: }\lambda_U\text{ is an isomorphism}\}$
and $\text{Coref}(X)=\{M\in\mathcal{D}(A)\text{: }\rho_M\text{ is an isomorphism}\}$
of $X$-reflexive and $X$-coreflexive objects are thick
subcategories of $\mathcal{D}(B^{op})$ and $\mathcal{D}(A)$, respectively (see the last three lines of the introduction). On the other hand, note that $\text{per}(B^{op})=\text{thick}_{\mathcal{D}(B^{op})}(Be_i\text{: }i\in I)$ and that $\text{per}(A)=\text{thick}_{\mathcal{D}(A)}(\epsilon_jA\text{: }j\in J)$ (see  Theorem \ref{thm.dg modules versus dg functors}, Remark \ref{rem.projective versus representable} and \cite[Section 5]{K1})

\bigskip

1) When  $X$ is left faithfully balanced, we then have that $\text{per}(B^{op})\subseteq\text{Ref}(X)$
.
Using now Lemma \ref{lem.equivalence reflexive-coreflexive},
we conclude that  $\mathbb{R}\text{Hom}_{B^{op}}(?,X)$ and $\mathbb{R}\text{Hom}_{A}(?,X)$
define quasi-inverse dualities between $\text{per}(B^{op})$
and the image of $\text{per}(B^{op})$ by $\mathbb{R}\text{Hom}_{B^{op}}(?,X)$.
This image is precisely
$\text{thick}_{\mathcal{D}(A)}(\mathbb{R}\text{Hom}_{B^{op}}(Be_i,X)\text{: }i\in I)$.
In order to prove assertion 1, it remains to check that there is an isomorphism
$\mathbb{R}\text{Hom}_{B^{op}}(Be_i,X)\cong e_iX$ in $\mathcal{D}(A)$,
for each $i\in I$. To see that, note that, due to the homotopically
projective condition of $Be_i$ (see Example \ref{ex.representables are homot.projective}),
if $\Pi :\mathcal{D}(B^{op})\longrightarrow\mathcal{H}(B^{op})$ is
the homotopically projective resolution functor, then $\Pi (Be_i)\cong Be_i$
in $\mathcal{H}(B^{op})$. We then have that
$\mathbb{R}\text{Hom}_{B^{op}}(Be_i,X)=\overline{HOM}_{B^{op}}(Be_i,X)$.
But the map $\Psi :\overline{HOM}_{B^{op}}(Be_i,X)\longrightarrow e_iX$,
given by $\Psi (f)=f(e_i)$ is an isomorphism of right dg $A$-modules.
Indeed, by Corollary \ref{cor.right A-structure on HOMB(U,X)},
we have that $\Psi (fa)=(fa)(e_i)=(-1)^{|a| |e_i|}f(e_i)a=f(e_i)a=\Psi (f)a$
since $|e_i|=0$, which immediately implies that $\Psi$ is an isomorphism in
$GR-A$ and $Gr-A$. On the other hand, if $\delta$ is the differential of
$\overline{HOM}_{B^{op}}(Be_i,X)$ and $d_{e_iX}=(d_X)_{| e_iX}$ is the
differential of $e_iX$, then we have
$(d_{e_iX}\circ\Psi )(f)=d_{e_iX}(f(e_i))=d_X(f(e_i))$ while we have
$$(\Psi\circ\delta )(f)=\Psi (\delta (f))=
\Psi (d_X\circ f-(-1)^{|f|}f\circ d_{Be_i})=
[d_X\circ f-(-1)^{|f|}f\circ d_{Be_i}](e_i)=d_X(f(e_i))$$
since $d_{Be_i}(e_i)=0$ because the differential of $B$ vanishes on $e_i$.
Therefore $\Psi$ is an isomorphism of right dg $A$-modules,
thus ending the proof of assertion 1.

\bigskip

2) Assertion 2 is proved as assertion 1 by exchanging the
roles of $A$ and $B^{op}$. Due to the fact that $\epsilon_jA$
is a homotopically projective right dg $A$-module, in a way
analogous to that of the previous paragraph, one checks that
the map
$\Phi :\mathbb{R}\text{Hom}_A(\epsilon_jA,X)=\overline{HOM}_A(\epsilon_jA,X)\longrightarrow X\epsilon_j$,
given by $\Phi (g)=g(\epsilon_j)$ is an isomorphism of left dg $B$-modules.

\bigskip

3)  For simplicity, put
$$F=\overline{HOM}_{A^{op}}(?,A): (A-Dg)^{op}\longrightarrow Dg-A$$ and
$$G=\overline{HOM}_A(?,A):(Dg-A)^{op}\longrightarrow A-Dg.$$
We then have $G\circ F^o:A-Dg\longrightarrow A-Dg$ and want to
get information about the unit
$\lambda :1_{\mathcal{D}(A^{op})}
\longrightarrow\mathbb{R}G\circ \mathbb{L}(F^o)=\mathbb{R}G\circ (\mathbb{R}F)^o$
(see Remark \ref{rem.not left derived of contravariant}).
For this, we consider the unit
$\tilde{\lambda}:1_{A-Dg}\longrightarrow G\circ F^o$ of the adjunction $(F^o,G)$.
By Proposition \ref{prop.dg transformation yields triang transformation},
we have an induced natural transformation of triangulated functors $\tilde{\lambda}:1_{\mathcal{D}(A^{op})}\longrightarrow\mathbb{L}(G\circ F^o)$
and, by Proposition \ref{prop.derived functor of composition}(2),
we get another natural transformation of triangulated functors
$\delta :\mathbb{L}(G\circ F^o)\longrightarrow\mathbb{R}G\circ (\mathbb{R}F)^o$.
It is not hard to see that $\lambda$ is the composition $1_{\mathcal{D}(A^{op})}
\stackrel{\tilde{\lambda}}{\longrightarrow}\mathbb{L}(G\circ F^o)
\stackrel{\delta}{\longrightarrow}\mathbb{R}G\circ (\mathbb{R}F)^o$.
If $j\in J$ is arbitrary, then $\Pi_A(A\epsilon_j)\cong A\epsilon_j$
in $\mathcal{H}(A^{op})$ since $A\epsilon_j$ is homotopically projective.
Moreover, by the proof of assertion 1 (with $A$ and $\epsilon_j$ instead
of $B$ and $e_i$),  we know that
$G^o(A\epsilon_j)=\overline{HOM}_{A^{op}}(A\epsilon_j,A)\cong\epsilon_jA$ in $Dg-A$.
It then follows from Proposition \ref{prop.derived functor of composition}(2)
that $\delta_{A\epsilon_j}$ is an isomorphism. Moreover, by
Proposition  \ref{prop.dg transformation yields triang transformation},
we know that if $$\tilde{\lambda}_{A\epsilon_j}:A\epsilon_j
\longrightarrow (G\circ F^o)(A\epsilon_j)=
\overline{HOM}_A(\overline{HOM}_{A^{op}}(A\epsilon_j,A),A)$$ is a
quasi-isomorphism (e.g. an isomorphism in $A-Dg$), then also
$$\tilde{\lambda}_{A\epsilon_j}:A\epsilon_j\longrightarrow\mathbb{L}(G\circ F^o)$$
is an isomorphism, and this will imply that $$\lambda_{A\epsilon_j}=
\delta_{A\epsilon_j}\circ\tilde{\lambda_{A\epsilon_j}}:
A\epsilon_j\longrightarrow [\mathbb{R}G\circ (\mathbb{R}F)^o](A\epsilon_j)$$
is an isomorphism in $\mathcal{D}(A^{op})$ and, hence, that $X={}_AA_A$ is
left homologically faithfully balanced.

We are led to give an explicit definition of
$\tilde{\lambda}_U:U\longrightarrow \overline{HOM}_A(\overline{HOM}_{A^{op}}(U,A),A)$,
for any left dg $A$-module. If we consider the natural isomorphism
$$\xymatrix{
\text{Hom}_{(Dg-A)^{op}}(\overline{HOM}_{A^{op}}(U,A),M)\ar@{=}[d]&&\text{Hom}_{A-Dg}(U,\overline{HOM}_A(M,A))\ar@{=}[d]\\
\text{HOM}_{A}(M,\overline{HOM}_{A^{op}}(U,A))\ar[rr]^{\xi=\xi_{U,M}}&&\text{HOM}_{A^{op}}(U,\overline{HOM}_A(M,A))}$$
%
(see the proof of Theorem \ref{thm.adjunction X-reflexive}),
the standard theory of adjunction gives that \newline
$\tilde{\lambda}_U=\xi_{U,\overline{HOM}_{A^{op}}(U,A)}(1_{\overline{HOM}_{A^{op}}(U,A)})$.
We then get that
$$\tilde{\lambda}_U(u)(\alpha )=\xi (1_{\overline{HOM}_{A^{op}}(U,A)})(u)(\alpha )=(-1)^{|\alpha | |u|}\alpha (u),$$
for all homogeneous elements $u\in U$ and $\alpha\in\overline{HOM}_{A^{op}}(U,A)$.
Abusing the notation and denoting by  $(?)^*$ both functors
$\overline{HOM}_A(?,A)$ and $\overline{HOM}_{A^{op}}(?,A)$,
we get a natural transformation
$\tilde{\lambda} :1_{A-Dg}\longrightarrow (?)^{**}$, where
$\tilde{\lambda}_U:U\longrightarrow U^{**}$ is given by the rule
$\tilde{\lambda}_U(\alpha )(u)=(-1)^{|\alpha| |u|}\alpha (u)$.

Consider the case $U=A\epsilon_j$ and let us consider now the
isomorphisms $\Psi: (A\epsilon_j)^*=\overline{HOM}_{A^{op}}(A\epsilon_j,A)\longrightarrow\epsilon_jA$
and $\Phi :(\epsilon_jA)^*=\overline{HOM}_A(\epsilon_jA,A)\longrightarrow A\epsilon_j$
given in the proofs of Assertions 1 and 2, when $X=A$. The composition
$$A\epsilon_j\stackrel{\Phi^{-1}}{\longrightarrow}(\epsilon_jA)^*\stackrel{\Psi^*}{\longrightarrow}(A\epsilon_j)^{**}$$
is then an isomorphism of left dg $A$-modules.
Note that $\Phi^{-1}(u)=f_u$, where $f_u(a)=ua$ for all $a\in\epsilon_jA$.
We then have that
$$(\Psi^*\circ\Phi^{-1})(u)=\Psi^*(f_u)=(-1)^{|\Psi| |f_u|}f_u\circ\Psi =f_u\circ\Psi,$$
using the action of the functor $(?)^*=\overline{HOM}_{A^{op}}(?,A)$ on homogeneous
morphisms (see the proof of Proposition \ref{prop.bifunctor})
and the fact that $|\Psi|=0$. Taking into account the comments after Lemma \ref{lem.left versus right dg module}, we  then have that
$$[(\Psi^*\circ\Phi^{-1})(u)](\alpha )=(f_u\circ\Psi )(\alpha )=
f_u(\alpha (\epsilon_j))=u\alpha (\epsilon_j)=
(-1)^{|u| |\alpha|}\alpha (u\epsilon_j)=(-1)^{|u| |\alpha|}\alpha (u),$$
for each $\alpha\in (A\epsilon_j)^*=\overline{HOM}_A(A\epsilon_j,A)$.
It follows that $\tilde{\lambda}_{A\epsilon_j}=\Psi^*\circ\Phi^{-1}$,
and hence that $\lambda_{A\epsilon_j}$ is an isomorphism, for all $j\in J$.

By a left-right symmetric argument, one checks that ${}_AA_A$ is also
right homologically faithfully balanced. The part of assertion 3 concerning duality follows from assertions 1 and 2 and from the first paragraph of this proof.
\end{proof}

\medskip
 We end with a result that has proved very useful in \cite{SZ2}:

\begin{Prop} \label{prop.transformation derived functors}
Let $\iota :A\longrightarrow B$ be a homomorphism of dg algebras with enough
idempotents such that $B=\iota (A)B\iota (A)$, and let us consider the
 dg functors: $$F:(Dg-A)^{op}
\stackrel{\overline{HOM}_A(?,A)}{\longrightarrow} A-Dg\stackrel{\iota^*}{\longrightarrow}B-Dg $$
 and
$$G:(Dg-A)^{op}\stackrel{\iota^{* o}}{\longrightarrow} (Dg-B)^{op}
\stackrel{\overline{HOM}_B(?,B)}{\longrightarrow}B-Dg,$$
where we denote by $\iota^*$ both extension of scalars functors $B\otimes_A?:A-Dg\longrightarrow B-Dg$ and $?\otimes B:Dg-A\longrightarrow Dg-B$. There is homological natural transformation of dg functors
$\eta:F\longrightarrow G$ whose triangulated version, when evaluated at compact objects,  gives
a natural isomorphism
$$\eta: [(B\otimes_A^\mathbb{L}?)\circ\mathbb{R}\text{Hom}_A(?,A)]_{| per(A)^{op}}\stackrel{\cong}{\longrightarrow}[\mathbb{R}\text{Hom}_B(?,B)\circ
(?\otimes_A^\mathbb{L}B)]_{\text{per} (A)^{op}}$$ of triangulated functors
$\text{per}(A)^{op}\longrightarrow\mathcal{D}(B)$.
\end{Prop}

\begin{proof}
All throughout the proof we fix a distinguished family of orthogonal
idempotents $(e_i)_{i\in }$ in $A$. Note that, after deleting the terms
which are zero, $(\iota (e_i))_{i\in I}$ is also a distinguished family
of orthogonal idempotents in $B$.
For each right dg $A$-module $M$, we define $$\eta_M:F(M)=B\otimes_A\overline{HOM}_A(M,A)\longrightarrow\overline{HOM}_B(M\otimes_AB,B)=G(M)$$
by the rule $\eta_M(b\otimes f)(m\otimes b')=b\iota (f(m))b'$,
for all homogeneous elements $b,b'\in B$, $f\in\overline{HOM}_A(M,A)$
and $m\in M$. We first check that $\eta:=\eta_M$ is well-defined.
Note that if $a\in A$ is a homogeneous element and $b,b',f,m$ are
homogeneous elements above, then we have
$$\eta (b\iota(a)\otimes f)(m\otimes b')=b\iota (a)\iota (f(m))b'$$
while $$\eta (b\otimes af)(m\otimes b')=b\iota((af)(m))b'=b\iota (af(m))b'=b\iota (a)\iota (f(m))b',$$
bearing in mind that the structure of right and left $A$-module
on $B$ is given by $a\cdot b\cdot a'=\iota (a)b\iota (a')$. Moreover,
if $b_1,b_2\in B$ are homogeneous elements, then we have that
$$\eta (b\otimes f)((m\otimes b_1)b_2)=\eta (b\otimes f)(m\otimes (b_1b_2))=b\iota (f(m))(b_1b_2)$$
and $$\eta (b\otimes f)(m\otimes b_1)b_2=(b\iota (f(m))b_1)b_2,$$
so that $\eta (b\otimes f)$ is a homogeneous morphism
$M\otimes_AB\longrightarrow B$ in $Dg-B$. On the other hand, if
$b=\sum_{i\in F}e_ib$, for a finite subset $F\subset I$, we clearly
have that $\text{Im}(\eta (b\otimes f))\subseteq\oplus_{i\in F}e_iB$,
thus showing that $\eta (b\otimes f)\in\overline{HOM}_B(M\otimes B,B)$
(see the initial paragraph of Section \ref{section.contravariant HOM}).
Therefore $\eta=\eta_M$ is a morphism
$F(M)=B\otimes_A\overline{HOM}_A(M,A)\longrightarrow\overline{HOM}_B(M\otimes_AB,B)$
in $Gr-K$. In order to check that it is actually a morphism in $B-Gr$, we need to check that $\eta (b_1(b\otimes f))=(-1)^{|b_1| | \eta|}b_1\eta (b\otimes f)=b_1\eta (b\otimes f)$. But this is clear since
$$\eta (b_1(b\otimes f))(m\otimes b')=\eta((b_1b)\otimes f)(m\otimes b')=(b_1b)\iota (f(m))b'$$ while
$$[b_1\eta (b\otimes f)](m\otimes b')=b_1\eta (b\otimes f)(m\otimes b')=b_1(b\iota (f(m))b'),$$
for each homogeneous element $b_1\in B$.

We now prove the naturality of $\eta$. Recall from the proof of
Proposition \ref{prop.tensor functor} that the action of
$B\otimes_A?:A-Dg\longrightarrow B-Dg$ on homogeneous morphisms is
given by the rule
$(B\otimes ?)(\alpha ):b\otimes x\rightsquigarrow (-1)^{|\alpha | |b|}b\otimes\alpha (x)$,
for all homogeneous elements $\alpha\in\text{HOM}_{B^{op}}(X,Y)$, $x\in X$ and $b\in B$.
Let $\varphi :M\longrightarrow N$ be a homogeneous morphism in $Dg-A$. Then
$$F(\varphi )=[(B\otimes_A?)\circ\overline{HOM}_A(?,A)](\varphi )=(B\otimes_A?)(\varphi^*).$$
This is a morphism $$F(N)=B\otimes_A\overline{HOM}_A(N,A)\longrightarrow
B\otimes_A\overline{HOM}_A(M,A)=F(M)$$ which takes
$b\otimes g\rightsquigarrow (-1)^{|\varphi^*| |b|}b\otimes\varphi^*(g)=(-1)^{|\varphi| |b|}b\otimes\varphi^*(g)$,
for all homogeneous elements $b\in B$ and $g\in\overline{HOM}_A(M,A)$.
By the definition of $\varphi^*$ (see the proof of Proposition \ref{prop.bifunctor}),
we then get that
$$F(\varphi )(b\otimes g)=
(-1)^{|\varphi| |b|}(-1)^{|\varphi| |g|}b\otimes (g\circ\varphi)=
(-1)^{|\varphi| (|g|+|b|)}b\otimes (g\circ\varphi)$$ On the other hand, we have that
$$G(\varphi )=[\overline{HOM}_B(?,B)\circ (?\otimes_AB)](\varphi )=(\varphi\otimes 1_B)^*.$$
This is a morphism
$$G(N)=\overline{HOM}_B(N\otimes_AB,B)\longrightarrow\overline{HOM}_B(M\otimes_AB,B)=G(M)$$
which takes
$u\rightsquigarrow (\varphi\otimes 1_B)^*(u)=(-1)^{|\varphi| |u|}u\circ (\varphi\otimes 1_B)$,
for each homogeneous element $u\in\overline{HOM}_B(N\otimes_AB,B)$. We then have
the following equalities, for all homogeneous elements $b\in B$ and $g\in\overline{HOM}_A(N,A)$:
$$[G(\varphi )\circ\eta_N](b\otimes g)=G(\varphi )(\eta_N(b\otimes g))
=(-1)^{| \varphi|(|b|+|g|)}\eta_N(b\otimes g)\circ (\varphi\otimes 1_B)$$
and
$$[\eta_M\circ F(\varphi )](b\otimes g)=(-1)^{| \varphi|(|b|+|g|)}\eta_M(b\otimes (g\circ\varphi )).$$
But, for all homogeneous elements $m\in M$ and $b'\in B$, we also have
$$[\eta_N(b\otimes g)\circ (\varphi\otimes 1_B)](m\otimes b')=
\eta_N(b\otimes g)(\varphi (m)\otimes b')=b\iota ((g\circ\varphi)(m)))b'=
\eta_M(b\otimes (g\circ\varphi))(m\otimes b'). $$
It follows that $G(\varphi )\circ\eta_N=\eta_M\circ F(\varphi )$, so
that $\eta$ is a natural transformation of dg functors.

We next prove that hat $\eta$ is homological, i.e., that
$d_{G(M)}\circ\eta_M-\eta_M\circ d_{F(M)}=0$. We  denote by
$d=d_{F(M))}:B\otimes_A\overline{HOM}_A(M,A)\longrightarrow B\otimes_A\overline{HOM}_A(M,A)$
and $\delta =d_{G(M)} :\overline{HOM}_B(M\otimes_AB,B)\longrightarrow\overline{HOM}_B(M\otimes_AB,B)$
the respective differentials. We need to check that
$\delta(\eta (b\otimes f))=\eta (d(b\otimes f))$, for all homogeneous
elements $b\in B$ and $f\in\overline{HOM}_A(M,A)$. For this, we shall
apply both members of this desired equality to a tensor $m\otimes b'$,
where $m\in M$ and $b'\in B$ are homogeneous elements. We then have:
\begin{eqnarray*}
[\delta \lefteqn{(\eta (b\otimes f))](m\otimes b')=}\\
&=&[d_B\circ\eta (b\otimes f))-(-1)^{|b|+|f|}\eta (b\otimes f))\circ d_{M\otimes_AB}](m\otimes b')\\
&=&d_B(b\iota (f(m))b')-(-1)^{|b|+|f|}\eta (b\otimes f))(d_M(m)\otimes b'+(-1)^{|m|}m\otimes d_B(b'))\\
&=&d_B(b)\iota (f(m))b'+(-1)^{|b|}bd_B(\iota (f(m))b')-(-1)^{|b|+|f|}[b\iota (f(d_M(m)))b'\\
&&+(-1)^{|m|}b\iota (f(m))d_B(b')]\\
&=&d_B(b)\iota (f(m))b'+(-1)^{|b|}b[d_B(\iota (f(m)))b'+(-1)^{|f|+|m|}\iota (f(m))d_B(b')]\\
&&-(-1)^{|b|+|f|}b\iota (f(d_M(m)))b'-(-1)^{|b|+|f|+|m|}b\iota (f(m))d_B(b')\\
&=&d_B(b)\iota (f(m))b'+(-1)^{|b|}bd_B(\iota (f(m)))b'-(-1)^{|b|+|f|}b\iota (f(d_M(m)))b'\;\;\;\;(*)
\end{eqnarray*}
while we also have
\begin{eqnarray*}
\eta \lefteqn{(d(b\otimes f))(m\otimes b')=}\\
&=&[\eta (d_B(b)\otimes f)+(-1)^{|b|}\eta (b\otimes d_{H}(f)))](m\otimes b')\\
&=&d_B(b)\iota (f(m))b'+(-1)^{|b|}b\iota (d_H(f)(m))b'\\
&=&d_B(b)\iota (f(m))b'+(-1)^{|b|}b\iota ((d_A\circ f)(m)-(-1)^{|f|}(f\circ d_M)(m))b' \hspace*{0.5cm} (**)
\end{eqnarray*}
The expression (*) and (**) are equal because $\iota :A\longrightarrow B$ is a morphism of
dg algebras with enough idempotents and, hence, one has that $\iota (d_A(a))=d_B(\iota (a))$,
for all $a\in A$.

For the final assertion, we start by pointing out that $\iota^*:Dg-A\longrightarrow Dg-B$
(resp. $\iota^*:A-Dg\longrightarrow B-Dg$) preserves homotopically projective dg modules.
Indeed if $P\in Dg-A$ is homotopically projective and $Y\in Dg-B$ is acyclic, then dg
adjunction gives an isomorphism $\text{HOM}_B(\iota^*(P),Y)\cong\text{HOM}_A(P,\iota_*(Y))$
in $Dg-K$. Since  $\iota_*$ preserves acyclic dg modules, we conclude that the last
dg $K$-module is acyclic and, hence, $\iota^*(P)$ is homotopically projective.

By Proposition \ref{prop.dg transformation yields triang transformation},
we have  an induced natural transformation of triangulated functors
$\eta:\mathbb{R}(\iota^*\circ\overline{HOM}_A(?,A))
\longrightarrow\mathbb{R}(\overline{HOM}_B(?,B)\circ\iota^{*o})$.
On the other hand, by Proposition \ref{prop.derived functor of composition},
we have natural transformations of triangulated functors $u:\mathbb{L}\iota^*\circ\mathbb{R}\text{Hom}_A(?,A)
\longrightarrow\mathbb{R}(\iota^*\circ\overline{HOM}_A(?,A))$ and
$v:\mathbb{R}(\overline{HOM}_B(?,B)\circ\iota^{*o})
\longrightarrow\mathbb{R}\text{Hom}_B(?,B)\circ (\mathbb{L}\iota^*)^o$.
The composition
$$\mathbb{L}\iota^*\circ\mathbb{R}\text{Hom}_A(?,A)
\stackrel{u}{\rightarrow}\mathbb{R}(\iota^*\circ\overline{HOM}_A(?,A))
\stackrel{\eta}{\rightarrow}\mathbb{R}(\overline{HOM}_B(?,B)\circ\iota^{*o})
\stackrel{v}{\rightarrow}\mathbb{R}\text{Hom}_B(?,B)\circ (\mathbb{L}\iota^*)^o $$
is then the desired natural transformation, which we want to prove that is an
isomorphism when evaluated at any $M\in\text{per}(A)$. Since $\text{per}(A)=\text{thick}_{\mathcal{D}(A)}(e_iA\text{: }i\in I)$
it is enough  to prove that
$(v\circ\eta\circ u)_{e_iA}=v_{e_iA}\circ\eta_{e_iA}\circ u_{e_iA}$ is an
isomorphism, for all $i\in I$.

Proposition \ref{prop.derived functor of composition}(4)  together with the
previous to the last paragraph tell us that $v$ is a natural isomorphism. On
the other hand,   $\Pi_A(e_iA)\cong e_iA$ in $\mathcal{H}(A)$ since $e_iA$ is
homotopically projective, for all $i\in I$. But $\overline{HOM}_A(e_iA,A)\cong Ae_i$
and, by the previous to the last paragraph again, also $[\iota^*\circ\overline{HOM}_A(?,A)](e_iA)$
is homotopically projective. Proposition \ref{prop.derived functor of composition}(3) then
gives  that $u_{e_iA}$ is an isomoprhism and
Proposition \ref{prop.dg transformation yields triang transformation}
implies that, in order to check that $\eta_{e_iA}$ is an isomorphism in $\mathcal{D}(B^{op})$
and hence end the proof, it is enough to prove that
$$\eta_{e_iA}:F(e_iA)=B\otimes_A\overline{HOM}_A(e_iA,A)
\longrightarrow\overline{HOM}_B(e_iA\otimes B,B)=G(e_iA)$$ is an isomorphism of left dg $B$-modules.

We proceed to prove  this fact. Recall from the proof of
Proposition \ref{prop.RHom dualities} that we have isomorphisms
$\Phi_A:\overline{HOM}_A(e_iA,A)\stackrel{\cong}{\longrightarrow}Ae_i$ and
$\Phi_B:\overline{HOM}_B(e_iB,B)\stackrel{\cong}{\longrightarrow}Be_i$, in
$A-Dg$ and $B-Dg$, respectively, mapping $f\rightsquigarrow f(e_i)$ in both
cases. Note that $\Phi_B^{-1}:Be_i\longrightarrow\overline{HOM}_B(e_iB,B)$
is given by the rule $\Phi_B^{-1}(b)(e_i)=b$ or, equivalently, by the rule
$\Phi_B^{-1}(b)(b')=bb'$, for all homogeneous elements $b\in Be_i$ and
$b'\in e_iB$. Note also that since $\Phi_A$ has degree zero, the induced
isomorphism $(B\otimes ?)(\Phi_A):B\otimes_A\overline{HOM}_A(e_i,A)
\stackrel{\cong}{\longrightarrow}B\otimes_AAe_i$ takes $b\otimes f
\rightsquigarrow b\otimes f(e_i)$, so that $(B\otimes_A?)(\Phi_A)=1_B\otimes\Phi_A$.
We now claim that $\eta_{e_iA}:B\otimes_A\overline{HOM}_A(e_iA,A)
\longrightarrow\overline{HOM}_B(e_iA\otimes_AB,B)$ can be decomposed
as the composition of morphisms in $B-Dg$

\begin{center}
$B\otimes_A\overline{HOM}_A(e_iA,A)
\stackrel{1_B\otimes\Phi_A}{\longrightarrow}B\otimes_AAe_i
\stackrel{\mu '}{\longrightarrow}Be_i
\stackrel{\Phi_B^{-1}}{\longrightarrow}\overline{HOM}_B(e_iB,B)
\stackrel{\mu^*}{\longrightarrow}\overline{HOM}_B(e_iA\otimes_AB,B)$,
\end{center}
where $\mu':B\otimes_AAe_i\longrightarrow Be_i$ and
$\mu :e_iA\otimes B\longrightarrow B$ are the multiplication maps,
$b\otimes a\rightsquigarrow b\iota (a)$ and $a\otimes b\rightsquigarrow \iota (a)b$.  Indeed we have:
\begin{eqnarray*}
(\mu^*\circ\Phi_B^{-1}\circ\mu'\circ (1_B\otimes\Phi ))(b\otimes f)
&=&(\mu^*\circ\Phi_B^{-1}\circ\mu') (b\otimes f(e_i))\\
&=&(\mu^*\circ\Phi_B^{-1})(b\iota (f(e_i)))\\
&=&\Phi_B^{-1}(b\iota (f(e_i)))\circ\mu
\end{eqnarray*}
since $\mu$ has zero degree. When we take homogeneous elements
$x\in e_iA$ and $b'\in B$ and make act the last morphism on $x\otimes b'$, we get
\begin{eqnarray*}
(\mu^*\circ\Phi_B^{-1}\circ\mu'\circ (1_B\otimes\Phi ))(b\otimes f)(x\otimes b')
&=&\Phi_B^{-1}(b\iota (f(e_i))) (\iota (x)b')\\
&=&b\iota (f(e_i))\iota (x)b'=b\iota (f(e_ix))b'\\
&=&b\iota (f(x))b'=\eta_{e_iA}(b\otimes f)(x\otimes b'),
\end{eqnarray*}
using the fact that $\iota :A\longrightarrow B$ is an algebra homomorphism
and $f:e_iA\longrightarrow A$ is a morphism of right $A$-modules.

The proof is hence reduced to check that $\mu':B\otimes_AAe_i\longrightarrow Be_i$
and $\mu :e_iA\otimes B\longrightarrow e_iB$ are isomorphisms in $B-Dg$ and $Dg-B$,
respectively. But this is clear. Their inverses map $b\rightsquigarrow b\otimes e_i$
and  $b\rightsquigarrow e_i\otimes b$, respectively.
\end{proof}

{\footnotesize

\noindent
Manuel Saor\'\i{}n
\newline Departemento de Matem\'aticas,
\newline Universidad de Murcia, Aptdo. 4021
\newline 30100 Espinardo, Murcia,
\newline Spain
\newline{\tt msaorinc@um.es}

}


\begin{thebibliography}{88}

\bibitem{AF}
Frank W. Anderson and Kent R. Fuller, {\sc Rings and categories of modules}, 2nd edition, Grad. Texts Maths {\bf 13}, Springer-Verlag (1992). 

\bibitem{B}
Thomas B\"uhler, {\em Exact categories}, Expo. Math. {\bf 28}(1)  (2010), 1-69



\bibitem{D1}
Yuriy A. Drozd, {\em Tame an wild matrix problems}, in Representations and quadratic forms, Kiev 1979; english translation: Amer. Math. Soc.  Transl. {\bf 128}(2) (1986), 31-55.

\bibitem{D2}
Yuriy A. Drozd, {\em Tame and wild matrix problems}, in Representation Theory II, Proc. Confer. Ottawa 1979; V. Dlab and P. Gabriel (edts.). Springer Lect. Notes Math. {\bf 832} (1980), 242-258. 

\bibitem{G}
Pierre Gabriel, {\em Des cat\'egories ab\'eliennes}, Bull. Soc. Math. France {\bf 90} (1962), 323-448.

\bibitem{GOR}
Natalia S. Golovaschuk, Serge Ovsienko and Andrei V. Rojter, {\emph On the schurian DGC}, Matrix problems, IM AN USSR, Kiev (1977), 162-165.


\bibitem{H}
Dieter Happel, {\sc Triangulated Categories in the Representation Theory
of Finite Dimensional Algebras}, London Math. Soc. Lect. Note Ser. {\bf 119}.
Cambridge University Press 1988.

\bibitem{HS}
Peter Hilton, Urs Stammbach, {\sc A Course in Homological Algebra}, 2nd edition.
Grad. Texts Math. {\bf 4}, Springer-Verlag (1971).

\bibitem{KS}
Masaki Kashiwara and Pierre Shapira, {\sc Categories and sheaves}, Grundl. Math. Wiss {\bf 332}, Springer-Verlag (2006)

\bibitem{K1} Bernhard Keller, {\em Deriving DG categories}, Ann.
Sci. \'Ecole Norm. Sup {\bf 27} (1994), 63-102.

\bibitem{K3} Bernhard Keller, {\em Derived categories and their uses}, Handbook of Algebra, vol. I, North-Hollad (1996), 671-701.

\bibitem{K2} Bernhard Keller, {\em On differential graded categories}, In:
International Congress of Mathematics, vol. II. Eur. Math. Soc. Zurich (2006), 151-190.

\bibitem{Kl-R}
Mark M. Kleiner and Andrei V. Rojter, {\em Representations of differential graded categories}, in Proceed. 1st International Conference on Representations of Algebras, Ottawa 1974, Springer Lect. Notes Math. {\bf 488} (1975), 316-339.

\bibitem{M}
Barry Mitchell, {\em Rings with several objects}, Adv. Math. {\bf 8} (1972), 1-161

\bibitem{NVO}
Constantin Nastasescu and Freddy Van Oystaeyen, {\sc Graded Ring Theory}, North-Holland (1982).

\bibitem{Neeman}
Amnon Neeman, {\sc Triangulated Categories}, Princeton University Press 2001.



\bibitem{NS-Japan}
Pedro Nicol\'as and Manuel Saor\'in, {\em Classical derived functors as fully
faithful embeddings}. Proc. 46th Japan Symp. Ring Theory and Repres. Theory
(edited by I. Kikumasa). Yamaguchi University (2014), 137-187.

\bibitem{NS-GTT}
 Pedro Nicol\'as and Manuel Saor\'in, {\em Generalized tilting theory}, Preprint available at https://arxiv.org/abs/1208.2803

\bibitem{O}
Serge Ovsienko, {\em Bimodule and matrix problems}, in 'Computational Methods for Representations of Groups and Algebras', Proceed. Euroconfer. Essen 1977, P. Dr\"axler, C.M. Ringel and G.O. Michler (edts.), Birkh\"auser Progress in Maths. {\bf 173} (1999), 325-357.

\bibitem{SZ2}
Manuel Saor\'\i n and Alexander Zimmermann, {\em Symmetry of the definition of degeneration in triangulated  categories}, Preprint


\bibitem{Ta}
Goncalo Tabuada, {\em Une structure de cat\'egorie de mod\`eles de Quillen sur la cat\'egorie des dg-cat\'egories}, Compt. Rend. Acad. Sci. Paris, s\'er {\bf I 340} (2005), 15-19.

\bibitem{T}
Bertrand To\"en, {\em The homotopy theory of dg-categories and derived Morita theory}, Invent. Math. {\bf 167} (2007), 615-667.

\bibitem{Verdier}
Jean-Louis Verdier,  {\em Des cat\'egories d\'eriv\'ees des cat\'egories abeliennes}.
Ast\'erisque {\bf 239}, Soc. Math. France (1996).



\bibitem{Wis}
Robert Wisbauer, {\sc Foundations of Module and Ring Theory}, Gordon and Breach Science Publishers (1991).

\bibitem{Yek}
Amnon Yekutieli, {\sc Derived categories - A textbook}. Available at https://arxiv.org/abs/1610.09640 and https://www.math.bgu.ac.il/~amyekut/teaching/2016-17/der-cats-IV/public53.pdf

\bibitem{reptheobuch}
Alexander Zimmermann, {\sc Representation Theory; A homological algebra point of view},
Springer Verlag (2014).



\end{thebibliography}
\end{document}